%% file: main.tex
\documentclass[
11pt, 
oneside, 
english, 
singlespacing, 
headsepline, 
]{MastersDoctoralThesis} 

\usepackage[utf8]{inputenc} 
\usepackage[T2A,T1]{fontenc}
\usepackage{xcolor}

\usepackage{caption}
\usepackage{subcaption}

\usepackage[backend=biber,  
sorting=none, 
style=numeric-comp,
natbib=true,
maxbibnames=10]{biblatex}

\usepackage[autostyle=true]{csquotes} 

\usepackage{pdfpages}

\usepackage{todonotes}
\usepackage{listings}

\definecolor{codegreen}{rgb}{0,0.6,0}
\definecolor{codegray}{rgb}{0.5,0.5,0.5}
\definecolor{codepurple}{rgb}{0.58,0,0.82}
\definecolor{backcolour}{rgb}{0.921, 0.929, 0.937} 

\lstdefinestyle{mystyle}{
backgroundcolor=\color{backcolour},   
commentstyle=\color{codegreen},
keywordstyle=\color{magenta},
numberstyle=\tiny\color{codegray},
stringstyle=\color{codepurple},
basicstyle=\ttfamily\footnotesize,
breakatwhitespace=false,         
breaklines=true,                 
captionpos=b,                    
keepspaces=true,                 
numbers=left,                    
numbersep=5pt,                  
showspaces=false,                
showstringspaces=false,
showtabs=false,                  
tabsize=2
}

\usepackage{makecell} 

\usepackage{mypackage}

\thesistitle{Spline-based Reconstruction of Periodic Signals with Sparse Innovations} 
\supervisorI{Dr. Matthieu \textsc{Simeoni}} 
\supervisorII{Dr. Julien \textsc{Fageot}}
\examiner{} 
\degree{Engineer} 
\author{Adrian \textsc{Jarret}} 
\addresses{} 

\subject{Biological Sciences} 
\keywords{} 
\university{\href{https://www.universite-paris-saclay.fr/en}{Université Paris-Saclay}} 
\universitygerman{\href{https://www.zhaw.ch/de/hochschule/}{Z{\"u}rcher Hochschule f{\"u}r Angewandte Wissenschaften}}
\department{\href{https://www.epfl.ch/en/}{Ecole Polytechnique Fédérale de Lausanne}} 
\group{\href{https://www.epfl.ch/labs/lcav/}{Laboratoire des Communications AudioVisuelles}} 
\faculty{\href{https://www.epfl.ch/schools/ic/}{School of Computer and Communication Sciences
}} 

\AtBeginDocument{
\hypersetup{pdftitle=\ttitle} 
\hypersetup{pdfauthor=\authorname} 
\hypersetup{pdfkeywords=\keywordnames} 
}

\begin{document}

\frontmatter 

\pagestyle{plain} 



\begin{titlepage}
\begin{center}

\begin{figure}
\makebox[\textwidth][c]{
\begin{subfigure}{0.6\textwidth}
    \centering
    \includegraphics[width=0.8\linewidth,center]{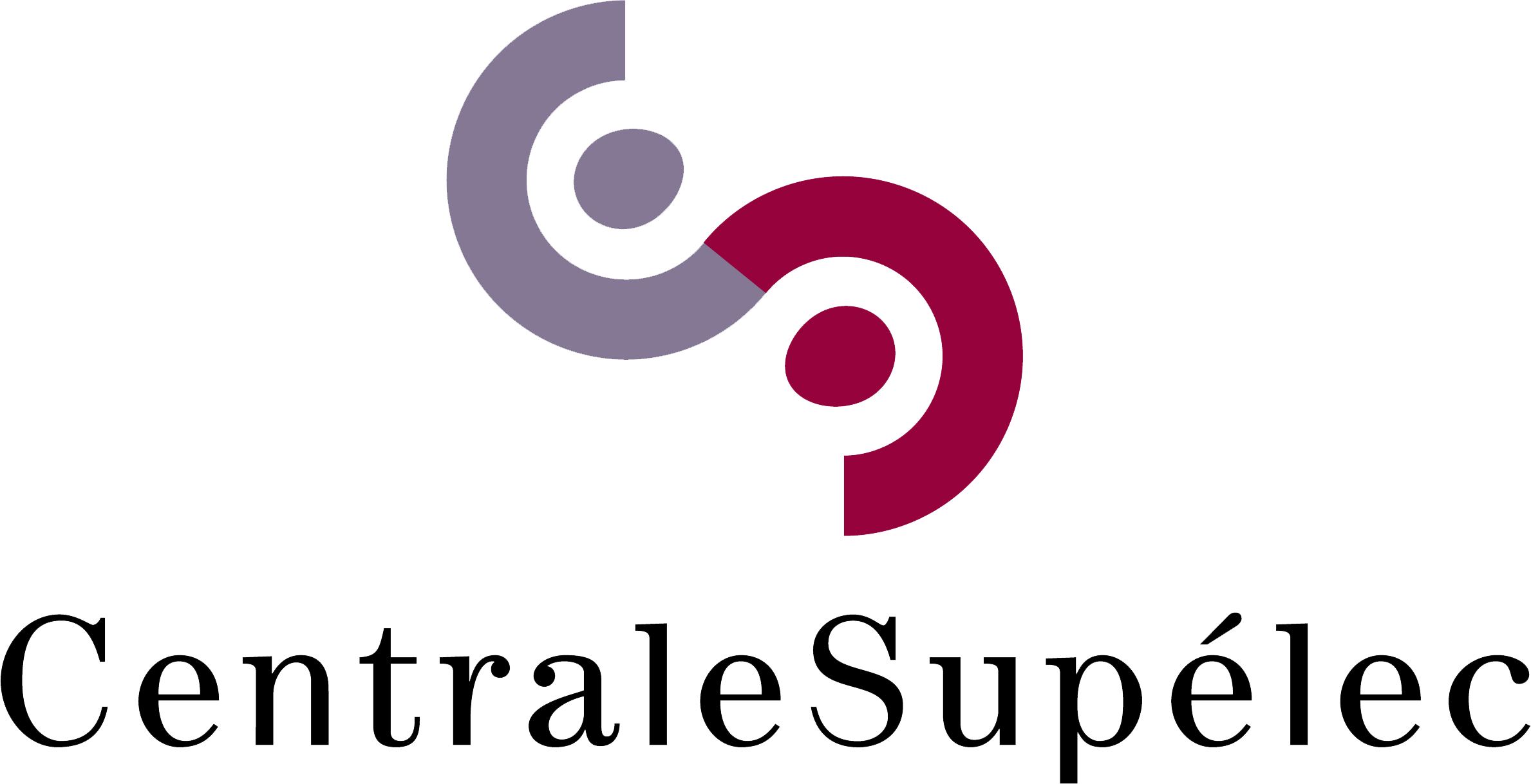}
\end{subfigure}
\begin{subfigure}{0.6\textwidth}
    \centering
    \includegraphics[width=0.8\linewidth, center]{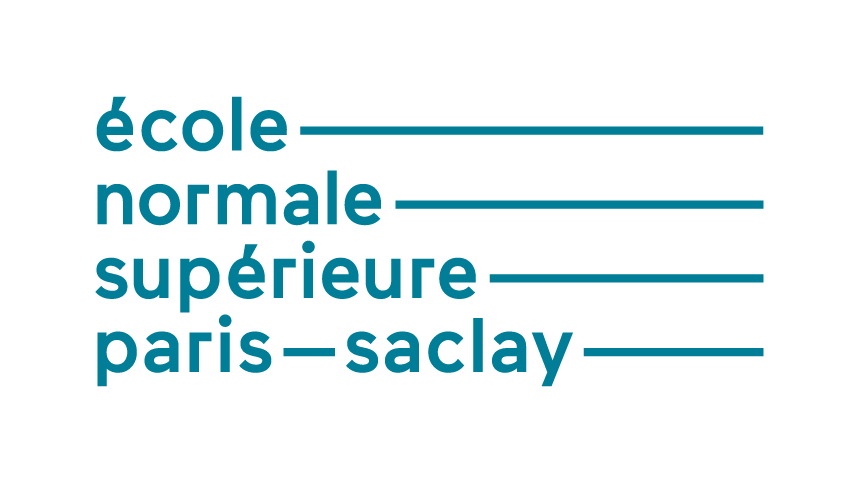}
\end{subfigure}
}
\end{figure}

\vspace*{.06\textheight}
{\scshape\LARGE \univname\par}\vspace{1.5cm} 
\textsc{\Large Master Internship Report}\\[0.5cm] 

\HRule \\[0.4cm] 
{\huge \bfseries \ttitle\par}\vspace{0.4cm} 
\HRule \\[1.5cm] 
 
\begin{minipage}[t]{0.4\textwidth}
\begin{flushleft} \large
\emph{Author:}\\
\authorname 
\end{flushleft}
\end{minipage}
\begin{minipage}[t]{0.4\textwidth}
\begin{flushright} \large
\emph{Supervisors:} \\
\href{https://people.epfl.ch/matthieu.simeoni}{\supnameI}\\ 
\href{http://bigwww.epfl.ch/fageot/index.html}{\supnameII}
\end{flushright}
\end{minipage}\\[2cm]
 
\vfill

\large \textit{Internship taking place in the }\\[0.3cm] 
\groupname \\[0.2cm]
\textit{at} \\[0.2cm]
\deptname  \\[.8cm] 

\begin{figure}[H]
    \centering
    \includegraphics[width=0.3\linewidth,center]{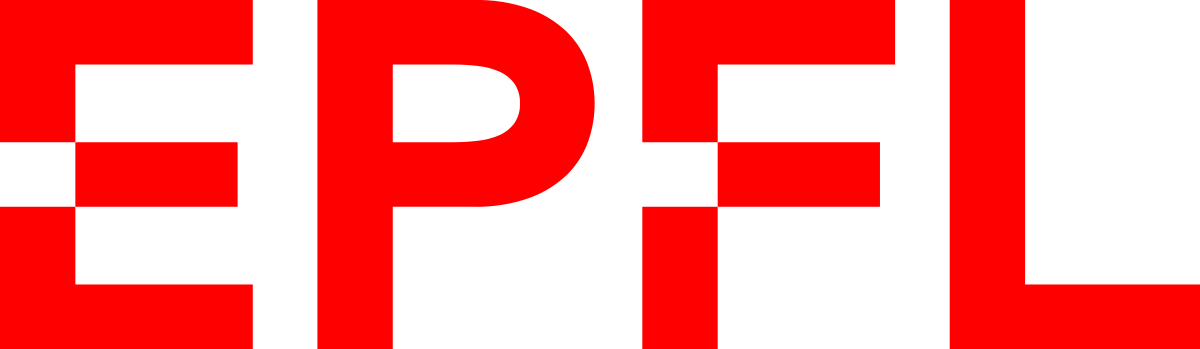}
\end{figure}

\vfill

{\large November 27, 2020}\\[4cm] 

\vfill
\end{center}
\end{titlepage}


\begin{abstract}

This report summarizes the research work achieved during this master internship. Optimization-based problems have become of great interest for signal approximation purposes, as they achieved good accuracy results while being extremely flexible and versatile. In this work, we put our focus on the context of periodic signals sampled with spatial measurements. The optimization problems are penalized thanks to the \emph{total variation} norm, using a specific class of (pseudo-)differential operators to use well-chosen reconstruction functions. We introduce three algorithms and their adaptation to this specific use case. The first one is a discrete grid-based $\ell_1$ method, the second is called CPGD, which relies on the estimation of discrete innovations within the FRI framework, and the last one is the \emph{Frank-Wolfe} algorithm. We put the emphasis on that later algorithm as we underline its greedy behavior. We consider a refined version of this algorithm, that leads to very encouraging outputs.

We pursue the algorithmic treatment with a theoretical analysis of the optimization problem, based on convexity and duality theory. We draw away two conditions that, if verified by the differential operator, would inform on the shape of the solutions and even ensure uniqueness, in some cases. Moreover, we demonstrate the general expression of periodic one-dimensional exponential splines, as a first step to verify the latter conditions for exponential operators.
\end{abstract}






\chapter*{Introduction}

    Signal reconstruction relying on optimization problems has been of great interest in the field of signal processing and compressed sensing lately. Beyond obtaining accurate experimental results, the underlying optimization part features both interesting theoretical aspects and a wide range of varying parameters. Inspired from domains as broad as radioastronomy, environmental data monitoring or acoustic sensing, this project focuses on the specific context of periodic signal reconstruction with non-uniform spatial sampling.
    
    Historically, such problems have been treated by way of discretization strategies, leading to powerful frameworks and approximations that are still used nowadays. One of the main singularities of this work is that it tries to tackle the signal reconstruction directly in a continuous setting. It actually corresponds to the most natural way of formulating the signal approximation problems of interest. It induces more flexibility in the reconstruction procedure, however it also requires a more conscious and careful statement of the problem and of the involved solutions.
    
    In such a context, penalization is necessary for \enquote{well-posedness} of the problem, but more than that it is also crucial to orient the solution.In such a context, the \enquote{well-posedness} of the problem is usually obtained through penalization, as it has proven to perform really well over the numerous use cases. Past the practical advantages, penalization is also crucial here to orient the solutions. The design of the regularization component is very subtle and completely transforms the shape of the reconstructed signals. In this project, we exclusively study \emph{total variation} penalization, which is the continuous counterpart of the well-known $\ell_1$ norm (also known as \emph{LASSO regression} or \emph{basis pursuit} depending on the application domain). It happens that this framework is efficient to promote sparse solutions, and thus is extensively used to this end. Beyond the fact that sparsity is usually a property of great interest in approximation (according to Occam's razor principle), in our context the use of the total variation norm also allows to decorrelate the measurement procedure and the general form of the solutions. This property is extremely helpful with our irregular spatial measurements, as the recovered solutions become \emph{independent} (in some sense) from the sampling strategy. Additionally, we make use of (pseudo-)differential operators in the penalization framework. These operators basically model the shape of the output by controlling the smoothness and the spatial extension of any solution.
    
    The periodicity and the spatial sampling are directly inspired from actual data sensing, as any directional signal can be seen as a spatially periodic signal around the measurement device. The theoretic analysis of such periodic functions and operators is thoroughly achieved thanks to the Fourier analysis, and happens to be simpler than the non-periodic case. Indeed, it nearly reduces to the analysis of compactly-supported functions and the behaviors at infinity are no longer a trouble. The periodic context has mostly gained attention recently, and now belongs to the modern fields of study.
    
    In a few words, this work focuses on the reconstruction task of an input signal through the prism of optimization-based methods. The problems are expressed in a continuous periodic setting and the data is obtained through arbitrary spatial sampling. The core of this framework, which enables its adaptability and precision, lies in the penalization strategy, thanks to the total variation norm coupled to specific well-chosen operators. The idea of this document is to give the broadest understanding of the reconstruction framework.
    
    The contributions of this work are of two kinds:
    \begin{itemize}
        \item We provide the implementation of the three studied strategies : \emph{equispaced grid}, \emph{CPGD} and \emph{Frank-Wolfe}, within a unified algorithmic framework. They cover the context of periodic one-dimensional reconstruction, with a wide range of penalization operators. The framework allows to easily generate synthetic samples, run and compare the methods with any set of custom reconstruction parameters. As far as the Frank-Wolfe strategy is concerned, this project has been through deep analysis of its greedy behavior and studies a refined version, namely the complete reweighting, about which very few literature is available.
        \item For the theoretical aspects, this project has taken great care of carefully explaining the dual analysis of the generalized continuous basis pursuit problem. We have obtained a generic approach to prove that, depending on the involved operator, the solutions are splines and might be unique. Our work paves the way for a systematic study of the (non-)uniqueness of TV-based optimization for spatial sampling and pseudo-differential regularization operators. We start to investigate the case of exponential operators, and successfully prove that the corresponding dual certificates always admit a finite number of saturation points. This result relies the general expression of the periodic one-dimensional exponential splines that we demonstrate.
    \end{itemize}
    
    In Chapter \ref{chap:preli}, we present the mathematical framework to properly establish the optimization problems. In Chapter \ref{chap:algo}, we go over three different algorithmic strategies to actually solve the aforementioned optimization problems. We precisely describe their mechanics and the adaptations to our context. Chapter \ref{chap:experiences} illustrates these algorithms with practical experiments and use cases. Finally, Chapter \ref{chap:dual_analysis} is part of an ongoing theoretical effort to deduce properties on the solution set of the studied problems. We strongly believe that it is possible to characterize the shape of the solutions as well as to prove uniqueness in some cases, as experiments let it appear.
    
\section*{Related Works}
    
    \paragraph{Discrete $\ell_1$ methods}
    Inverse problems have been expressed and treated for a long time in a discrete setting, involving either vector variable as well as infinite dimensional variables. The use of $\ell_1$-norm penalization was initially justified for sparsity promoting, in comparison with the frequently used Tikhonov regularization ($\ell_2$-norm). The statistics community refers to the obtained problem as \emph{LASSO}, for \emph{least absolute shrinkage and selection operator} \cite{tibshirani96regression}, meanwhile it is commonly referred to under the name of \emph{basis pursuit} (BP) within the signal processing world \cite{chen98basis}. These problems and their resolution have given birth to the field of compressed sensing. More recently, theoretical studies have been proposed \cite{tibshirani2013uniqueness, ali2019generalized} to address the question of uniqueness, in the case of regular $\ell_1$ penalization as well as when a penalization operator is used (generalized LASSO). Fast and accurate algorithmic strategies, relying on proximal methods, have been proposed in the last decade to numerically approximate solutions (primal dual splitting algorithm, FISTA, ...) \cite{condat:hal-00609728, amir2009fista}. Although recent, they have turned the $\ell_1$ penalized problems into an extensively used tool.

    \paragraph{Optimization over measure spaces}
    When aiming at approximating real life signals, one might express optimization problems in the continuous setting, with data living in a infinite dimensional space (spaces of functions, for instance). It allows more flexibility in the reconstruction, compared to the finite dimensional case, as the reconstructed data are no longer constrained to live within the available dimensions of the problem. Using continuous data might also lead to more natural expression of the optimization problems, with straightforward objective cost functionals. It often amounts to recover data whose shape is a stream of Dirac $w = \sum \alpha_k \delta_{\mathbf{x}_k}$. It enters the framework of so-called functions with \emph{finite rate of innovations} (FRI), which is the continuous counterpart of the discrete notion of sparsity. Some reconstruction algorithms propose to leverage this sparsity property to directly estimate the finite number of reconstruction positions with help of Prony's methods \cite{blu2008fri, simeoni2020cpgd}.
    In functional spaces, the generalization of $\ell_1$-norm is the \emph{total variation} norm, which is commonly used to penalize such optimization problems, as in \cite{COCV_2013__19_1_190_0}. It is used for sparse spikes deconvolution or super-resolution problems, as in \cite{duval2014exact, candes2012mathematical}. There are sometimes referred to as \emph{continuous basis pursuit} \cite{ekanadham2011cbp, duval2017cbp} or \emph{generalized minimal extrapolation} \cite{decastro2012exact}. These problems can numerically be solved or approximated by means of versions of the \emph{matching pursuit} algorithm (that generalizes in a continuous setting) \cite{knudson2014inferring} or with the \emph{conditional gradient method} \cite{COCV_2013__19_1_190_0, denoyelle2018sliding}.
    
    \paragraph{Spline-based approximation}
    Piecewise constant functions are the most well-known splines, but it is actually possible to define splines more generally. A spline is essentially a piecewise-defined function, whose parts are shifted elements from the null space of a (pseudo)-differential operator. The piecewise constant spline is obtained with the derivative operator, and using other operators (with satisfying properties) would lead to other type of spline. The operator acts as a generalization of the derivative and allows to hand-tailor the obtained spline.
    It has been proven that introducing such a generalized derivative operator into the total variation penalization term of an optimization problem conduces to solutions expressed in terms of splines \cite{Unser_2017}. More interestingly, the reconstruction splines are the splines issued from the operator used in the regularization term, independently of the sampling procedure. This results certifies that it is very consistent to try to approximate any input signal with a sum of splines. Fundamentally, this new penalization paradigm substitutes splines to the previously used Dirac impulses in the reconstruction strategy. It turns out to be really convenient when the data does not fit well with Dirac reconstruction, in the case of spatial sampling for instance. An interesting case  over the real line is studied in \cite{debarre2020sparsest} and it has been generalized to the periodic setting (over the multidimensional torus) \cite{fageot2020tvbased} or over the sphere \cite{Simeoni:275337}.


\mainmatter 

\pagestyle{thesis} 


\include{Chapters/chapter1}
\include{Chapters/chapter2} 
\include{Chapters/chapter3}
\include{Chapters/chapter4}

\chapter*{Conclusion and acknowledgments}
The explicit objectives of this internship were to understand, develop and adapt the three algorithms, namely the equispaced knots, Cadzow Plug-and-Play Gradient Descent and Frank-Wolfe, to the context of reconstruction of periodic signal with sparse innovations and non-uniform spatial sampling. These goals have successfully been achieved. We even came out with an improved version of the Frank-Wolfe algorithm that outperforms the regular one in all the tested scenarios. On the way, we have developed our understanding of the mechanics at play. It has led us to really interesting theoretical questions, involving duality and convex optimization theory.

On that topic, we have designed a demonstration framework that would allow to both characterize the shape of the elements of the solution set and to ensure uniqueness of the solution. We have succeeded to start applying this framework to the case of exponential operators. In particular, we came up with the general expression of the Green's function associated to the one-dimensional periodic exponential operator. However, this work remains uncompleted, and various ideas are on the verge to be considered in that direction. In further works, it would be really meaningful to focus on other operators such as Sobolev or Wendland. We strongly believe that many results would logically come out of the current achieved work on the topic. \\

As far as I am concerned, this internship has been extremely enriching and educative. Although essentially theoretic, it has brought me to discover compress sensing, which is a field I had very limited experience with. Beyond mathematical material and programming skills, these six months of research work have also taught me how to pursue an investigation work in an academic environment, by planning my working objectives and making sure I keep the motivation and interests high enough to make progress. This is a significant personal improvement compared to my previous research experience two years ago. The maturity I have acquired from that former internship has definitely helped me to make this current experience more productive and rewarding.

Doing a PhD is the logical follow-up for this work, and this is what I am willing to enroll myself in. The academic research environment fits well with my personality and my way of learning and discovering. I am also interested in continuing this project. A lot of progress can be made in very different directions, and I would really enjoy turning it into to some concrete contributions. I have already applied to a doctoral school at EPFL. The choices of an advisor and a topic are really crucial and I would like to be conscientious in making them. That is why I will also make use of the following weeks to browse thesis proposals and to carefully explore other opportunities. Nonetheless, I have really appreciated working at EPFL, and specifically in LCAV, with my supervisors and their colleagues from other laboratories and I am sure it would provide a great environment to produce a high quality successful PhD thesis.
\\

During this internship, I have been closely supervised by both Matthieu and Julien. I would like to deeply thank them for the help, patience and attention they provided about my work, and I am sincerely looking forward to keeping this fruitful working environment with them.






\printbibliography


\end{document}

%% file: Chapters/chapter1.tex
\chapter{Problem statement} 

\label{chap:preli}

An extensive presentation of the mathematical context for periodic functional recovery with total variation penalization has been achieved by Fageot and Simeoni in their paper on periodic splines \cite{fageot2020tvbased}. In this section, we summarize the key points of the mathematical framework, and we refer to their work for a more exhaustive treatment.


\section{Search Space and Splines}

\subsection{A few words of topology}
We denote by $\mathcal{S(\mathbb{T})}$ the space of infinitely differentiable $T$-periodic functions over $\mathbb{R}$. $\mathcal{S'(\mathbb{T})}$ is the topological dual of $\mathcal{S(\mathbb{T})}$, which is the space of periodic generalized functions. The duality product between $\mathcal{S(\mathbb{T})}$ and $\mathcal{S'(\mathbb{T})}$ is usually noted as:
\begin{equation}
    \langle \cdot , \cdot \rangle : \quad \mathcal{S'(\mathbb{T})} \times \mathcal{S(\mathbb{T})} \to \mathbb{R}, \qquad (w, f) \mapsto \int_\mathbb{T}{f \mathrm{d}w}.
\end{equation}

The space of continuous periodic functions over $\mathbb{T}$ is denoted as $\mathcal{C}(\mathbb{T})$. We define the \emph{supremum} of \textit{infinite} norm over $\mathcal{C}(\mathbb{T})$ as $\norm{\varphi}_\infty = \sup_{t\in\mathbb{T}}{|\varphi(t)|}$. Endowed with this norm, $\left(\mathcal{C}(\mathbb{T}), \norm{\cdot}_\infty\right)$ forms a Banach space whose dual space is isomorphic to the space of periodic finite Radon measures $\mathcal{M}(\mathbb{T})$. Thus $\mathcal{M}(\mathbb{T})$ is a Banach space when it is endowed with the dual norm (also coined as \emph{total variation} norm):
\begin{equation}
    \forall w \in \mathcal{M}(\mathbb{T}), \quad \norm{w}_{TV} = \sup_{\varphi \in \mathcal{C}(\mathbb{T}), \norm{\varphi}_\infty \leq 1} {\langle w, \varphi \rangle}
\end{equation}
with ${\langle w, \varphi \rangle = \int_\mathbb{T}{\varphi \ \mathrm{d}w}}$.
As $\mathcal{S}(\mathbb{T})$ is dense in $\mathcal{C}(\mathbb{T})$ with respect to the supremum norm $\norm{\cdot}_\infty$, it is possible to extend the total variation norm to $\mathcal{S}'(\mathbb{T})$ and thus we obtain
\begin{equation}
    \mathcal{M}(\mathbb{T}) = \left\{ w \in \mathcal{S}'(\mathbb{T}), \sup_{\varphi \in \mathcal{S}(\mathbb{T}), \norm{\varphi}_\infty \leq 1} {\langle w, \varphi \rangle} < \infty \right\}.
\end{equation}

The dual space of $\mathcal{M}(\mathbb{T})$ with respect to the strong topology is complex and strictly larger than $\mathcal{C}(\mathbb{T})$. However, when $\mathcal{M}(\mathbb{T})$ is endowed with its \emph{weak-$^*$ topology}, its dual space is actually $\mathcal{C}(\mathbb{T})$. Then, when $\mathcal{C}(\mathbb{T})$ is respectively endowed with its \textit{weak} topology, they play symmetrical role and they are both the dual space of each other. In the following, we will keep with this \enquote{reflexive} topology, that is more convenient than the one induced by the TV-norm.\\

It is interesting to point out that the TV-norm is actually a generalization of the $\mathrm{L}^1$-norm for integrable functions. Indeed, for any $f\in L^1(\mathbb{T})\cap\mathcal{S}(\mathbb{T})$, the two norms coincide : $\norm{f}_1 = \int_\mathbb{T}{|f(t)|\ \mathrm{d}t} = \norm{f}_{TV}$. Moreover, let us introduce the $T$-periodic Dirac impulse $\Sha = \sum_{i\in\mathbb{Z}}\delta(\cdot - iT)$, coined as Dirac comb. The two following relations hold
\begin{equation*}
    \norm{\Sha}_{TV} = 1
\end{equation*}
\begin{equation*}
    \norm{\sum_{k=1}^K {\alpha_k \Sha(\cdot - t_k)}}_{TV} = \norm{\boldsymbol{\alpha}}_1
\end{equation*}
for $K\in\mathbb{N}^*$ and any $\boldsymbol{\alpha}\in\mathbb{R}^K$, $(t_1, \dots, t_K)\in\mathbb{T}^K$.

\subsection{Pseudo-Differential Operators}

We consider the \emph{linear, shift-invariant} operators that are \emph{continuous} from $\mathcal{S}'(\mathbb{T})$ to itself denoted as $\mathcal{L}_{\mathrm{SI}}\left(\mathcal{S}'(\mathbb{T})\right)$. It is possible to leverage the periodic property of the objects to perform an extensive study of these operators, as it is done in \cite{Badoual_2018}.

The same way as any generalized periodic function can be uniquely decomposed over the Fourier basis $e_k : x \mapsto e^{ikx} \in \mathcal{S}\left(\mathbb{T}\right)$ for $k\in\mathbb{Z}$, one may uniquely decompose any operator $\mathcal{D}\in\mathcal{L}_{\mathrm{SI}}\left(\mathcal{S}'(\mathbb{T})\right)$ over the same Fourier basis. That way, any operator is equivalently characterized by its Fourier sequence. This allows a theoretical expression of the operators as well as explicit numerical approximations of such operators.

\begin{figure}[H]
\makebox[\textwidth][c]{
\begin{subfigure}{0.6\textwidth}
    \centering
    \includegraphics[width=0.95\linewidth,center]{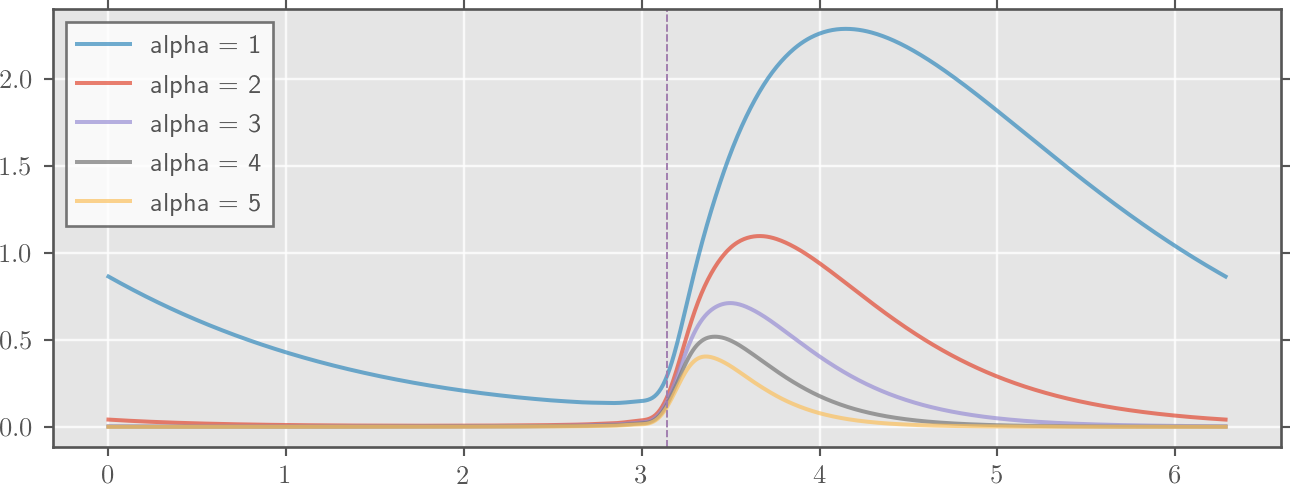}
    \caption{Exponential operator $\mathcal{D} = \left(\mathrm{D} + \alpha \mathrm{Id}\right)^2$}
    \label{fig:green_exponential}
\end{subfigure}
\begin{subfigure}{0.6\textwidth}
    \centering
    \includegraphics[width=0.95\linewidth, center]{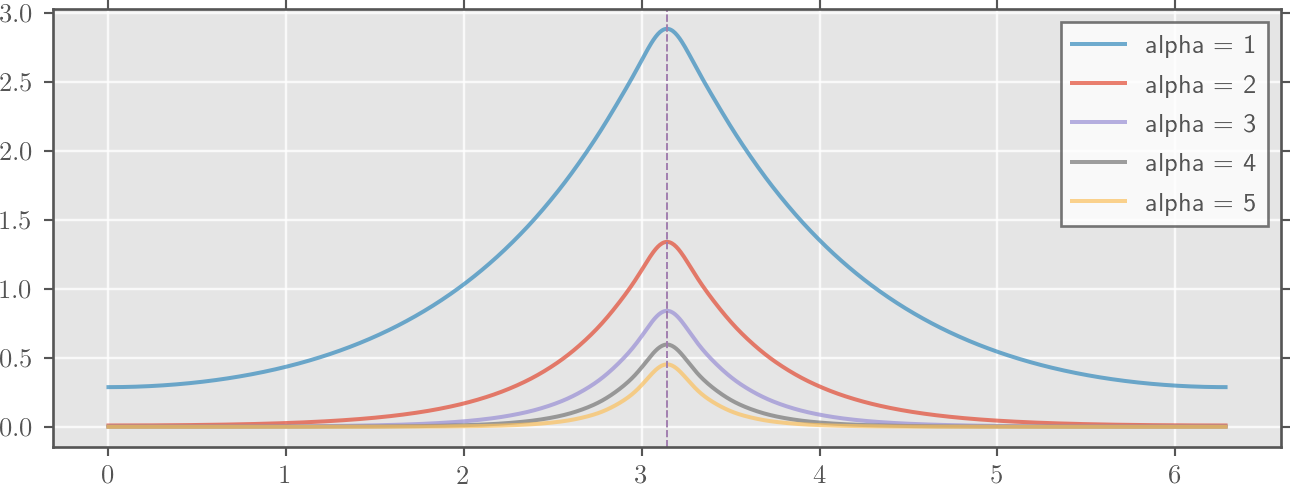}
    \caption{Sobolev operator $\mathcal{D} = \left(\alpha^2 \mathrm{Id} - \mathrm{D}^2 \right)$}
    \label{fig:green_sobolev}
\end{subfigure}
}
\caption{Green's functions with different operator parameter $\alpha$ and a shifted knot at $t_0 = \pi$}
\label{fig:green_alpha}
\end{figure}

\begin{figure}[H]
\makebox[\textwidth][c]{
\begin{subfigure}{0.6\textwidth}
    \centering
    \includegraphics[width=0.95\linewidth,center]{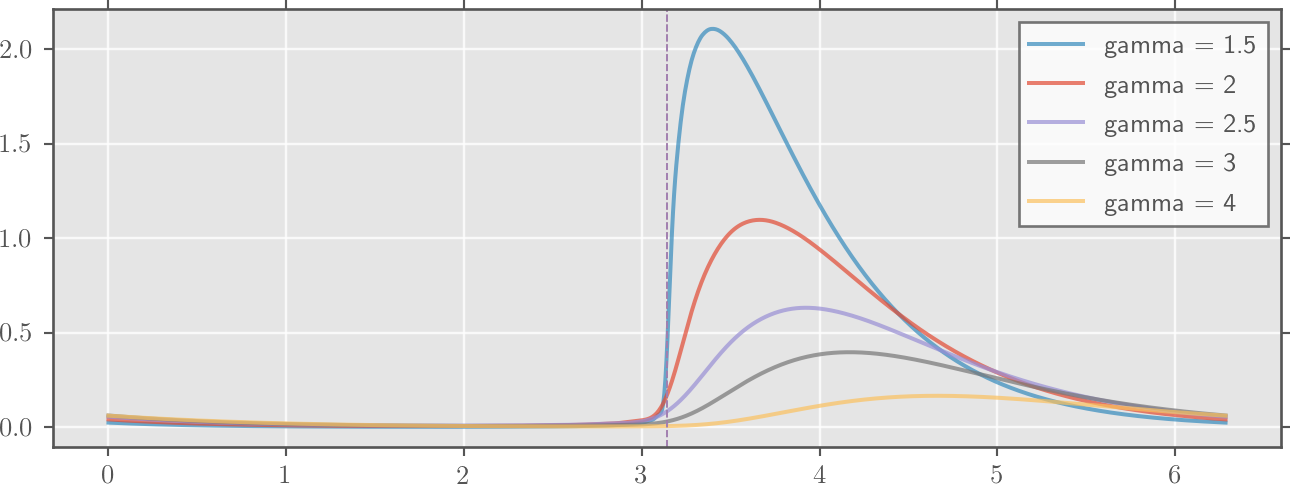}
    \caption{Exponential operator $\mathcal{D} = \left(\mathrm{D} + 2 \mathrm{Id}\right)^\gamma$}
    \label{fig:green_exponential_gamma}
\end{subfigure}
\begin{subfigure}{0.6\textwidth}
    \centering
    \includegraphics[width=0.95\linewidth, center]{Figures/greens_functions/sobolev_green.png}
    \caption{Sobolev operator $\mathcal{D} = \left(2^2 \mathrm{Id} - \mathrm{D}^2 \right)^{\gamma/2}$}
    \label{fig:green_sobolev_gamma}
\end{subfigure}
}

\caption{Green's functions with different operator exponent $\gamma$ and a shifted knot at $t_0 = \pi$}
\label{fig:green_gamma}
\end{figure}

In this project, we only consider invertible operators, such that they have trivial null space and admit an inverse operator. We can cite the exponential operator $\mathcal{D} = \left(\mathrm{D} + \alpha \mathrm{Id}\right)^\gamma$ or the Sobolev operator $\mathcal{D} = \left(\alpha^2 \mathrm{Id} - \mathrm{D}^2\right)^{\gamma/2}$ both of order $\gamma>1$. Any invertible operator $\mathcal{D}\in\mathcal{L}_{\mathrm{SI}}\left(\mathcal{S}'(\mathbb{T})\right)$ belongs to the so-called \emph{spline-admissible operators} \cite[Definition~2]{fageot2020tvbased}, and thus we define
\begin{definition}[Periodic $\mathcal{D}$-spline]
Any function $f\in\mathcal{S}'(\mathbb{T})$ is said to be a $\mathcal{D}$-spline if there exists $K\geq 0$ weights $\left(\alpha_1, \dots, \alpha_K\right) \in \mathbb{R}^K$ and $K$ positions $\left( t_1, \dots, t_K \right) \in \mathbb{T}^K$ such that
\begin{equation}
    \mathcal{D}f = \sum_{k=1}^K {\alpha_k \Sha(\cdot - t_k)}.
\end{equation}
The pairs $\left(\alpha_k, t_k \right)$ are called the \emph{innovations} of the spline while the $(t_k)$ are the \emph{knots}.
\end{definition}
\noindent As we work with invertible operators, any $\mathcal{D}$-spline $f$ can be rewritten as
\begin{equation}
    f = \mathcal{D}^{-1}\left\{ \sum_{k=1}^K {\alpha_k \Sha(\cdot - t_k)} \right\} = \sum_{k=1}^K {\alpha_k \psi_\mathcal{D}(\cdot - t_k)}
\end{equation}
where $\psi_\mathcal{D} = \mathcal{D}^{-1}\left\{\Sha\right\}$ is called the \emph{Green's function} of the operator $\mathcal{D}$, defined as follows:
\begin{definition}[Green's function]
Let $\mathcal{D}$ be a \textit{spline admissible} operator. Any function $\psi\in\mathcal{S}'\left(\mathbb{T}\right)$ such that $\mathcal{D}\psi = \Sha$ is called a \emph{Green's function} of $\mathcal{D}$. In the case of an invertible operator $\mathcal{D}$, one obtains the expression $\psi_\mathcal{D} = \mathcal{D}^{-1}\left\{\Sha\right\}$.
\end{definition}
Some examples of Green's functions with exponential and Sobolev operators are given in figures \ref{fig:green_alpha} and \ref{fig:green_gamma}. The reader should refer to \cite[Section~5]{fageot2020tvbased} for more classes of spline admissible operators.

\section{Spatial Sampling for Recovery Problems}
This project aims at performing signal approximation from (possibly noisy) measurements of the signal. The measurements we are concerned about are direct measurements, which is also known as \emph{spatial sampling}.

Let $f_0$ be a $T$-periodic real-valued function $f_0 : \mathbb{T} \rightarrow \mathbb{R}$.
Let $\mathbf{y}_0 = (y_1, \dots, y_L)$ a $L$-tuple of (potentially noisy) samples of $f_0$ such that
$$\mathbf{y}_0 = \boldsymbol{\Phi}(f_0) + \mathbf{w} = (f_0(\theta_1)+ w_1, \dots, f_0(\theta_L)+w_L)$$
for $(\theta_1, \dots, \theta_L)$ in $ \left[0, T\right]$ the spatial samples and $\mathbf{w}$ a measurement noise ( a Gaussian white noise for instance, in our context). We aim at finding a function $f\in\mathcal{S}'(\mathbb{T})$ that minimizes $E(\mathbf{y}_0, \mathbf{\Phi}(f))$ with $E$ a convex, proper, lower semi-continuous cost functional (quadratic cost $E(\mathbf{x}, \mathbf{y}) = \norm{\mathbf{x} -\mathbf{y}}^2_2$ for instance). As $\mathcal{S}'(\mathbb{T})$ is infinite dimensional, this problem needs a regularization. Indeed, without penalization, one could always consider adding an element of the null space of the measurement operator to any solution to create a new solution with the exact same functional cost. A solution commonly proposed is to use the TV norm in combination with a (pseudo-)differential operator \cite{debarre2020sparsest, Simeoni:275337}. The TV norm is known to promote sparsity in the solution. With the differential operator, the sparsity is to be searched in the generalized derivative of the solution instead of in the solution itself.

For any $\lambda > 0$ and $\mathcal{D}$ a periodic linear sample-admissible operator with trivial null space, we will focus on the following problem:
\begin{equation}
    \mathrm{Find:}\qquad f^* \quad \in \argmin_{f \in \mathcal{M}(\mathbb{T})} \norm{\mathbf{y}_0 - \mathbf{\Phi}(f)}_2^2 + \lambda \norm{\mathcal{D}f}_{TV}
    \label{eq:ill-posed-pb}
\end{equation}
The parameter $\lambda$ moderates the relative significance of both these terms.

However, written as \eqref{eq:ill-posed-pb}, two details need to be carefully addressed. First, $\norm{\mathcal{D}f}_{TV}$ might not be finite for all $f$ in $\mathcal{M}(\mathbb{T})$. It is actually possible to show that the maximal search space is 
\begin{equation}
    \mathcal{M_\mathcal{D}}(\mathbb{T}) = \left\{ f \in \mathcal{S'(\mathbb{T})}, ~ \norm{\mathcal{D}f}_{TV} < \infty\right\}
\end{equation}
which reduces the search space for solutions.
Secondly, the operator $\mathcal{D}$ needs to be \emph{sampling-admissible}, in the sense defined by \cite[Section~6.2]{fageot2020tvbased}. It amounts for the Green's function $\psi_\mathcal{D} = \mathcal{D}^{-1} \Sha \in \mathcal{S}'(\mathbb{T})$ of the operator to be continuous (which is not true for any operator). The authors of \cite{fageot2020tvbased} present some sufficient conditions for sample-admissibility, in particular focusing on the spectral growth of the operator (when it exists). Sampling-admissibility is necessary in our context as it allows working with spatial measurements. Note that it is possible to use less restrictive operators with other forms of measurements (Fourier sampling for instance).

The well-posed problem with a sampling-admissible operator $\mathcal{D}$ becomes
\begin{equation}
    \min_{f \in \mathcal{M_\mathcal{D}}(\mathbb{T})} \norm{\mathbf{y}_0 - \mathbf{\Phi}(f)}_2^2 + \lambda \norm{\mathcal{D}f}_{TV}.
    \label{eq:continuous-pb}
\end{equation}
This problem is referred to under many different names and none of them tend to make a consensus for the moment. We will retain \emph{generalized Beurling LASSO} (B-LASSO), but one might encounter \emph{functional penalized basis pursuit}, \emph{generalised continuous basis pursuit} \cite{ekanadham2011cbp} or even \emph{sparse spikes deconvolution} \cite{duval2014exact}.\\

As the linear measurements $f \mapsto \Phi_\ell(f) = f(\theta_\ell) = \int_\mathbb{T}{f(t)\ \mathrm{d}\Sha(t-\theta_\ell)}$ are linearly independent, we know from representer theorem \cite[Theorem~4]{fageot2020tvbased} that the solution set of \eqref{eq:continuous-pb} is a non empty convex set spanned by the extreme points of the form:
\begin{equation}
    f_{ext} = \sum_{k=1}^{K} \beta_k \psi_\mathcal{D}(\cdot - t_k) = \mathcal{D}^{-1}\left(\sum_{k=1}^{K} \beta_k \Sha(\cdot - t_k)\right)
    \label{eq:spline-sol}
\end{equation}
with a positive integer $K \leq L$ (the number of samples) and distinct $t_k \in \mathbb{T}$ called \textit{knots} of the spline.

Given this expression, one can make the link with the definition of \emph{sampling admissibility}. Indeed, to perform spatial sampling, it is necessary for the sampled function to be continuous, for consistency (equality between the right hand side and left hand side limits at each point). Thus, as $\psi_\mathcal{D}$ is continuous for a sampling admissible operator, such extreme points solutions \eqref{eq:spline-sol} can actually be sampled.

%% file: Chapters/chapter2.tex
\chapter{Algorithmic treatment} 

\label{chap:algo}

The generalized Beurling LASSO problem is now well defined and we know that its solution set is non empty and convex. We also know that the extreme points of this set are splines with sparse innovations. This knowledge comes really handful in a computational approximation framework, as the reconstruction task turns into estimating the knots and the intensities of those innovations. Different algorithmic strategies are presented in the following chapter, based on this consideration. They rely on known algorithms but our particular context, specifically the use of differential operators and the periodic setting, has required subtle adaptation work. This is the main part of the contribution to the project during this internship. In everything that follows, the period is set to $T=2\pi$. All the programming is based on Python programming language and the splines are simulated thanks to the package \emph{Periodispline} \cite{periodispline} available on GitHub.


\section{Equispaced Grid Approximation}

In all areas of approximation, the always increasing computing power of modern devices keeps the discretization frameworks both computationally possible but also remarkably accurate. The first step in our research is logically to bring these methods and our problem face to face. The periodicity is rather convenient here as the space to discretize is known and bounded to the period.


\subsection{Reconstruction with fixed knots}
\label{sec:fixed_knots}

Let us first present a closed problem that will help understand the discretization scheme and be used again later on. Let suppose we know a set of $I>0$ reconstruction knots $\left(k_1, \dots, k_I\right)\in\mathbb{T}^I$. We recall that the extreme points of the solution set are expressed as \eqref{eq:spline-sol}, such that any solution $f$ here would be
\begin{equation*}
    \forall t\in\mathbb{T}, \qquad f(t) = \sum_{i=1}^{I} \beta_i \psi_\mathcal{D}(t - k_i)
\end{equation*}
with unknown weights $\boldsymbol{\beta}\in\mathbb{R}^I$.

We want to plug such a solution into the problem \eqref{eq:continuous-pb}. Remembering that the Green's function satisfies $\mathcal{D} \psi_\mathcal{D} = \Sha$, the regularization term $\norm{\mathcal{D}f}_{TV}$ becomes 
\begin{equation*}
    \norm{\mathcal{D}f}_{TV} = \norm{\sum_{i=1}^{I} \beta_i \mathcal{D}\psi_\mathcal{D}(\cdot - k_i)}_{TV} = \norm{\sum_{i=1}^{I} \beta_i \Sha(\cdot - k_i)}_{TV} = \norm{\boldsymbol{\beta}}_1.
\end{equation*}
The sampling operator evaluation turns into
\begin{equation*}
    \boldsymbol{\Phi}(f) = \left(\begin{array}{c}
     f(\theta_1)  \\
     \dots \\
     f(\theta_L)
\end{array}\right) 
= \left(\begin{array}{c}
     \sum_{k=1}^{K} \beta_k \psi_\mathcal{D}(\theta_1 - t_k)  \\
     \dots \\
     \sum_{k=1}^{K} \beta_k \psi_\mathcal{D}(\theta_L - t_k)
\end{array}\right)
= \mathbf{H}\boldsymbol{\beta}
\end{equation*}
with $\mathbf{H} \in \mathbb{R}^{L \times I}$ the matrix whose entries are
\begin{equation}
     \mathbf{H}\left[\ell, i\right] = \psi_\mathcal{D}(\theta_\ell - k_i).
     \label{eq:definition_H}
\end{equation}

Then, given a set of reconstruction knots $(k_i)_{i=1, \dots, I}$, the B-LASSO problem is transformed into a regular LASSO problem whose unknown is the vector of weights $\boldsymbol{\beta}\in\mathbb{R}^I$. We obtain the following convex problem:
\begin{equation}
    \min_{\boldsymbol{\beta} \in \mathbb{R}^I} \norm{\mathbf{y} - \mathbf{H} \boldsymbol{\beta}}_2^2 + \lambda \norm{\boldsymbol{\beta}}_{1}
    \label{eq:discrete-pb}
\end{equation}
This problem can be solved quite efficiently nowadays thanks to recent proximal algorithms, even for large values of $I$. We use the accelerated proximal gradient descent algorithm, also known as FISTA \cite{amir2009fista}, for computational treatment.


\subsection{Discretization of the period}

We have seen that when the reconstruction knots are known, it is straightforward to determine the weights and thus solve the problem. However, there is no prior for the knots to be known beforehand. Let us remind that we are looking for solution in the shape of \eqref{eq:spline-sol}. One shall notice that it is equivalent to consider a knot associated to a null weight or no knot at all. Based on this remark, an idea is to discretize the whole interval $\mathbb{T} = \left[0, 2\pi \right]$ with a large number $N$ of equispaced knots $\hat{t}_n = (n-1) T / N$ for $n=1, \dots, N$ and to search for the $N$ associated weights $\beta_n$. The final number of knots will be the number of non zero coefficients of the weight vector $\boldsymbol{\beta}$. This method has proven to be efficient in the non periodic setting, even with non invertible differential operator \cite[Section~V.B]{gupta2018continuous}. \\

Let suppose the source signal is a function $f_0 = \mathcal{D}^{-1}m_0$ with $m_0$ being a sparse signal with $K$ knots. Sparse signal refers to a stream of Dirac defined as linear combination of periodic Dirac impulse, such that $m_0 = \sum_{k=1}^K {\alpha_k \Sha(\cdot - t_k)}$. We want $N$ to be significantly larger than $K$ to make sure that the grid is fine enough and consequently is able to precisely approximate the solution. We obtain the same problem as \eqref{eq:discrete-pb} with a matrix $\mathbf{H}\in\mathbb{R}^{L\times N}$ where $N\gg K$. The explicit expression of this matrix is:
\begin{equation}
    \mathbf{H} = \begin{pmatrix}
    \psi_\mathcal{D}(\theta_1 - \hat{t}_1) & \cdots & \psi_\mathcal{D}(\theta_1 - \hat{t}_N) \\
    \vdots & \ddots & \\
    \psi_\mathcal{D}(\theta_L - \hat{t}_1) & \cdots & \psi_\mathcal{D}(\theta_L - \hat{t}_N)
    \end{pmatrix}.
\end{equation}
Let $\mathbf{y}_0 = \mathbf{\Phi}(f_0)$, the problem of interest inow becomes:
\begin{equation}
    \min_{\boldsymbol{\beta} \in \mathbb{R}^N} \norm{\mathbf{y}_0 - \mathbf{H} \boldsymbol{\beta}}_2^2 + \lambda \norm{\boldsymbol{\beta}}_{1}
    \label{eq:grid_pb}
\end{equation}
with a vector $\boldsymbol{\beta}$ of size $N$. It completely fits with the application setting of the \emph{Accelerated Proximal Gradient Descent algorithm} (APGD), and more precisely to the special case of FISTA (see \cite[Algorithm~7.6]{Simeoni:275337}).

\subsubsection{APGD}
Consider the problem
\begin{equation}
    \mathrm{Find :}\quad\mathbf{u}^\star\ \in \quad \argmin_{\mathbf{u}\in\mathbb{R}^N}{\mathcal{F}(\mathbf{u}) + \mathcal{G}\left(\mathbf{u}\right)}
\end{equation}
with $\mathcal{F} : \mathbb{R}^N \to \mathbb{R}$ a \emph{convex} and \emph{differentiable} function with $\gamma$-\emph{Lipschitz continuous gradient} and $\mathcal{G} : \mathbb{R}^N \to \mathbb{R}$ a \emph{proper}, \emph{convex} and \emph{lower semi-continuous} function with an easy to compute proximal operator. The difficulty come from the fact that $\mathcal{G}$ is not differentiable. The APGD algorithm works as follows:
\begin{algorithm}[H]
\caption{Accelerated Proximal Gradient Descent $(\boldsymbol{\beta}_0, \tau, \delta, \varepsilon$)}
\label{alg:apgd}
\begin{algorithmic}[1]
\State Initialization : $(\mathbf{z}_0, \mathbf{x}_0) \leftarrow \left(\boldsymbol{\beta}_0, \boldsymbol{\beta}_0 \right)$
\For {$n \geq 1$}
    \State $\mathbf{z}_n = \mathbf{prox}_{\tau \mathcal{G}}\left(\mathbf{x}_{n-1}-\tau \nabla\mathcal{F}(\mathbf{x}_{n_1})\right)$
    \State $\mathbf{x}_n = \mathbf{z}_n + \frac{n-1}{n+\delta}(\mathbf{z}_n - \mathbf{z}_{n-1})$
    \If{$\norm{\mathbf{x}_n - \mathbf{x}_{n-1}}_2 \leq \varepsilon \norm{\mathbf{x}_{n-1}}_2$}
        \State $\mathbf{x}_n$ solution of the problem. Stop.
    \Else
        \State $n\leftarrow n+1$
    \EndIf
\EndFor \\
\Return $(\mathbf{x}_n)_{n\in\mathbb{N}}$
\end{algorithmic}
\end{algorithm}
\noindent In our case, $\mathcal{F}$ is the data fidelity term and $\mathcal{G}$ is the regularization term. Let us get into details about the steps and parameters involved.

\subsubsection{Initialization}
The initialization step is straightforward. If the value of the first iterate is not specified, it is set to $\boldsymbol{\beta}_0 = \mathbf{0}_N$. The two iterable variables $\mathbf{x}_n$ and $\mathbf{z}_n$ are then initialized with value $\boldsymbol{\beta}_0$. During our experiments, we have used $\mathbf{0}_N$ as initial iterate.

\subsubsection{Proximal step}
The APGD algorithm makes use of the proximal operator associated to the $1$-norm. It actually admits a closed-form expression and the following equality holds \cite{Simeoni:275337}:
\begin{equation}
\begin{aligned}
    \forall \alpha\in\mathbb{R}, \, \forall \mathbf{x}\in\mathbb{R}^N,\quad \mathbf{prox}_{\alpha\norm{\cdot}_1} (\mathbf{x}) =& \argmin_{\mathbf{z}\in\mathbb{R}^N} \frac{1}{2\alpha} \norm{\mathbf{x}-\mathbf{z}}_2^2 + \norm{\mathbf{z}}_1\\
    =& \left(\mathrm{soft}_\alpha(x_i)\right)_{i=1, \dots, N}
\end{aligned}
\end{equation}
with $\mathrm{soft}_\alpha$ being the soft thresholding function
\begin{equation}
    \forall \alpha\in\mathbb{R}, \, \forall x \in\mathbb{R},\quad \mathrm{soft}_\alpha(x) = \max(|x| - \alpha, 0)\ \mathrm{sgn}(x).
    \label{eq:soft_definition}
\end{equation}
The function $\mathrm{sgn}$ is the \emph{sign function} such that $x\cdot\mathrm{sgn}(x)>0$ for $x\in\mathbb{R}^*$, $\mathrm{sgn}(0)=0$. These simple expressions lead to an easy computation of the proximal step.

\subsubsection{Step size}
The data fidelity term $\mathcal{F} : \boldsymbol{\beta} \rightarrow \norm{\mathbf{y}_0-\mathbf{H}\boldsymbol{\beta}}^2_2$ is a smooth function with $\gamma$-Lipschitz continuous derivative. We can easily show that
\begin{equation*}
    \gamma = 2\, \norm{\mathbf{H}^T \mathbf{H}}_2 = 2\, \norm{ \mathbf{H}}^2_2 = 2 \,\lambda_{max}(\mathbf{H}^T \mathbf{H})
\end{equation*}
where $\lambda_{max}(\mathbf{M})$ is the maximum eigenvalue of the matrix $\mathbf{M}$.
Then, the step size $\tau$ is chosen to maximize the convergence speed according to $\tau \leq \frac{1}{\gamma}$. In our case, we chose the greatest value possible, so that $\tau = \frac{1}{\gamma}$.

\subsubsection{Iteration}
Then, we obtain the following expression for the iteration step:
\begin{equation}
    \left\{
    \begin{aligned}
        \mathbf{z}_{n+1} \ &\leftarrow \ \mathbf{soft}_{\tau\lambda}(\mathbf{x}_n - 2 \tau \mathbf{H}^T(\mathbf{H}\mathbf{x}_n - \mathbf{y})) \\
        \mathbf{x}_{n+1} \ &\leftarrow \ \mathbf{z}_{n+1} + \frac{n-1}{n + \mathrm{\delta}} (\mathbf{z}_{n+1} - \mathbf{z}_n)
    \end{aligned}
    \right. 
    \label{eq:proximal_algo}
\end{equation}
The parameter $\delta$ is an empirical acceleration parameter, whose value is usually taken around $\delta=75$. This choice is motivated by the good results obtained in \cite{liang2018activity}, which highlights significant acceleration for $\delta\in\left[50, 100\right]$.

\subsubsection{Stopping criterion}
When the norm of the difference between two iterates of $\mathbf{x}_n$ becomes small enough compared to the norm of the previous iterate, one considers that the process has converged and can stop. This is called the \emph{relative improvement}, ruled by a parameter $\varepsilon$ that is usually set quite low (for instance we use $\varepsilon \approx 10^{-4}$).

\begin{remark}[Value of $\lambda$]
In such a problem, the \textit{regularization parameter} $\lambda$ is at the heart of sparse recovery. Theoretically, as long as $\lambda > 0$, the problem is well-posed, and its value is not considered any further. However, when running actual code, setting the value of $\lambda$ is crucial to exploit sparsity promoting algorithms to their best potential. Practically, $\lambda$ (or more precisely a $\lambda$ related quantity) answers to question: \enquote{To what extent am I ready to pay (in sparsity) to allow more accurate approximation ?} Low values mean sparsity is cheap, whereas high values lead to extremely sparse solution, up to a critical point where even a unique parameter solution costs too much and the only affordable trade-off is a null approximation.

In the current case of equations \eqref{eq:proximal_algo}) starting with an initial value $\mathbf{x}_0 = \mathbf{0}$, one obtains: $$\mathbf{z}_1 = \mathbf{soft}_{\tau \lambda}(2\tau\mathbf{H}^T\mathbf{y})$$ which leads to $\mathbf{0}$ when $$\lambda \geq  2\max_{i=1, \dots, N}|\mathbf{H}^T \mathbf{y}|_i = 2\norm{\mathbf{H}^T \mathbf{y}}_\infty.$$  By linear combinations, the next iterates would then be stuck at $\mathbf{0}$. Even with another non null initial value $\mathbf{x}_0$, we notice the same collapsing behavior of the iterates towards $\mathbf{0}$ after a certain number of iterations.

Consequently, it seems interesting to scale the value of $\lambda$ to $\norm{\mathbf{H}^T \mathbf{y}}_\infty$, which is a practical measurement that depends both on the sampling operator (that includes sampling positions and reconstruction spline) and the sampled values. More precisely, $\lambda$ should verify :
\begin{equation}
    \lambda \in \left[0, \lambda_{max}\right], \qquad \text{with} \qquad \lambda_{max} = 2 \norm{\mathbf{H}^T \mathbf{y}}_\infty
\end{equation}
We usually set $\lambda = \sigma \norm{\mathbf{H}^T \mathbf{y}}_\infty$ with $\sigma\in\mathbb{R}$ close to $0.1$ .
\end{remark} 

\begin{remark}[Sparsity index]
The solution obtained following this procedure does not necessarily meet the sparsity index guaranteed by the representer theorem, \emph{i.e.} sparsity index of the solution lower than $L$. \cite[Theorem 7]{gupta2018continuous} ensures that one can fulfill this sparsity constraint using the simplex algorithm over the linear problem
\begin{equation}
    \min_{\boldsymbol{\beta}\in\mathbb{R}^N} \norm{\boldsymbol{\beta}}_1 \quad \text{subject to} \quad H \boldsymbol{\beta} = H \boldsymbol{\beta ^*}
\end{equation}
with $\boldsymbol{\beta^*}$ solution of the FISTA algorithm. It is interesting to point out, however we will not focus on this point any further, as we are more interested in the quality on the approximation than the sparsity of the solutions.
\end{remark}


\section{Finite Rate of Innovation}

Another relevant type of algorithms for the approximation of problem \eqref{eq:continuous-pb} are referred to as \emph{finite rate of innovation} algorithms (FRI) \cite{blu2008fri}. In the definition, the FRI is a way of measuring the quantity of information contained per time interval for a signal. For instance, to approximate a signal such as \eqref{eq:spline-sol}, it is necessary to estimate the $t_k$ and $\beta_k$ values, which leads to a FRI of $\rho = \frac{2K}{T}$, with $T=2\pi$ the period of the signal. Practically, for a given periodic signal, it is necessary to perform at least the FRI number $\gamma$ of observations per unit of time during a period to collect enough information and reconstruct the signal.

Algorithms coined as \textit{FRI} decouple the problem of reconstructing a spline into two steps. The first step consists in estimating the knots of the spline, the second one in estimating the weights given the knots position. As we have seen, this second step can efficiently be carried out thanks to FISTA as it amounts to the reconstruction problem with fixed knots \eqref{eq:discrete-pb}. Let us then focus to the estimation of the positions of the knots. We apply a refined version of FRI presented in \cite{simeoni2020cpgd} called CPGD for \emph{Cadzow Plug and Play Gradient Descent}, that have been specifically designed to reduce the noise dependency of the original FRI methods.

The FRI framework is based on the Fourier analysis of the input signal. In order to retrieve the positions of the knots, the algorithms first estimate the Fourier series coefficients of the input signal (up to some order). We first present how spatial samples of the signal can lead to these coefficients, then we analyze how the knowledge of the Fourier series allow the computation of the knots.

\subsection{Linking Fourier coefficients to the samples}

In order to comply with the formulation of the problem solved by CPDG algorithm, we need to link the spatial samples $\mathbf{y}$ to the truncated Fourier series of a Dirac stream linked to the input signal. Let us reformulate our problem so that this requirement is fulfilled. 

Let suppose the source signal $f_0$ is a $K_0$ sparse signal such that
\begin{equation}
    \mathbf{y}_0 = \mathbf{\Phi}(f_0) = \mathbf{\Phi}(\mathcal{D}^{-1}m_0) \qquad \mathrm{and}\qquad m_0 = \sum_{k=1}^{K_0}{\beta_{0, k} \Sha( \cdot - t_{0,k})}
    \label{eq:source_sparse}
\end{equation}
that we want to approximate with a $K$ sparse signal
\begin{equation}
    f = \sum_{k=1}^K{\beta_k \psi_\mathcal{D}(\cdot - t_k)} = \mathcal{D}^{-1}\left(\sum_{k=1}^{K}{\beta_k \Sha( \cdot - t_k)}\right) = \mathcal{D}^{-1}m
    \label{eq:sol_sparse}
\end{equation}
with the corresponding stream of periodic Dirac $m$.
It is possible to compute the infinite Fourier transform series of $m$, noted $\hat{m} = \mathcal{F}(m)$, as:
\begin{equation}
     \hat{m}_n = \left(\mathcal{F}(m)\right)_n = \sum_{k=1}^K \beta_k\exp(-j2\pi n\frac{t_k}{T}), \quad n\in\mathbb{Z}.
\end{equation}
Computationally, we will consider the truncated Fourier series $\mathbf{\hat{m}}$ up to coefficient $M\in\mathbb{N}$ as an approximation of $\mathcal{F}(m)$:
 \begin{equation}
    \mathbf{\hat{m}} = \mathcal{P}_M \mathcal{F} (m) \in \mathbb{C}^{2M+1}
    \label{eq:stream_dirac}
 \end{equation}
with:
\begin{equation}
    \mathcal{P}_M : \left\{\begin{aligned}
        \mathbb{C}^\mathbb{Z} \quad &\to \quad \mathbb{C}^{2M+1} \\
        \quad (a_n)_{n\in\mathbb{Z}} \quad &\mapsto \quad (a_{-M}, \dots, a_M)
    \end{aligned} \right.
\end{equation}
With the Fourier basis noted as ${e}_n(t) = \exp(j2\pi n t/T)$, the approximated reverse operation of building back the signal $m$ using its truncated Fourier transform series becomes
\begin{equation}
    m \approx \mathcal{F}^{\,\dag}_M (\mathbf{\hat{m}}) := \sum_{n=-M}^M{\hat{m}_n} {e}_n.
\end{equation}
Naturally, the greater $M$ the better the approximation holds.
Then, one can approximate the $\mathbf{y}$ noiseless samples of the signal as:
 \begin{equation}
     \mathbf{y} = \Phi(f) = \Phi \mathcal{D}^{-1} (m) \approx \Phi \mathcal{D}^{-1} \mathcal{F}^{\,\dag}_M(\mathbf{\hat{m}}).
 \end{equation}
For a given $M\in\mathbb{N}$, calculations lead to
\begin{align*}
    \Phi \mathcal{D}^{-1} \mathcal{F}^{\,\dag}_M(\mathbf{\hat{m}}) ={ }& \Phi \mathcal{D}^{-1}\left(\sum_{n=-M}^M{\hat{m}_n} {e}_n \right) \\
    ={ }& \Phi \left( \sum_{n=-M}^M{\hat{m}_n} \widehat{\mathcal{D}^{-1}}\left[n\right]{e}_n \right) \\
    ={ }& \Phi \left( \sum_{n=-M}^M{\hat{m}_n} \widehat{\psi_\mathcal{D}}\left[n\right]{e}_n \right) \\
    ={ }& \left( \sum_{n=-M}^M{\hat{m}_n} \widehat{\psi_\mathcal{D}}\left[n\right]{e}_n(\theta_\ell) \right)_{l=1, \dots, L}
\end{align*}
using the fact that the Fourier sequence of the Green's function $\psi_\mathcal{D}$ of a spline admissible invertible operator $\mathcal{D}$ coincides with the Fourier sequence of $\mathcal{D}^{-1}$. Eventually, the following approximation holds :
\begin{equation}
    \mathbf{y} \approx \mathbf{G} \mathbf{\hat{m}}
    \label{eq:y_Ghatx}
\end{equation}  
with $\mathbf{G}\in\mathbb{C}^{L\times 2M+1}$ the following matrix
\begin{equation}
    \mathbf{G}\left[l, n \right] = \widehat{\psi_\mathcal{D}}\left[n\right]{e}_n(\theta_\ell) .
\end{equation}
Note that $\mathbf{G}$ can be rewritten as a product of diagonal and Vandermonde matrices such that:
\begin{equation}
    \mathbf{G} = \textbf{Diag}\left[ {e}_{-M}(\theta_\ell) \right]_\ell \textbf{Vand}^{2M+1}\left[{e}_{1}(\theta_\ell)\right]_\ell \textbf{Diag}\left[\widehat{\psi_\mathcal{D}}\left[n\right]\right]_n
    \label{eq:def_G}
\end{equation}
or, with explicit matrix notations :
\begin{equation*}
    \mathbf{G} =
     \begin{bmatrix}
        \mathbf{e}_{-M}(\theta_1) & & 0\\
        & \ddots & \\
        0 & & \mathbf{e}_{-M}(\theta_L)
    \end{bmatrix}
    \begin{bmatrix}
        1 & {e}_{1}(\theta_1) & {e}_{1}(\theta_1)^2 & \dots & {e}_{1}(\theta_1)^{2M}\\
        1 & {e}_{1}(\theta_2) & {e}_{1}(\theta_2)^2 & \dots & {e}_{1}(\theta_2)^{2M}\\
        \vdots & \vdots & \vdots & \ddots &\vdots \\
        1 & {e}_{1}(\theta_L) & {e}_{1}(\theta_L)^2 & \dots & {e}_{1}(\theta_L)^{2M}
    \end{bmatrix}
    \begin{bmatrix}
        \widehat{\psi_\mathcal{D}}\left[-M\right] & & 0\\
        & \ddots & \\
        0 & & \widehat{\psi_\mathcal{D}}\left[M\right]
    \end{bmatrix}
\end{equation*}

Finally, the noisy version of equation \eqref{eq:y_Ghatx} linearly links the measurements $\mathbf{y}$ to the Fourier series coefficients of the stream of Dirac \eqref{eq:stream_dirac} and transforms the problem to the CPGD context presented in \cite{simeoni2020cpgd}:
\begin{equation}
    \mathbf{y} \approx \mathbf{G} \mathbf{\hat{x}} + \mathbf{w}
    \label{eq:noisy_FRI}
\end{equation}
with $\mathbf{w}$ a random centered noise.


\subsection{Estimating knots positions}

\subsubsection{Annihilating filter}
In order to estimate the knots position, the FRI framework leverages a tool called annihilating filter. As presented in \cite{simeoni2020cpgd}, this filter is represented by a vector $\mathbf{h}\in\mathbb{C}^{K+1}$ such that:
\begin{itemize}
    \item if we consider the $0$-padded sequence $h \in \mathbb{C}^\mathbb{Z}$ of $\mathbf{h}$, the convolution product $\hat{x} * h$ with the (infinite) sequence of Fourier coefficients of $m$ is null: $(\hat{m} * h)_n = 0\text{,}\quad n\in\mathbb{Z}$,
    \item the roots of the \textit{z}-transform of $\mathbf{h}$ are $u_k = \exp(-j 2\pi t_k/T)$, in one-to-one relation with the knots positions $t_k$.
\end{itemize}
From the convolution product $(\hat{m} * h) \in \mathbb{C}^\mathbb{Z}$, one can extract the $2M-K+1$ coefficients that would result from the finite-sequences convolution product between $\mathbf{\hat{m}}$ and $\mathbf{h}$, which can be computed as:
\begin{equation}
    (\mathbf{\hat{m}} * \mathbf{h}) = \mathrm{T}_K(\mathbf{\hat{m}}) \mathbf{h}
\end{equation}
with the \textit{Toeplitzification operator} $\mathrm{T}_P$ defined for any integer $0\leq P \leq M$ as:
\begin{equation}
    \mathrm{T}_P \colon \left\{ \begin{aligned}
        \quad \mathbb{C}^{2M+1}\quad & \to \mathbb{C}^{({2M}-P+1)\times(P+1)} \\
        \mathbf{v}\quad & \mapsto \left[ \mathrm{T}_P(\mathbf{v})\right]_{i, j} := x_{-M + P + i - j}
    \end{aligned} \right.
\end{equation}
It has been shown (in \cite{blu2008fri}) that, as long as $M \geq K$, the matrix $\mathrm{T}_K(\mathbf{\hat{m}})$  has rank $K$, such that the equation
\begin{equation}
    \mathrm{T}_K(\mathbf{\hat{m}})\, \mathbf{\tilde{h}} = \mathbf{0} \ \in\mathbb{C}^{2M-K+1}
    \label{eq:total_least_squares}
\end{equation}
always admits non trivial solutions for $\mathbf{\tilde{h}}\in\mathbb{C}^{K+1}$. Actually, when $P$ is such that $K\leq P \leq M$, the matrix $\mathrm{T}_P(\mathbf{\hat{m}}) \in \mathbb{C}^{(2M-P+1)\times(P+1)}$ is still of rank $K$, which leads to the conclusion that
\begin{equation}
    \mathrm{T}_P(\mathbf{\hat{m}})\, \mathbf{\tilde{h}} = \mathbf{0} \ \in\mathbb{C}^{2M-P+1}
    \label{eq:total_least_squares_P}
\end{equation}
still admits non trivial solutions for $\mathbf{\tilde{h}}\in\mathbb{C}^{K+1}$ (note the change of $K$ to $P$).

\subsubsection{Fourier series estimation}
Then, as the $\mathbf{G}$ matrix presented earlier \eqref{eq:def_G} is injective, the authors propose to approximate the $\mathbf{\hat{m}}$ Fourier coefficients of $m$ as the solution of the minimization problem:
\begin{equation}
    \min_{\mathbf{z}\in\mathbb{C}^{2M+1}}  \norm{\mathbf{G}\mathbf{z} - \mathbf{y}}_2^2 \quad \text{subject to} \quad \left\{\text{rank}\  \mathrm{T}_K(\mathbf{z}) \leq K\right\}
    \label{eq:rank_constrained_pb_K}
\end{equation}

However, as a sort of oversampling to reduce noise influence, this problem is solved using a slightly more restrictive constraint. Using the fact that $\mathrm{T}_P(\mathbf{\hat{m}})$ has also rank $K$ for $K \leq P \leq M$, the constraint is replaced by $\left\{\mathrm{rank}\  \mathrm{T}_M(\mathbf{z}) \leq K\right\}$ (using $P=M$, bigger matrix), which yields the best empirical performances, leading to:
\begin{equation}
    \min_{\mathbf{z}\in\mathbb{C}^{2M+1}}  \norm{\mathbf{G}\mathbf{z} - \mathbf{y}}_2^2 \quad \text{subject to} \quad \left\{\mathrm{rank}\  \mathrm{T}_M(\mathbf{z}) \leq K\right\}
    \label{eq:rank_constrained_pb}.
\end{equation}

\subsubsection{Knots recovery}
So far, the rank constraint has been used to ensure the existence of a non trivial vector $\mathbf{h}$ but this vector is not calculated when solving \eqref{eq:rank_constrained_pb}. Indeed, it is not necessary to solve the problem and retrieve the truncated Fourier series coefficients $\mathbf{\hat{m}}$. Nonetheless, it takes part in the calculation of the knots locations (in relation with the roots of the $\textit{z}$-transform). Consequently, once $\mathbf{\hat{m}}$ has been estimated as the solution of problem \eqref{eq:rank_constrained_pb}, $\mathbf{h}$ is found as the solution of equation \eqref{eq:total_least_squares} (using \emph{Total Least Squares} for example). Note that here, $P$ is set to the value $K$ which leads to $K+1$ coefficients of $\mathbf{h}$ and so $K$ locations of the knots as expected. Solving the equation with a different value of P would lead to an inaccurate number of solutions.


\subsection{Solving the Fourier Series inverse problem}

Due to its rank constraint, the problem \eqref{eq:rank_constrained_pb} is quite specififc and needs to be treated carefully. It is reformulated in an unconstrained way using the indicator function $\iota$ whose values live in $\left\{0, +\infty\right\}$. Let $\mathcal{H}_K$ be the set of complex matrices of rank at most $K$ :
\begin{equation}
    \mathcal{H}_K := \left\{\mathbf{M}\in\mathbb{C}^{(2M+1-P)\times(P+1)}\ |\ \text{rank}\ \mathbf{M}\leq K \right\}
\end{equation}
The unconstrained version of \eqref{eq:rank_constrained_pb} becomes
\begin{equation}
    \label{eq:unconstrained_pb}
    \min_{\mathbf{z}\in\mathbb{C}^(2M+1)}  \norm{\mathbf{G}\mathbf{z} - \mathbf{y}}_2^2 + \iota_{\mathcal{H}_K}(\mathrm{T}_M(\mathbf{z}))
\end{equation}
 which is the sum of a \textit{convex} and \textit{differentiable} term (data fidelity term) and a \textit{non-convex} (as $\mathcal{H}_K$ is not convex) and \textit{non differentiable term} (constraint term).\\

As the data fidelity term is $\gamma$-\textit{Lipschitz continuously differentiable}, it is possible to optimize \eqref{eq:unconstrained_pb}) using \textit{proximal gradient descent} (PGD). The Lipschitz constant is given by :
\begin{align}
    \gamma =& 2\norm{\mathbf{G}^* \mathbf{G}} \nonumber \\
    =& 2 \lambda_{max}\left( \mathbf{G}^* \mathbf{G}  \right)
\end{align}
that can be computed as  the greatest eigenvalue of the matrix. The non convexity of the constraint term does not ensure that the proximal operator would return unique solutions, but it is not a problem, as any proximal point can be used interchangeably at this step.

However, the proximal step implies at each iteration to solve the following non convex optimization problem:
\begin{equation}
    \boldsymbol{\tilde{z}} \in \mathbf{prox}_{\tau B} (\boldsymbol{u}) = \argmin_{ \boldsymbol{z} \in \mathbb{C}^2M+1}{ \left\{ \frac{1}{2\tau}\norm{\boldsymbol{u}-\boldsymbol{z}}^2 +  \iota_{\mathcal{H}_K}(\mathrm{T}_P(\boldsymbol{z})) \right\} }, \qquad \boldsymbol{u}\in\mathbb{C}^{2M+1}
\end{equation}
with $B(\cdot) = \iota_{\mathcal{H}_K}(T_P(\cdot))$.
The authors of \cite{simeoni2020cpgd} introduce the diagonal matrix $\mathbf{\Gamma}\in\mathbb{C}^{(2M+1)\times (2M+1)}$ as :
 \begin{equation*}
     \Gamma_{i, i} = \min{(i, \,P+1,\, 2M+2-i)}, \quad i=1, \dots, 2M+1 .
 \end{equation*}
and the weighting matrix $\mathbf{W}\in\mathbb{C}^{(2M+1-P)\times (P+1)}$ as:
\begin{equation*}
    \mathbf{W} = \mathrm{T}_P \mathrm{Diag}\left(\mathbf{\Gamma}^{-1/2}\right)
\end{equation*}
such that the proximal operator becomes:
\begin{equation}
    \mathbf{prox}_{\tau H}(\boldsymbol{x}) = \mathrm{T}_P^\dag \Pi^\mathbf{W}_{\mathbb{T}_P \cap \mathcal{H}_K} \mathrm{T}_P(\boldsymbol{x})
\end{equation}
where $\Pi^\mathbf{W}_{\mathbb{T}_P \cap \mathcal{H}_K}$ is the projection onto $\mathbb{T}_P \cap \mathcal{H}_K$ with respect to the $\mathbf{W}$-weighted Frobenius norm, $\mathbb{T}_P$ being the space of Toeplitz matrices of size $(2M+1-P)\times(P+1)$. This operator does not admit any simple closed-form expression, however it can efficiently be approximated by the \textit{method of alternating projections} (MAP):
\begin{equation}
    \Pi^\mathbf{W}_{\mathbb{T}_P \cap \mathcal{H}_K} \approx \left[\Pi_{\mathbb{T}_P} \Pi^\mathbf{W}_{\mathcal{H}_K}\right]^n
\end{equation}
for some $n\in\mathbb{N}$ big enough.
We finally detail the expressions of the involved projection operators:
\begin{itemize}
    \item the projection onto $\mathbb{T}_P \subset \mathbb{C}^{(2M+1-P)\times(P+1)}$ of rectangular Toeplitz matrices adimts a simple closed-form expression:
    \begin{equation}
        \Pi_{\mathbb{T}_P} = \mathrm{T}_P \mathrm{T}^\dag_P = \mathrm{T}_P \mathbf{\Gamma}^{-1} \mathrm{T}_P^*
    \end{equation}
    \item the projection onto $\mathcal{H}_K$ of matrices of rank at most $K$ with respect to the $\mathbf{W}$-weighted norm is actually more complex, but can be approximated to the first order with the regular projection onto $\mathcal{H}_K$, such that :
    \begin{equation}
        \Pi^\mathbf{W}_{\mathcal{H}_K} (\mathbf{X}) \approx \Pi_{\mathcal{H}_K} (\mathbf{X})  = \mathbf{U}\mathbf{\Lambda}_K\mathbf{V}^H, \quad X\in\mathbb{C}^{(2M+1-P)\times(P+1)}
    \end{equation}
    where $\mathbf{X} = \mathbf{U}\mathbf{\Lambda}\mathbf{V}^H$ is the singular value decomposition of $\mathbf{X}$ and $\mathbf{\Lambda}_K$ is the diagonal matrix of sorted singular values truncated to the $K$ strongest ones.
\end{itemize}

All the operators involved at any step of the algorithm are either expressed or approximated with simple closed-form expressions. The final version is then fully usable and could be summarized as in \cref{alg:cpgd}.
\begin{algorithm}
\caption{CPGD $(K, P, M, \varepsilon$)}
\label{alg:cpgd}
\begin{algorithmic}[1]
\State Solve \eqref{eq:unconstrained_pb} using PGD :
\State Initialize : $(\mathbf{z}_0, \mathbf{x}_0) \leftarrow \left(\boldsymbol{0}_{2M+1}, \boldsymbol{0}_{2M+1} \right)$
\For {$n \geq 1$}
    \State $\mathbf{z}_n = \mathbf{prox}_{\tau \mathcal{G}}\left(\mathbf{x}_{n-1}- 2 \tau \mathbf{G}^*(\mathbf{G x}_{n_1} - \mathbf{y}_0)\right)$
    \State $\mathbf{x}_n = \mathrm{T}_M^\dag \Pi^\mathbf{W}_{\mathbb{T}_P \cap \mathcal{H}_K} \mathrm{T}_M(\boldsymbol{\mathbf{z}_n})$
    \If{$\norm{\mathbf{x}_n - \mathbf{x}_{n-1}}_2 \leq \varepsilon \norm{\mathbf{x}_{n-1}}_2$}
        \State $\mathbf{x}_n$ solution of the problem. Stop.
    \Else
        \State $n\leftarrow n+1$
    \EndIf
\EndFor
\State Find $\mathbf{h}\in\mathbb{C}^{K+1}$ by solving \eqref{eq:total_least_squares_P} with TLS
\State Find the roots $(u_k)_{k=1, \dots, K}$ of the $z$-transform of $\mathbf{h}$
\State Deduce the $K$ knot positions $(t_k)_{k=1, \dots, K}$
\end{algorithmic}
\end{algorithm}

\begin{remark}[On the number M of Fourier coefficients]
    It is important to note that the maximum number of innovations that can be recovered depends on the practical implementation of the Green functions. Indeed, CPGD aims at finding $\mathbf{\hat{m}}$ which is composed of $2M+1$ coefficients of the Fourier series of the input signal. We used the python package Periodispline \cite{periodispline} which also leverages the Fourier series decomposition to define and evaluate the Green's functions. In most of the cases, there is no reason for the Fourier series of the considered Green's functions to be bandlimited: they theoretically admit an infinite number of coefficients. However, this is not the case for their actual implementation. For finite storage reason, the series of a given Green's function $\psi_\mathcal{D}$ is truncated to the first main consecutive coefficients, up to a \textit{cutoff frequency} (that preserves most of the energy of the function, set as a given percentage of the maximal absolute value frequency). Let $M_c$ be the index of the highest frequency stored for the function's representation series. The extreme point solution \eqref{eq:spline-sol} is expressed as a linear combination of shifted versions of the reconstruction Green's function, and thus is implemented with a finite bandwidth in the Fourier frequency space. It leads to the conclusion that one shall not expect the CPGD algorithm to accurately reconstruct more than $2M_c+1$ coefficients of the periodic Dirac stream $x$, practically summarized as $M \leq M_c$.
    \noindent Stated another way, one of the two following statements shall hold:
    \begin{itemize}
        \item  Given a reconstruction Green's function $\psi_\mathcal{D}$ (or equivalently a regularization operator $\mathcal{D}$), the maximum frequency index of the Fourier series to reconstruct $M$ is bounded by the \textit{cutoff frequency} of $\psi_\mathcal{D}$.
        \item Given an objective reconstruction frequency index M, one shall make sure that enough coefficients $M_c \geq M$ are considered to evaluate the reconstruction Green function (it is possible to increase the number of coefficients $M_c$ used by changing the precision in the calculation of $\psi_\mathcal{D}$).
    \end{itemize}
\end{remark}


\subsection{Recovering the weights of the knots}

After estimating the $K$ positions of the reconstruction knots only the associated weights remain unknown. It is the exact same situation of \cref{sec:fixed_knots} which is already streamlined. The same FISTA algorithm will lead to the weights $\boldsymbol{\beta}$ within few iterations. Indeed, when compared to the grid reconstruction, the number of weights to estimate has reduced from $N >> K_0$ to only $K$.

\begin{remark}[Choice of number K of knots]
    In the CPGD framework, the reconstruction method needs to be specified the number of knots $K$ to consider. The final solution returned by the algorithm has at most $K$ knots, but in practice it is made of less as the weights estimation step tends to set some knots to $0$. Indeed, this is formulated as a LASSO problem, and consequently returns sparse solutions. A solution to give more flexibility to the whole reconstruction procedure is to set this maximal number of reconstruction knots to its maximal value $K=P=M$. That way, the weighting step has available a larger number of knots and can be more selective in the choice of knots to keep for the final solution.
\end{remark}


\section{Frank-Wolfe Algorithm}

Another way of solving problem (\ref{eq:continuous-pb}) is to adapt an optimization algorithm called \textit{Frank-Wolfe} algorithm \cite{frankwolfe1956} or \textit{Conditional Gradient Method} \cite{COCV_2013__19_1_190_0}). Frank-Wolfe algorithm (FW) is built to efficiently solve a problem written as:
\begin{equation}
    \min_{m\in C} f(m)
    \label{eq:frankwolfe_pb}
\end{equation}
where $C$ is a weakly compact convex set of a Banach space and $f$ is a differentiable convex function. It is very interesting when the Banach structure of the search space does not allow the existence of any gradient method. In our particular case, $C$ will be a subspace of $\mathcal{M}$ but the objective function we used so far does not comply with the requirements of such a framework. However, with a little effort, it is possible to sort out an equivalent formulation that is FW-friendly. Let us detail the algorithm generic formulation in \cref{alg:fw}.

\begin{algorithm}
\caption{Frank-Wolfe Algorithm}
\label{alg:fw}
\begin{algorithmic}[1]
\For {$k=0, \cdots, n$}
    \State Minimize : $c^{[k]} \in \argmin_{c\in C} {f(m^{[k]}) + \mathrm{d}f(m^{[k]}) [c - m^{[k]}]}$.
    \If{$\mathrm{d}f(m^{[k]}) [c^{[k]} - m^{[k]}] = 0$}
    \State $m^{[k]}$ solution of \eqref{eq:frankwolfe_pb}. Stop.
    \Else
    \State Step size search: $\gamma^{[k]} \gets \frac{2}{k+2}$ or $\gamma^{[k]} \in \argmin_{\gamma\in [0, 1]}{f(m^{[k]} + \gamma(c^{[k]}-m^{[k]}))}$
    \State Update: $m^{[k+1]} \gets m^{[k]} + \gamma^{[k]}(c^{[k]}-m^{[k]})$
    \EndIf
\EndFor
\end{algorithmic}
\end{algorithm}


\subsection{Reformulation of the original problem}

The same way as before, consider a source sparse signal as defined in \eqref{eq:source_sparse}. Inspired from Q. Denoyelle's work in \cite{denoyelle2018sliding}, we want to express the measurement operator $\mathbf{\Phi}$ as an integral operator that applies to the stream of periodic Dirac (that we previously denoted as $m_0(t) = \sum_k{\beta_{0,k} \Sha(t-t_{0, k})}$). More precisely, we want to write:
\begin{equation}
    \mathbf{y}_0 = \mathbf{\Phi}(f_0) = \mathbf{\tilde{\Phi}}(m_0)
\end{equation}
thus:
\begin{equation}
    \mathbf{\tilde{\Phi}} = \mathbf{\Phi} \mathcal{D}^{-1}.
\end{equation}
Again, let us use a generic $K$-sparse solution function as \eqref{eq:sol_sparse}. The following holds:
\begin{align*}
    y =& \left\{\mathbf{\Phi} \mathcal{D}^{-1}\right\} (m) \\
    =& \mathbf{\Phi} \mathcal{D}^{-1}\left(\sum_{k=1}^K{\beta_k \Sha(\cdot - t_k)}\right) \\
    =& \mathbf{\Phi}(f) \\
    =& \mathbf{\Phi}\left(\sum_{k=1}^{K} \beta_k \psi_\mathcal{D}(\cdot - t_k)\right) \\
    =& \left(\sum_{k=1}^K \beta_k \psi_\mathcal{D}(\theta_\ell - t_k)\right)_{\ell=1, \dots, L}\\
    =& \left(\int_{\mathbb{T}} \psi_\mathcal{D}(\theta_\ell - t)\ \mathrm{d}m(t) \right)_{\ell=1, \dots, L} \\
    =& \left(\int_{\mathbb{T}} \varphi_\ell(t)\ \mathrm{d}m(t) \right)_{\ell=1, \dots, L}
\end{align*}
with
\begin{equation}
    \forall t \in \mathbb{T}, \qquad \varphi_\ell(t) = \psi_\mathcal{D}(\theta_\ell - t).
\end{equation}
Then the integral operator appears as:
\begin{equation}
    \forall \mu \in\mathcal{M}(\mathbb{T}), \qquad \mathbf{\tilde{\Phi}}(\mu) = \int_\mathbb{T} \boldsymbol{\varphi}(t)\ \mathrm{d}\mu(t)
\end{equation}
and the new measurement operator eventually becomes:
\begin{equation}
    \mathbf{\tilde{\Phi}} \colon \left\{ \begin{aligned}
        \quad \mathcal{M}(\mathbb{T}) \quad & \to \quad \mathbb{R}^L \\
        \mu(t) \quad & \mapsto \quad \left( \int_{\mathbb{T}} \varphi_\ell\ \mathrm{d}\mu \right)_{\ell=1, \dots, L}
    \end{aligned} \right.
\end{equation}
Particularly, if $\mu(t) = \sum_{i=1}^I{\alpha_i \Sha(t-t_i)}$, for any $I\in\mathbb{N}$, $(\alpha_i)_{i=1, \dots, I}\in\mathbb{R}^I$, then:
\begin{equation}
    \mathbf{\tilde{\Phi}}(\mu) = \left( \sum_{i=1}^I{\alpha_i \psi_\mathcal{D}(\theta_\ell-t_i)} \right)_{\ell=1, \dots, L}
\end{equation}
In order to comply with the context presented in \cite{denoyelle2018sliding}, it is necessary to make sure that $\boldsymbol{\varphi} \in \mathcal{C}^{(2)}$ (see \cite[Definition~4]{denoyelle2018sliding}).

\noindent Another few lines of calculations lead to its adjoint operator that will be useful later on:
\begin{equation}
    \mathbf{\tilde{\Phi}^*} \colon \left\{ \begin{aligned}
        \quad \mathbb{R}^L \quad & \to \quad \mathscr{C}_0(\mathbb{T})\\
        \mathbf{h} \quad & \mapsto \quad \sum_{\ell=1}^L {h_\ell\ \psi_\mathcal{D}(\theta_\ell - \cdot)}
    \end{aligned} \right.
\end{equation}
Indeed, we can write :
\begin{align*}
    \langle\mathbf{\tilde{\Phi}}\mu, \mathbf{h}\rangle_{\mathbb{R}^L} &=\ \langle \mu, \mathbf{\tilde{\Phi}^*h}\rangle_{\mathcal{M}, \mathscr{C}_0}\\
    \sum_{\ell=1}^L{ h_\ell \int_\mathbb{T}{ \psi_\mathcal{D}(\theta_\ell - t)}\  \mathrm{d}\mu(t)} &= \int_\mathbb{T}{\left\{\mathbf{\tilde{\Phi}^*h}\right\}(t)\ \mathrm{d}\mu(t)} \\
    \int_\mathbb{T}{\ \sum_{\ell=1}^L{ h_\ell \psi_\mathcal{D}(\theta_\ell - t)}\  \mathrm{d}\mu(t)} &= \int_\mathbb{T}{\left\{\mathbf{\tilde{\Phi}^*h}\right\}(t)\ \mathrm{d}\mu(t)} \\
\end{align*} \\

Remembering that $\mathcal{D}f = m$, problem \eqref{eq:continuous-pb} can be expressed in the equivalent following way:
\begin{equation}
    \min_{m\in\mathcal{M}(\mathbb{T})} \mathrm{P}_\lambda (m) := \norm{\mathbf{y}_0 - \mathbf{\tilde{\Phi}}(m)}_2^2 + \lambda \norm{m}_{TV}.
    \label{eq:fw_pb}
\end{equation}
This can be stated another way around. Thanks to the invertibility of $\mathcal{D}$, the set of solutions of (\ref{eq:fw_pb}) simply is $\mathcal{D}^{-1}$ applied to the set of solutions of (\ref{eq:continuous-pb}):
\begin{equation}
    \mathrm{Sol}\,(\ref{eq:fw_pb}) = \mathcal{D}^{-1}\{\mathrm{Sol}\,(\ref{eq:continuous-pb})\} .
\end{equation}
However, $\mathrm{P}_\lambda$ is not differentiable as it and thus this expression is not well suited for directly applying the Frank-Wolfe algorithm. Problem (\ref{eq:fw_pb}) can be turned into an equivalent differentiable problem by means of an epigraphical lift, as demonstrated in the following lemma.
\begin{lemma}
    The B-Lasso problem 
    \begin{equation*}
        \min_{m\in\mathcal{M}(\mathbb{T})} \mathrm{P}_\lambda (m) := \norm{\mathbf{y}_0 - \mathbf{\tilde{\Phi}}(m)}_2^2 + \lambda \norm{m}_{TV}
    \end{equation*}
    is equivalent to 
    \begin{equation}
    \min_{(s, m)\in C} \mathrm{\widetilde{P}}_\lambda (s, m) := \norm{\mathbf{y}_0 - \mathbf{\tilde{\Phi}}(m)}_2^2 + \lambda s.
    \label{eq:fw_pb_diff}
    \end{equation}
    with $C = \left\{(s,m) \in \mathbb{R}_+ \times \mathcal{M}(\mathbb{T})\ ;\ \norm{m}_{TV}\leq s \leq M = \norm{\mathbf{y}_0}^2_2/\lambda \right\}$.
\end{lemma}

\begin{proof}
Consider a solution $m^\star$ of $\mathrm{P}_\lambda$ on $\mathcal{M}(\mathbb{T})$, then
\begin{equation*}
    \mathrm{P}_\lambda(m^\star) \leq \mathrm{P}_\lambda (0) = \lambda M
\end{equation*}
with $M = \norm{\mathbf{y}_0}^2_2/\lambda$.
Moreover,
\begin{equation*}
    \mathrm{P}_\lambda(m^\star) \leq \mathrm{P}_\lambda (m) \leq \mathrm{\widetilde{P}}_\lambda (s, m)
\end{equation*}
for any $m\in\mathcal{M}$ and $s\geq \norm{m}_1$. The equality is reached when $m=m^\star$ and $s = \norm{m^\star}$ for instance.
\end{proof}
In their paper, Denoyelle et al. also prove that this latter problem \eqref{eq:fw_pb_diff} is well-posed, and $\mathrm{\widetilde{P}}_\lambda$ reaches a minimum on $C$, which concludes justifying the use of the algorithm.


\subsection{Practical implementation of the algorithm}

The Frank-Wolfe algorithm leverages the differentiable property of the objective function to iteratively minimize a first order linear approximation. It keeps track of a running iterate $m^{[k]}$ that best approximates $m$ after $k$ steps.

We apply the algorithm as previously described with the objective function $f$ being $\mathrm{\widetilde{P}}_\lambda$, thus the current iterate is actually an element $(s^{[k]}, m^{[k]}) \in C$. The next iterate $(s^{[k+1]}, m^{[k+1]})$ is consequently calculated as a convex combination of $(s^{[k]}, m^{[k]})$ and the solution $c^{[k]} = (t^{[k]}, u^{[k]})$ of the minimization of the first order approximation of $\mathrm{\widetilde{P}}_\lambda$ given by:
\begin{equation}
    \left(t^{[k]}, u^{[k]}\right) \in \argmin_{c = (t, u) \in C} {\mathrm{\widetilde{P}}_\lambda(t, u) + \mathrm{d}\mathrm{\widetilde{P}}_\lambda(s^{[k]}, m^{[k]}) \left(t-s^{[k]}, u - m^{[k]}\right)}.
    \label{eq:first_order}
\end{equation}
With $\gamma^{[k]}$ the combination weight, the following holds :
\begin{equation}
    m^{[k+1]} = m^{(k]} + \gamma^{[k]} \left(u^{[k]} - m^{[k]}\right) = \left(1-\gamma^{[k]}\right) m^{[k]} + \gamma^{[k]} u^{[k]}
\end{equation}
\begin{equation}
    s^{[k+1]} = s^{(k]} + \gamma^{[k]} \left(t^{[k]} - s^{[k]}\right) = \left(1-\gamma^{[k]}\right) s^{[k]} + \gamma^{[k]} t^{[k}    
\end{equation}


\subsubsection{New spike creation}
It is important to point out that the approximation step (\ref{eq:first_order}) aims at finding the minimum of a linear functional, thus this minimal point is necessary either $(0, 0)$ or located on the border of the set $C$. Extreme points of $C$ are of the form $(M, \pm M\ \Sha(\cdot - t))$ for $t\in\mathbb{T}$. Consequently, what the minimisation step does is basically adding a new periodic Dirac impulse to the approximation iterate, which is then reweigthed.

First, let us detail the calculation of the differential of $\mathrm{\widetilde{P}}_\lambda$. One can see the data fidelity term as the composition between the two functions:
\begin{equation*}
    R \colon \left\{ \begin{aligned}
        \quad \mathbb{R}^L \quad & \to \quad \mathbb{R} \\
        \mathbf{x} \quad & \mapsto \quad \norm{\mathbf{y}_0 - \mathbf{x}}^2_2
    \end{aligned} \right. 
    \quad \mathrm{and} \quad
    \mathbf{\tilde{\Phi}} \colon \left\{ \begin{aligned}
        \qquad \mathcal{M}(\mathbb{T}) \qquad & \to \quad \mathbb{R}^L \\
        \mu(t) \quad & \mapsto \quad \mathbf{\tilde{\Phi}}\mu
    \end{aligned} \right.
\end{equation*}
that we can easily determine the differential of. Indeed, $R$ is a quadratic form and $\mathbf{\tilde{\Phi}}$ is linear, thus is equal to its differential. We obtain, for any $\mathbf{a}\in\mathbb{R}^L$ and $m\in\mathcal{M}\left(\mathbb{T}\right)$:
\begin{equation*}
    \mathrm{d}R(\mathbf{a)} \colon \left\{ \begin{aligned}
        \quad \mathbb{R}^L \quad & \to \quad \mathbb{R} \\
        \mathbf{h} \quad & \mapsto \quad \langle \mathbf{y-a}, \mathbf{h} \rangle_{\mathbb{R}^L}
    \end{aligned} \right. 
    \quad \mathrm{and} \quad
    \mathrm{d}\mathbf{\tilde{\Phi}}(\mu) \colon \left\{ \begin{aligned}
        \qquad \mathcal{M}(\mathbb{T}) \qquad & \to \quad \mathbb{R}^L \\
        \mu'(t) \quad & \mapsto \quad \mathbf{\tilde{\Phi}}\mu'
    \end{aligned} \right.
\end{equation*}


Given $ R \circ \mathbf{\tilde{\Phi}} : \mathcal{M}(\mathbb{T}) \  \to \ \mathbb{R}$, and for any $\mu\in\mathcal{M}\left(\mathbb{T}\right)$, $\mathrm{d} \left( R \circ \mathbf{\tilde{\Phi}}\right) (\mu) : \mathcal{M}(\mathbb{T}) \  \to \ \mathbb{R}$, we obtain:
\begin{equation*}
    \mathrm{d} \left( R \circ \mathbf{\tilde{\Phi}}\right) (\mu) \cdot \mu' = \left( \mathrm{d}R \left( \mathbf{\tilde{\Phi}}m \right) \circ \mathrm{d}\mathbf{\tilde{\Phi}}(\mu)\right) \cdot \mu' = \langle \mathbf{y}_0 - \mathbf{\tilde{\Phi}}\mu, \mathbf{\tilde{\Phi}}(\mu') \rangle_{\mathbb{R}^L}
\end{equation*}
that can be nicely rewritten as:
\begin{equation}
    \mathrm{d} \left( R \circ \mathbf{\tilde{\Phi}}\right) (\mu) \cdot \mu' = \langle \mathbf{\tilde{\Phi}^*} \left( \mathbf{y}_0 - \mathbf{\tilde{\Phi}}\mu \right), \mu' \rangle_{\mathcal{C}_0 \left( \mathbb{T} \right), \mathcal{M}\left(\mathbb{T}\right)} = \int_{\mathbb{T}} \mathbf{\tilde{\Phi}^*} \left( \mathbf{y}_0 - \mathbf{\tilde{\Phi}}\mu \right) \mathrm{d}\mu'
\end{equation}

Finally, the closed form expression is:
\begin{equation}
    \mathrm{d} \mathrm{\widetilde{P}}_\lambda(s, \mu) :\quad (s', \mu') \longmapsto \int_\mathbb{T} 2\mathbf{\tilde{\Phi}^*}\left(\mathbf{\tilde{\Phi}}\mu - \mathbf{y}_0\right) \mathrm{d}\mu' + \lambda s'.
\end{equation}
Then, finding a minimizer of (\ref{eq:first_order}) having the shape of a single peak periodic Dirac amounts to finding a point
\begin{equation}
    t^{[k]} \in \argmin_{t\in\mathbb{T}} \left(\pm \frac{2}{\lambda} \left( \mathbf{\tilde{\Phi}^*} (\mathbf{\tilde{\Phi}}m^{[k]} - \mathbf{y}_0) \right) (t) + 1 \right)\lambda M
\end{equation}
or equivalently
\begin{equation}
    t^{[k]} \in \argmax_{t\in\mathbb{T}} \left(\left|\eta^{[k]}(t)\right|-1\right)
    \label{eq:max_certificate}
\end{equation}
with
\begin{equation}
    \eta^{[k]}(t) := \frac{2}{\lambda} \left( \mathbf{\tilde{\Phi}^*} (\mathbf{y}_0 - \mathbf{\tilde{\Phi}}m^{[k]} ) \right) (t)
\end{equation}
This function $\eta^{[k]}$ is called empirical dual certificate. Any given iterate $m^{[k]}$ being expressed as a sum of (at most, depends on some $0$-weighted Dirac) $k$ periodic Dirac, \textit{i.e.} $m^{[k]} = \sum_{i=1}^k \alpha_i^{[k]} \Sha( \cdot - t_i^{[k]})$, the function $\eta^{[k]}$ admits the closed form expression:
\begin{equation}
    \eta^{[k]}(t) = \frac{2}{\lambda} \sum_{\ell=1}^L\left( y_\ell - \sum_{j=1}^L \sum_{i=1}^k \alpha_i^{[k]}  \psi_\mathcal{D}(\theta_j - t_i^{[k]}) \right) \psi_\mathcal{D}(\theta_\ell - t).
\end{equation}

In practice, it is not sufficient to know the position of the next periodic Dirac impulse, it is also necessary to know if it is positively or negatively weighted. 
The first naive idea basically solves two different minimization problems, each one related to one of the two situations. The one leading to the lower value of the objective function determines the sign of the new impulse to add to the current iterate. Positive impulse leads to the problem:
\begin{equation}
    \argmin_{t\in\mathbb{T}} -\eta^{[k]}(t)
\end{equation}
By linearity of $\mathbf{\tilde{\Phi}}^*$, negatively weighted impulses lead to the opposite problem.

However, it is actually possible to complete this step by solving only one minimization problem. Let us consider any intermediate step index $k$ of the algorithm. We write :
\begin{equation}
    L^{[k]} = \mathrm{d} \mathrm{\widetilde{P}}_\lambda(s^{[k]}, m^{[k]})\ \colon \qquad \left(\mathbb{R}, \mathcal{M}\left(\mathbb{T}\right)\right) \quad\to\quad \mathbb{R}
\end{equation}
Let $\varepsilon \in \{-1, +1\}$, then, for any $t\in\mathbb{T}$, we have:
\begin{align}
    L^{[k]} (M, \varepsilon M \Sha_t)\ =&\ 2\varepsilon M \left\{ \mathbf{\tilde{\Phi}^*}\left(\mathbf{\tilde{\Phi}}m^{[k]} - \mathbf{y}_0\right) \right\} (t) + \lambda M \nonumber \\
    =&\ -\varepsilon M \lambda \eta^{[k]}(t) + \lambda M \nonumber \\
    =&\ \lambda M \left(1 - \varepsilon \eta^{[k]}(t) \right)
\end{align}
However, when $t = t_{k+1}$ is the solution of (\ref{eq:max_certificate}), there exist $\varepsilon^{[k+1]} \in \{-1, +1\}$ such that the pair $\left(M, \varepsilon^{[k+]} M \Sha_{t_{k+1}} \right)$ minimizes $L^{[k]}$ over $C$. $\varepsilon^{[k+1]}$ simply is the sign of the next periodic Dirac impulse. In particular, being the minimum:
\begin{equation*}
    L^{[k]}\left(M, \varepsilon^{[k+1]} M \Sha_{t_{k+1}} \right) \leq L^{[k]}(0, 0) = 0
\end{equation*}
Thus:
\begin{equation}
    \varepsilon^{[k+1]} \eta^{[k]}(t_{k+1}) \geq 1 > 0
\end{equation}
Which leads to :
\begin{equation}
    \varepsilon^{[k+1]} = \mathrm{sgn}(\eta^{[k]}(t_{k+1}))
\end{equation}
Eventually, the algorithm first numerically solves (\ref{eq:max_certificate}) to estimate the position of the next impulse, and then evaluates the empirical dual certificate at this latter position to determine the sign of the impulse to append to the current iterate.

\subsubsection{Combination weight}

The weight $\gamma^{[k]}$ of the convex combination can be expressed as  decreasing with $k$ in such a manner that the first iterates have more significance than the last ones in the final expression of the solution. For instance, one can choose:
\begin{equation*}
    \gamma^{[k]} = \frac{2}{k + 2}
\end{equation*}
The convergence is actually quite long and with very low accuracy.

Another more efficient way is to find $\gamma^{[k]}$ that minimizes the expression
\begin{equation}
    \gamma^{[k]} \in \argmin_{\gamma \in [0, 1]} \mathrm{\widetilde{P}}_\lambda\left(\norm{m^{[k]}}_{TV} + \gamma (M - \norm{m^{[k]}}_{TV}), m^{[k]} + \gamma (u^{[k]} - m^{[k]})\right)
    \label{eq:gamma_pb}
\end{equation}
This method conduces to better results in less iterations as the iterates are tailored accordingly to the problem, not with general coarse formula.

Additionally, we can show that it is possible to compute the exact value of $\gamma^{[k]}$ at each step. First let us emphasize a simple property that was not explicitly stated in \cite{denoyelle2018sliding}. Let us call the first variable $s$ of (\ref{eq:fw_pb_diff}) the regularization variable, and denote by $s^{[k]}$ the iterates when applying the FW algorithm. Then, the following holds:
\begin{lemma}
    For any iteration index $k\in\mathbb{N}$,
    \begin{equation}
        s^{[k]} = \norm{m^{[k]}}_{TV}
        \label{eq:regularization_norm}
    \end{equation}
\end{lemma}
\begin{proof}
If we consider two finite weighted sums of Dirac combs $g_1$ and $g_2$, with distinct knots in $\mathbb{T}$ and no common knot between the two, then $\norm{g_1 + g_2}_{TV} = \norm{g_1}_{TV} + \norm{g_2}_{TV}$. Then, for any $\gamma\in\left[0, 1\right]$, at any iteration step $k$, we can write
$$\norm{m^{[k]} + \gamma \left(u^{[k]} - m^{[k]}\right)}_{TV} = \left( 1-\gamma\right)\norm{m^{[k]}}_{TV} + \gamma \norm{u^{[k]}}_{TV}$$
(as long as the new Dirac impulse position does not belong to the previous set of positions, which is supposed to be guaranteed by the FW theory). Starting the iterates of FW with $(0, 0)$ at iteration $k=0$ (also true at iteration $k=1$ with $\left(\gamma^{[1]}M,\pm \gamma^{[1]}M\Sha_{x_1}\right)$), equation \eqref{eq:regularization_norm} is verified. Equation \eqref{eq:regularization_norm} being true at any given iteration $k$, then : $$s^{[k+1]} = \left( 1-\gamma^{[k]}\right)s^{[k]} + \gamma^{[k]}M = \left( 1-\gamma^{[k]}\right)\norm{m^{[k]}}_{TV} + \gamma^{[k]} \norm{u^{[k]}}_{TV} = \norm{m^{[k+1]}}_{TV},$$ it is still true at iteration $k+1$, which proves the property by induction.
\end{proof}

This result allows to explicitly compute the solution.
\begin{proposition}
The unique solution to \eqref{eq:gamma_pb} is given by
\begin{equation}
    \gamma^{[k]} = \frac{-2 \left[ \left(\mathbf{\tilde{\Phi}}m^{[k]} - \mathbf{y}_0 \right)^T\left(\mathbf{\tilde{\Phi}}u^{[k]}+\mathbf{y}_0 \right) + \norm{\mathbf{y}_0}_2^2 - \norm{\mathbf{\tilde{\Phi}}m^{[k]}}^2_2\right] -\lambda \left[M - \norm{m^{[k]}}_{TV}\right]}{2 \norm{\mathbf{\tilde{\Phi}}m^{[k]} - \mathbf{\tilde{\Phi}}u^{[k]}}^2_2}.
\end{equation}
\end{proposition}

\begin{proof}
Indeed, one can rewrite the expression of the objective function of \eqref{eq:gamma_pb} as follows:
\begin{align*}
    &\mathrm{\widetilde{P}}_\lambda\left(s^{[k]} + \gamma (M - s^{[k]}),\ m^{[k]} + \gamma (u^{[k]} - m^{[k]})\right) \\
    &=\mathrm{\widetilde{P}}_\lambda\left(\norm{m^{[k]}}_{TV} + \gamma (M - \norm{m^{[k]}}_{TV}),\ m^{[k]} + \gamma (u^{[k]} - m^{[k]})\right) \\
    &= \norm{ \mathbf{y}_0 - \gamma \mathbf{\tilde{\Phi}}u^{[k]} - (1-\gamma)\mathbf{\tilde{\Phi}}m^{[k]} }_2^2 + \lambda (\gamma M + (1-\gamma) \norm{m^{[k]}}_{TV}) \\
    &= \gamma^2 \left[ \norm{\mathbf{\tilde{\Phi}}u^{[k]}}_2^2 + \norm{\mathbf{\tilde{\Phi}}m^{[k]}}_2^2 - 2 \left(\mathbf{\tilde{\Phi}}u^{[k]}\right)^T \mathbf{\tilde{\Phi}}m^{[k]} \right] \\
    & \quad +2 \gamma \left[ - \norm{\mathbf{\tilde{\Phi}}m^{[k]}}_2^2 - \mathbf{y}_0^T \mathbf{\tilde{\Phi}}u^{[k]} + \mathbf{y}_0^T \mathbf{\tilde{\Phi}}m^{[k]} + \left(\mathbf{\tilde{\Phi}}u^{[k]}\right)^T \mathbf{\tilde{\Phi}}m^{[k]} \right] + \gamma\lambda \left[M - \norm{m^{[k]}}_{TV}\right]\\
    & \quad\quad + \left[ \norm{\mathbf{y}_0}^2_2 + \norm{\mathbf{\tilde{\Phi}}m^{[k]}}_2^2 - 2\mathbf{y}_0^T \mathbf{\tilde{\Phi}}m^{[k]} \right] + \norm{m^{[k]}}_{TV} \\
    &= \gamma^2 \norm{\mathbf{\tilde{\Phi}}m^{[k]} - \mathbf{\tilde{\Phi}}u^{[k]}}^2_2 \\
    &\quad +2\gamma \left[ \left(\mathbf{\tilde{\Phi}}m^{[k]} - \mathbf{y}_0 \right)^T\left(\mathbf{\tilde{\Phi}}u^{[k]}+\mathbf{y}_0 \right) + \norm{\mathbf{y}_0}_2^2 - \norm{\mathbf{\tilde{\Phi}}m^{[k]}}_2^2\right] + \gamma \lambda\left[M - \norm{m^{[k]}}_{TV}\right]\\
    & \quad \quad +\norm{\mathbf{\tilde{\Phi}}m^{[k]} - \mathbf{y}_0}_2^2 + \lambda\norm{m^{[k]}}_{TV}
\end{align*}
The expression is a convex second degree polynomial, the minimum of which can be computed thanks to the coefficients. It leads to the aforementioned expression.
\end{proof}

\subsubsection{Stopping criterion}
As stated earlier, the new spike creation implies solving problem (\ref{eq:max_certificate}) in order to estimate the position of the new Dirac impulse. Interestingly enough, the actual value of the maximum of the empirical dual certificate $\max_{t\in\mathbb{T}}{|\eta^{[k]}(t)|} = \norm{\eta^{[k]}}_\infty$ is very useful and so is also regarded. Denoyelle (still in \cite{denoyelle2018sliding}, remark 5) has proved that this value is related to the natural stopping criterion
\begin{equation}
    \left(s^{[k]}, m^{[k]}\right) \in \argmin_{c\in C} \left\{\mathrm{d} \mathrm{\widetilde{P}}_\lambda \left(s^{[k]}, m^{[k]}\right) \left[c\right] \right\}
    \label{eq:stopping_crit}
\end{equation}
in the way that, if $\left(s^{[k]}, m^{[k]}\right)$ satisfies (\ref{eq:stopping_crit}), then the following equation holds:
\begin{equation}
    \norm{\eta^{[k]}}_\infty = 1
\end{equation}
Thus, after calculating the new Dirac impulse position, verifying the stopping criterion becomes a trivial step. Given an arbitrary satisfying precision tolerance $\varepsilon$, the algorithm simply verifies if the lastly calculated value of the empirical certificate function is $\varepsilon$-close from $1$.

In practice, when the minimization step is successfully solved, the successive extremum values of the certificate are monotonously decreasing, although two behaviors have been observed. Either the sequence decays fast and distinctly reaches a value lower than $1$, either it has an asymptotic and slow convergence towards $1$ by upper value. In the first situation, the couple $(0, 0)$ is actually the solution of \eqref{eq:first_order}. Considering this couple as the new spikes and plugging it in the convex combination of solution has the effect of decreasing the current weights of the iterate and eventually increasing the supremum norm of the empirical certificate. It usually becomes greater than 1 and the convergence continues until it reaches the stopping criterion. In the second situation, the sopping criterion is naturally reached after a finite number of iterations.

\subsubsection{Choice of the regularization parameter}

Let us recall the expression of the empirical dual certificate:
\begin{equation*}
    \eta^{[k]}(t) = \frac{2}{\lambda} \left( \mathbf{\tilde{\Phi}^*} (\mathbf{y}_0 - \mathbf{\tilde{\Phi}}m^{[k]} ) \right) (t)
\end{equation*}
Then, at step $k=0$, $m^{[0]} = 0$ such that:
\begin{equation}
    \eta^{[0]}(t) = \frac{2}{\lambda} \mathbf{\tilde{\Phi}^*} (\mathbf{y}_0) (t)
\end{equation}
Consequently, in the same way as the equispaced knots strategy, the value of the regularization parameter $\lambda$ is also limited by an extreme value of
\begin{equation}
    \lambda_{max} := 2 \max_{t\in\mathbb{T}} | \mathbf{\tilde{\Phi}^*} (\mathbf{y}_0) (t) | = 2\norm{\mathbf{\tilde{\Phi}^*} (\mathbf{y}_0)}_\infty.
\end{equation}
Any value of $\lambda$ larger than $\lambda_{max}$ would force the algorithm to stop at the first iteration. In practice, we also set $\lambda = \sigma \norm{\mathbf{\tilde{\Phi}^*} (\mathbf{y}_0)}_\infty$ with $\sigma \approx 0.1$, which conduces to the best empirical results.

To put everything in a nutshell, the Frank-Wolfe algorithm approximates the sum of periodic Dirac $m(t) = \sum_{k=1}^K \beta_k \Sha(t-t_k)$ using an iterative method. Each iteration of the loop adds a new weighted periodic Dirac impulse to the current approximate, the weights being decreased along the iterations to ensure that the first retrieved Dirac impulses are given more significance in the reconstruction.


\subsection{All-knots simultaneous reweighting strategy}

An assessment we made during our experimentation is that the regular FW strategy is greedy. Indeed, each iteration step exhibits a new Dirac impulse to add to the current iterate by maximizing a given cost function over the problem (namely the empirical dual certificate). This greedy approach might be successful, and even efficient, when the problem is well presented (noiseless, with regular measurement points). However, it might also lead to very unsatisfactory approximations. Indeed, this greedy approach is a \enquote{forward} strategy, that is only able to produce more complexity in the solution by increasing the size of the reconstruction support. It is not able to perform retro actions on the previously evaluated weights or positions.

One way to circumvent this issue comes from the analysis of the procedure. The regular FW in practice retrieves the innovations by first finding the position, and then, separately, estimating the weight. While the position estimation step is really subtle, from deep dual analysis, the weight estimation appears quite simple. The convex combination relies on the fact that all the previous weights have been correctly estimated, and gives really small margin for the new impulse to be significant with respect to the others. Then, the same way it is done with the sliding version \cite[Algorithm~2]{denoyelle2018sliding} or in \cite{courbot:hal-02940848}, a solution is to re-estimate the weights of all the knots at the same time, once a new candidate knot position is obtained. More precisely, at any step $k$, the current iterate is $m^{[k]}$ and the new spike to be added is at position $u^{[k]} \in \mathbb{T}$ from \eqref{eq:first_order}. We solve the LASSO problem associated to the fixed knots $\mathrm{Supp}\left(m^{[k]} \cup u^{[k]}\right)$ to simultaneously estimate all the weights $\boldsymbol{\beta}^{[k]}$.

This strategy is very interesting, as it allows to \emph{a posteriori} cancel any previously placed spike, thanks to the sparsity promoting behavior of the LASSO. It gives the algorithm the right to make mistakes and correct them pretty efficiently. In term of evolution of the value of the objective function $\mathrm{\tilde{P}}$ along the iterations, it guarantees a decreasing behavior, as the same problem is solved step after step, with a larger support of recovery allowed (the set of possible reconstruction knots is increasing). To go even faster and make use of the largest information possible, the FISTA algorithm is initialized with the regular convex combination weighted iterate. It leads to significant gains, both in term of accuracy and duration of computation. The both strategies are really easy to compare as they share the same stopping criterion.

A consequent part of this internship has been focused on designing this reweighting strategy and analyzing its results. It leaded to a fine understanding of the regular FW mechanics as well as its weaknesses.

%% file: Chapters/chapter3.tex

\chapter{Practical Experiments} 

\label{chap:experiences}

One of the explicit goals of this internship was to obtain a working algorithmic framework that is efficiently able to run the three presented algorithms. Thanks to the rigorous analysis of their respective behavior, this objective has been fulfilled and we now dispose of a fresh implementation of each of these algorithms. In this chapter we thus present the common experimental setting that we have proposed for experimentation and show some examples of reconstruction.


\section{Experimental Setting}

Our experimental framework generates measurement data $\mathbf{y}_0$ from a synthetic input signal $f_0 = \mathcal{D}^{-1}m_0$. $m_0$ is a random sparse signal with $K_0$ knots and random associated weights. Then $\mathbf{y}_0$ is fed to the algorithms as the measurement data. The response from the algorithms is given as  list of $K$ innovations that characterizes the sparse reconstructed signal that approximates $f_0$ by solving problem \eqref{eq:continuous-pb}.

\subsubsection{Practical statement of the problem}

Let us dive into details on the input parameters involved.

\begin{itemize}
    \item \textbf{The sampled spline} $f_0$ is built from a differential operator $\mathcal{D}_0$ that needs to be specified as an input. In our experiment, we set its number of knots $K_0 = 4$, and let its innovation be taken randomly (uniform distribution for the knots over $\mathbb{T}$ and normal distribution for the weights). The number $L$ of measurements is set as $8 K_0 + 1$ (so that $M$ from CPGD is $M=4K_0$).
    \item \textbf{The reconstruction spline} $f$ also needs to be built from an operator $\mathcal{D}$. Actually, there is no reason for $\mathcal{D}$ to be equal to $\mathcal{D}_0$, as we suppose there is no prior knowledge on the input data. One might guess the type of operator $\mathcal{D}_0$ involved from the input signal, but the specific operator's parameters are impossible to estimate beforehand. In our experiences, we nevertheless let $\mathcal{D} = \mathcal{D}_0$, but it is important to keep in mind that they might differ.
    \item \textbf{The noise level} needs to be specified as a \emph{peak signal-to-noise ratio} (PSNR) in decibel. It is involved in the samples generation as $\mathbf{y}_0 = \mathbf{\Phi}(f_0) + \mathbf{w}$ with $\mathbf{w}$ drawn according to $\mathcal{N}(0, \omega^2)$ with
    $$\omega = \max_{\ell=1, \dots, L}{|\mathbf{\Phi}(f_0)|} \times \exp{\left(-\frac{PSNR}{10}\right)}.$$
    We set the PSNR to $20$dB for most of the experiences, but this value is often changed to study noise influence.
    \item \textbf{The regularization parameter} is chosen as a fraction of a maximal value $\lambda = \sigma\lambda_{max}$ previously calculated. This parameter needs to increase with the value of the noise, but experiments have shown that $\sigma \approx 0.1$ is usually a good rule of thumb for moderate noise levels of about $20$dB.
\end{itemize}

\subsubsection{Stopping criterion}

As we have already explained, the stopping criterion differs with the algorithm. However, the algorithm FISTA is used at many occurrences, either as the main algorithm (equispaced grid strategy) or as a step to estimate the weights (for CPGD and FW with reweighting). For all of these different cases, the same relative improvement tolerance $\varepsilon$ is set to be equal to $10^{-4}$.

For CPGD, the stopping criterion for the estimation of the truncated Fourier series coefficients is also a relative improvement precision (see algorithm \ref{alg:cpgd}). It is also set to $\varepsilon = 10^{-4}$.

As the FW algorithm is concerned, the stopping criterion is slightly different and hard to relate to the other ones. Indeed, it is based on the maximum absolute value reached by the empirical dual certificate. If this value is $\nu$-close from $1$ the algorithm stops. Actually, when $\nu$ is lower than $10^{-2}$, the regular FW takes significantly longer to converge (when it does). Consequentlly, we do not try to go lower and set $\nu = 10^{-2}$ as it already provides satisfactory results.

\subsubsection{Metrics}

Until now, the metrics have not been much into consideration. The question of a relevant metric really is task-dependent. As we are not concerned yet with any application, we made the choice of considering a bunch of different indicators that we think are relevant for the reconstruction. These indicators include:
\begin{itemize}
    \item the computation time,
    \item the final objective value of the optimization problem (we evaluate \eqref{eq:continuous-pb} with the returned solution),
    \item the final sparsity index of the solution,
    \item the quadratic error between the source spline and the reconstructed one (expressed as the root of the relative squared error), which is an \emph{a priori} quantity as it needs the full input signal to be computed,
    \item the quadratic error between the measurement sampled and the expected value of the reconstructed spline on the measurement points, which can be computed \emph{a posteriori} as it only needs the samples,
    \item the fact that the algorithm has explicitly converged or not.
\end{itemize}

The final value of the objective function has come of great help when comparing the different method. Indeed, it is a common quantity that is agnostic of the reconstruction algorithm, and thus allows to measure to what extent each method performs with respect to the minimization task.


\section{Explicit Results}

\subsection{Grid-based solutions}

\begin{figure}[H]
\makebox[\textwidth][c]{
\begin{subfigure}{0.6\textwidth}
    \centering
    \includegraphics[width=0.95\linewidth,center]{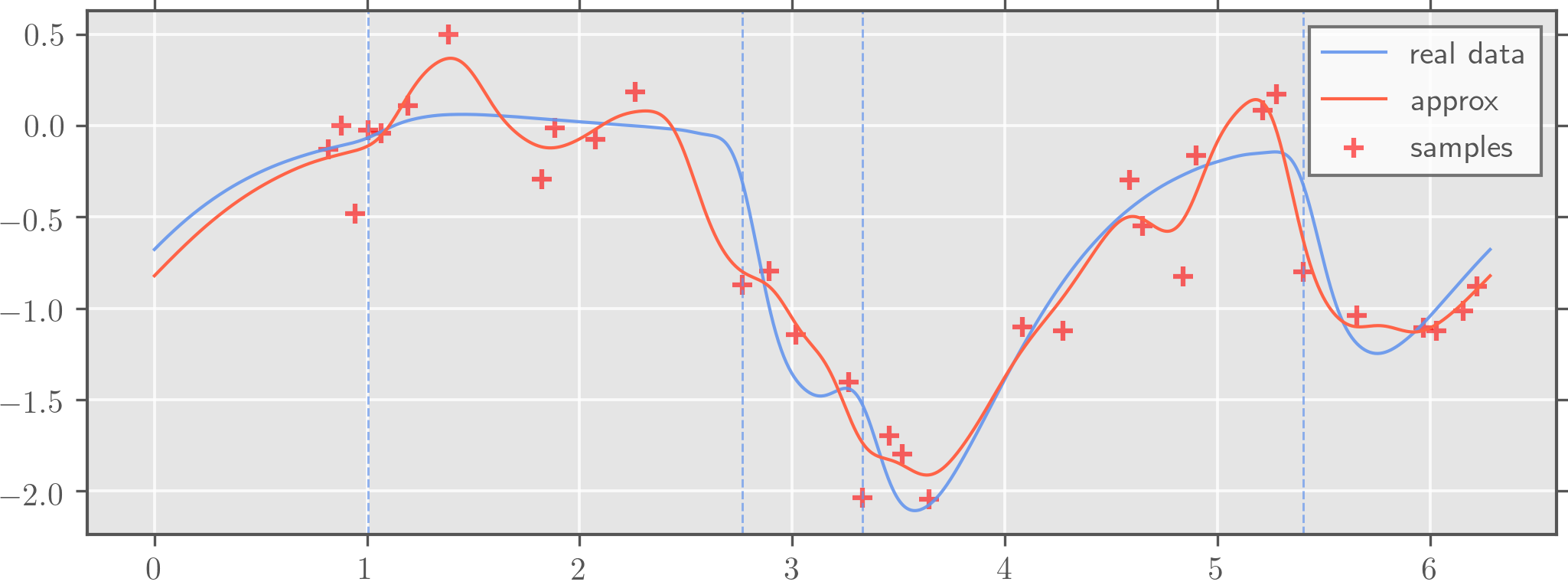}
    \caption{$\sigma$ = 0.01}
\end{subfigure}
\begin{subfigure}{0.6\textwidth}
    \centering
    \includegraphics[width=0.95\linewidth, center]{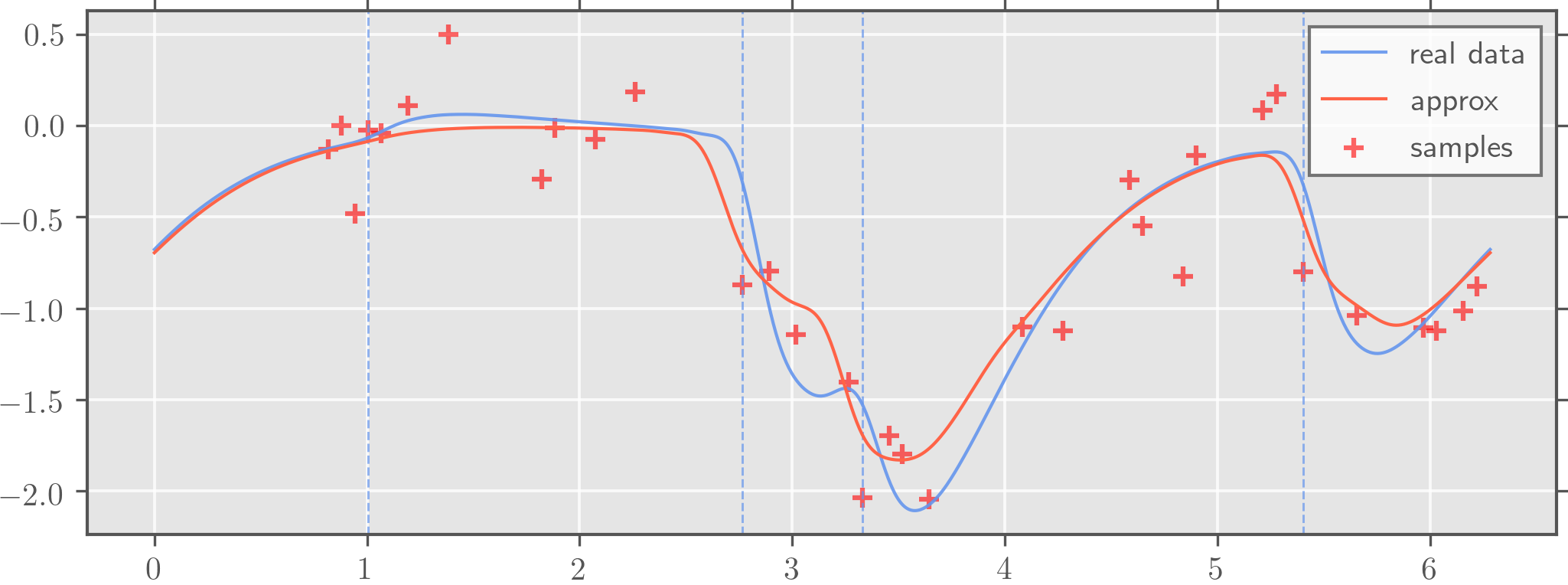}
    \caption{$\sigma$ = 0.1}
\end{subfigure}
}

\makebox[\textwidth][c]{
\begin{subfigure}{0.6\textwidth}
    \centering
    \includegraphics[width=0.95\textwidth]{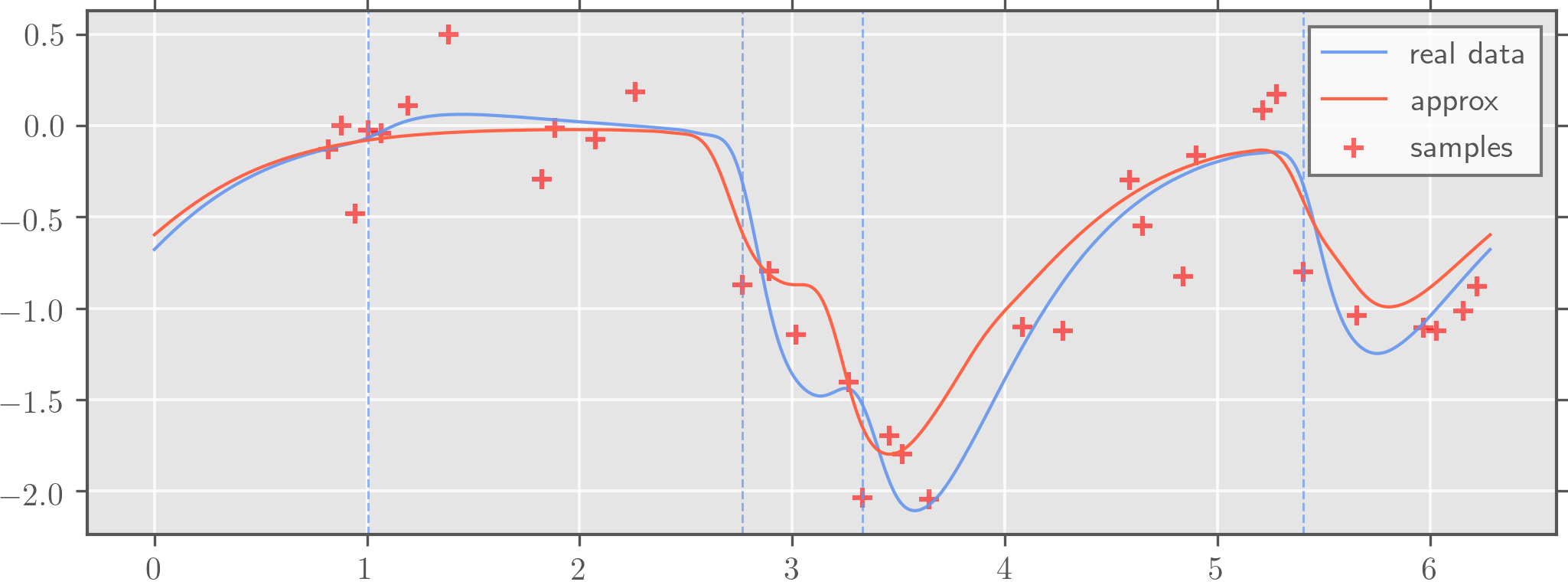}
    \caption{$\sigma$ = 0.2}
\end{subfigure}
\begin{subfigure}{0.6\textwidth}
    \centering
    \includegraphics[width=0.95\textwidth]{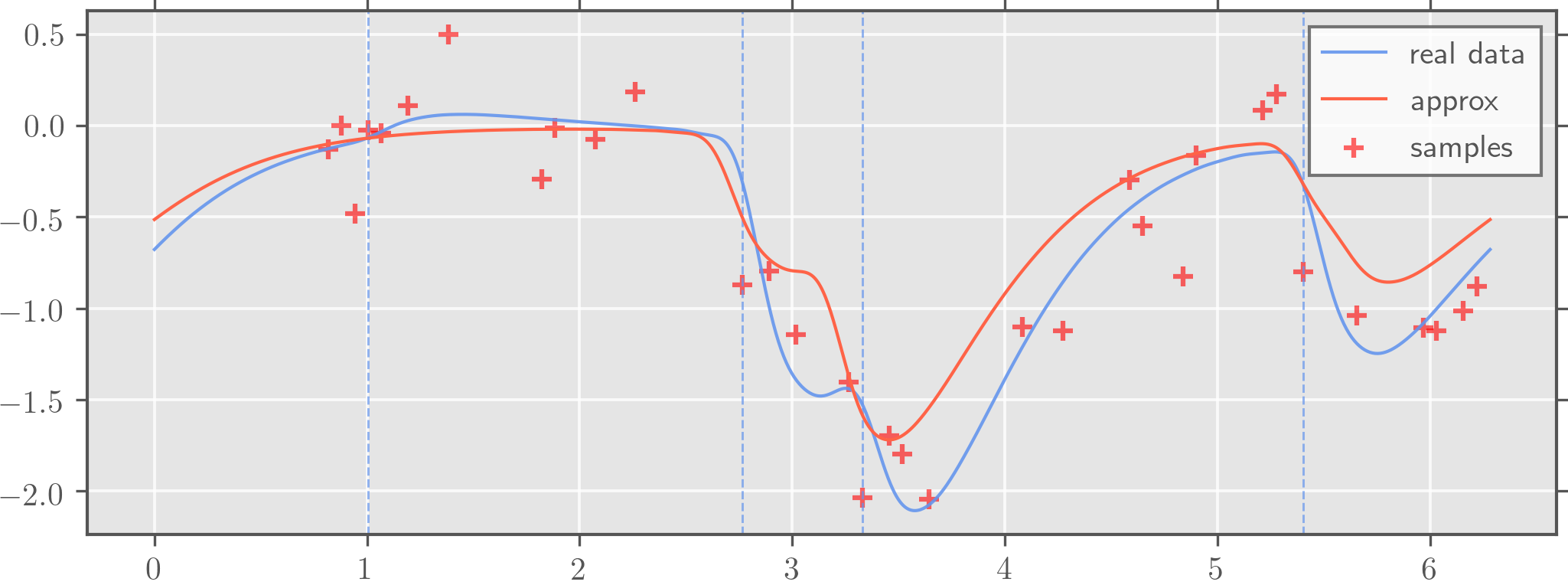}
    \caption{$\sigma$ = 0.3}
\end{subfigure}
}
\caption{Grid-based reconstruction with different values of $\sigma$}
\label{fig:grid_reco}
\end{figure}

Let us set the size of the reconstruction grid to $N=300$ and the penalization parameter $\lambda = \sigma \lambda_{max}$. The noise level is set to $20$ dB and the operator are exponential $\mathcal{D} = \mathcal{D}_0 = \left(\mathrm{D} + \alpha \mathrm{Id}\right)^\gamma $ with $\gamma=2$ and $\alpha=3$. We present different reconstructions with various values of $\sigma$ in \cref{fig:grid_reco} and \cref{tab:grid_reco}.

\begin{table}[H]
    \centering
    \begin{filecontents*}{Figures/reconstruction_grid/seed9970_psnr20.csv}
factors,iterations,duration,converged,objective_fun,rrse_splines,rrse_samples
0.01,2000,0.19648361206054688,False,1.3931745035488845,0.2798293959541258,0.21927235250588065
0.1,2000,0.22111201286315918,False,5.266016298356066,0.22981781116512362,0.3339665338358812
0.2,2000,0.19667792320251465,False,8.616309185678057,0.31827363882711557,0.39360570681195517
0.3,2000,0.18811678886413574,False,11.587176124386204,0.40003255114573516,0.46291663826728974
\end{filecontents*}
\csvreader[TableStyle]
               {Figures/reconstruction_grid/seed9970_psnr20.csv}{1=\factors,2=\iterations,3=\duration,4=\converged,5=\objective,6=\rrse,7=\rrsebis}%
      {\factors & \iterations & \num[round-precision=3]{\duration} & \converged & \num{\objective} & \num{\rrse} & \num{\rrsebis}}%
    \caption{Metrics for grid-based reconstruction}
    \label{tab:grid_reco}
\end{table}

As expected, when $\sigma$ tends to $0$, the reconstruction behaves as an interpolation framework, where the solution tries to go through every sampled point. Obviously, the error between the measurements and the predicted values of the spline at the same positions (referred to as \enquote{samples rrse} in the table) tends to be smaller in this case. However, it is not correlated with the error between the source spline and the approximation one (confer to \enquote{splines rrse} in the table).

Note that none of the cases have converged within $2000$ iterations. It is due to the conditionning number of the measurement matrix $\mathbf{H}$ being quite large, as the matrix is full and of size $N \times N$. Nevertheless, the non-convergence is not a real concern, as the improvement is already quite low and the solution does not vary in the following iterations.

\subsection{CPGD}

\begin{figure}[H]
\makebox[\textwidth][c]{
\begin{subfigure}{0.6\textwidth}
    \centering
    \includegraphics[width=0.95\linewidth,center]{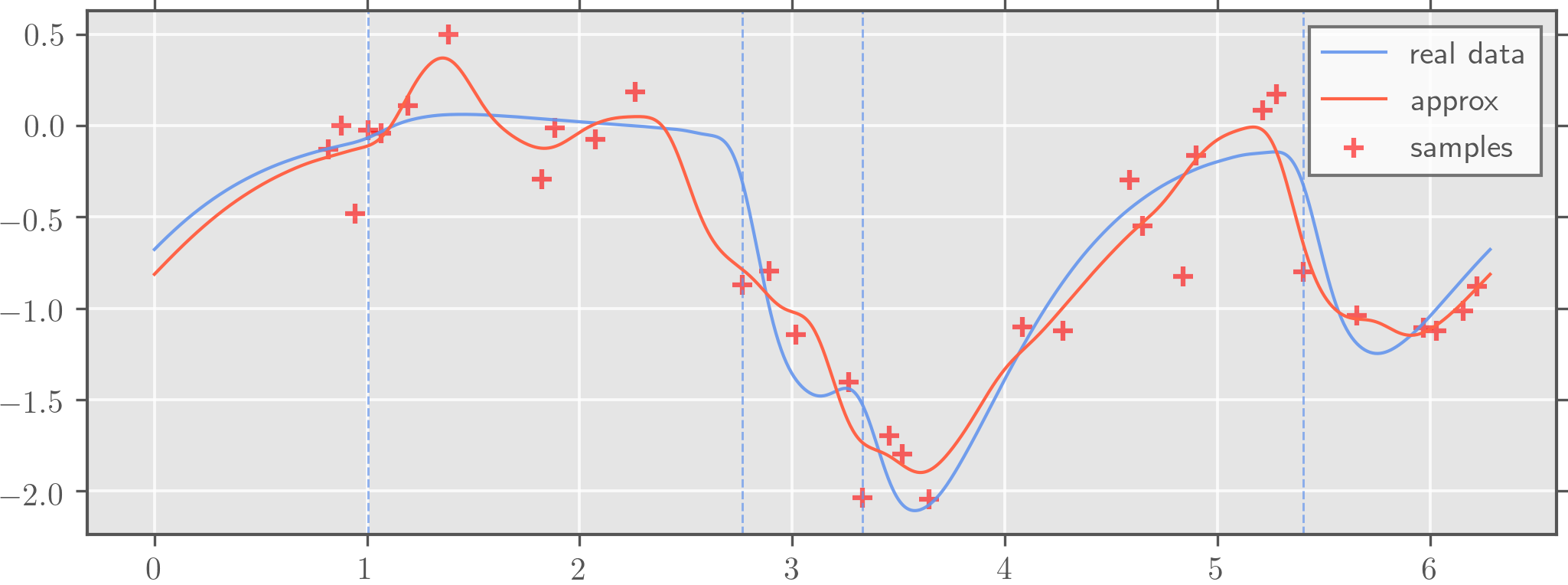}
    \caption{$\sigma$ = 0.01}
\end{subfigure}
\begin{subfigure}{0.6\textwidth}
    \centering
    \includegraphics[width=0.95\linewidth, center]{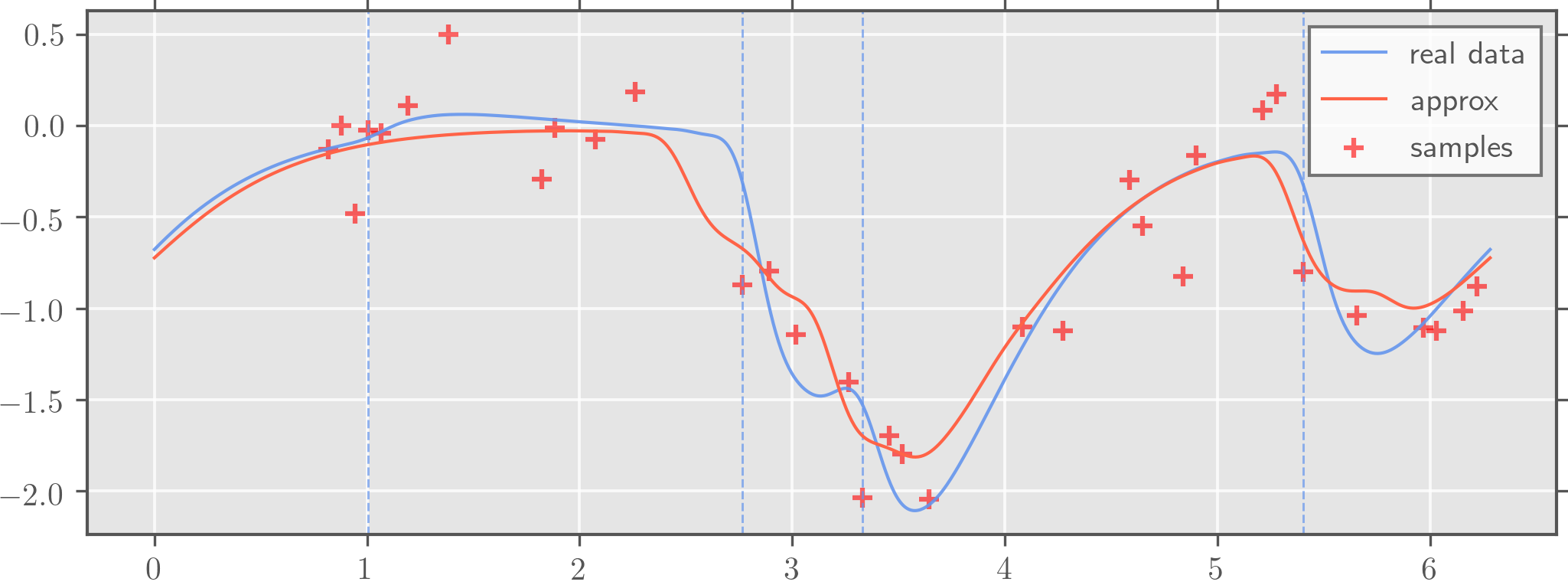}
    \caption{$\sigma$ = 0.1}
\end{subfigure}
}

\makebox[\textwidth][c]{
\begin{subfigure}{0.6\textwidth}
    \centering
    \includegraphics[width=0.95\textwidth]{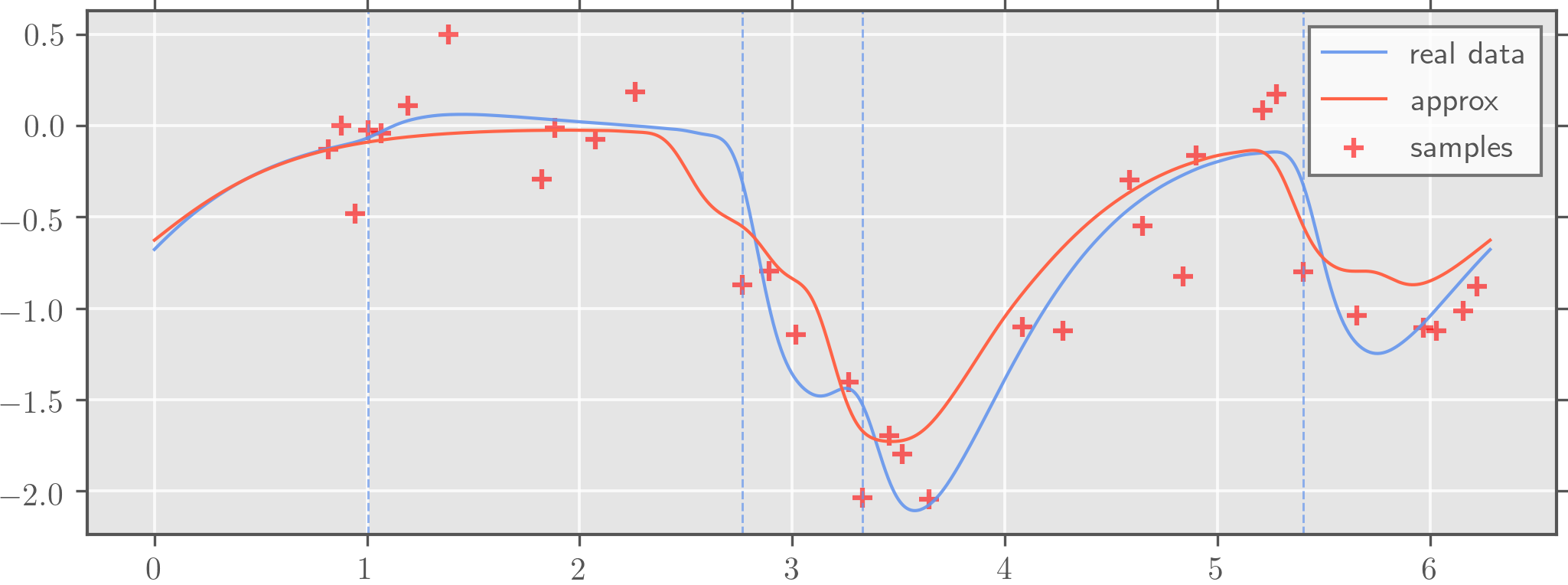}
    \caption{$\sigma$ = 0.2}
\end{subfigure}
\begin{subfigure}{0.6\textwidth}
    \centering
    \includegraphics[width=0.95\textwidth]{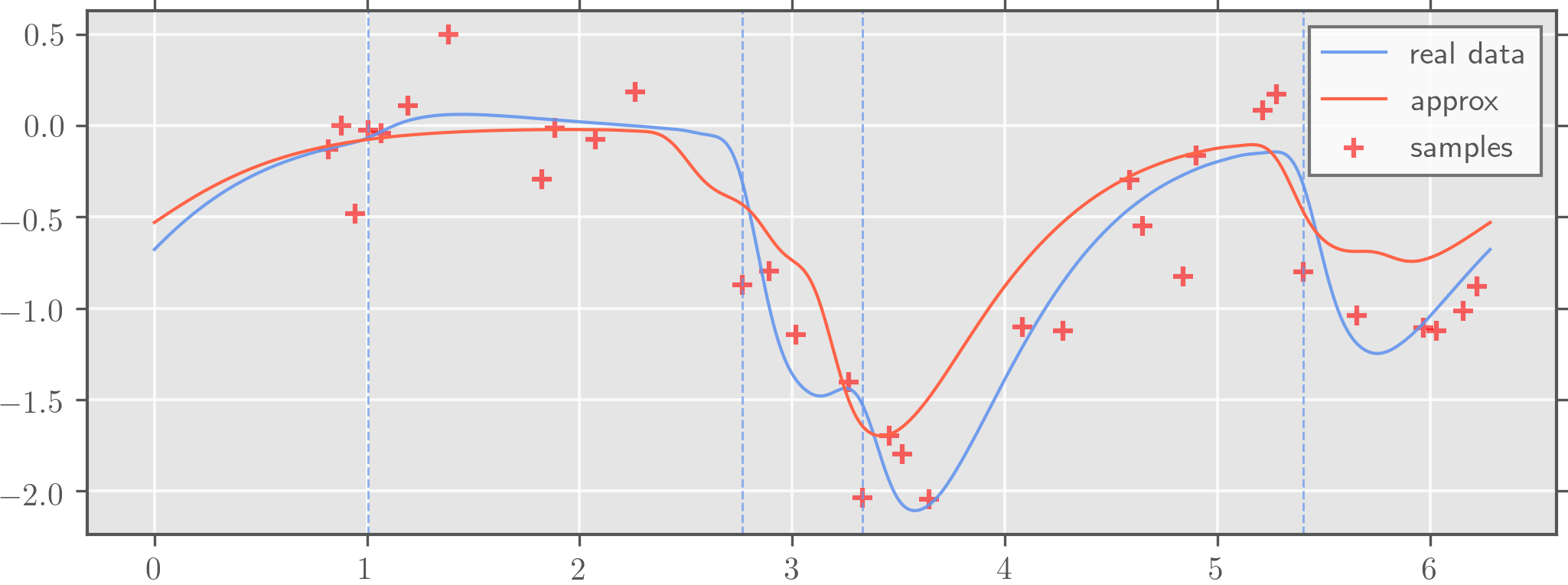}
    \caption{$\sigma$ = 0.3}
\end{subfigure}
}
\caption{CPGD reconstruction with different values of $\sigma$}
\label{fig:cpgd_reco}
\end{figure}

Let us run the CPGD algorithm with the same setting (same source spline and same input data). As previously mentioned, the parameters of CPGD are set to $M = 4K_0$ and $K = P = M$. We obtain the following results in \cref{fig:cpgd_reco} and \cref{tab:cpgd_reco}.

\begin{table}[H]
    \centering
    \begin{filecontents*}{Figures/reconstruction_cpgd/seed9970_psnr20.csv}
factors,iterations,duration,converged,objective_fun,rrse_splines,rrse_samples
0.01,500,70.61034321784973,False,1.5912881794097795,0.28578623733676467,0.2621963101904207
0.1,500,10.472101211547852,False,5.415039453056373,0.2774322541109543,0.3421569524387394
0.2,500,4.436116695404053,False,8.868348818214502,0.34042573027616,0.39929121749219143
0.3,500,4.632452964782715,False,11.899346469054334,0.438644980351382,0.4796239468270814
\end{filecontents*}
\csvreader[TableStyle]
               {Figures/reconstruction_cpgd/seed9970_psnr20.csv}{1=\factors,2=\iterations,3=\duration,4=\converged,5=\objective,6=\rrse,7=\rrsebis}%
      {\factors & \iterations & \num[round-precision=3]{\duration} & \converged & \num{\objective} & \num{\rrse} & \num{\rrsebis}}%
    \caption{Metrics for CPGD reconstruction}
    \label{tab:cpgd_reco}
\end{table}

For this exact same example, the CPGD algorithm overall performs worse then the equispaced knots reconstruction. In terms of accuracy and objective function value, it does nearly as good, however the computation times are really long. It is an encouraging result though, as this algorithm is quite new and has not been as well optimized as FISTA has been over the years.

\subsection{Frank-Wolfe algorithm}
\begin{figure}[H]
\makebox[\textwidth][c]{
\begin{subfigure}{0.6\textwidth}
    \centering
    \includegraphics[width=0.95\linewidth,center]{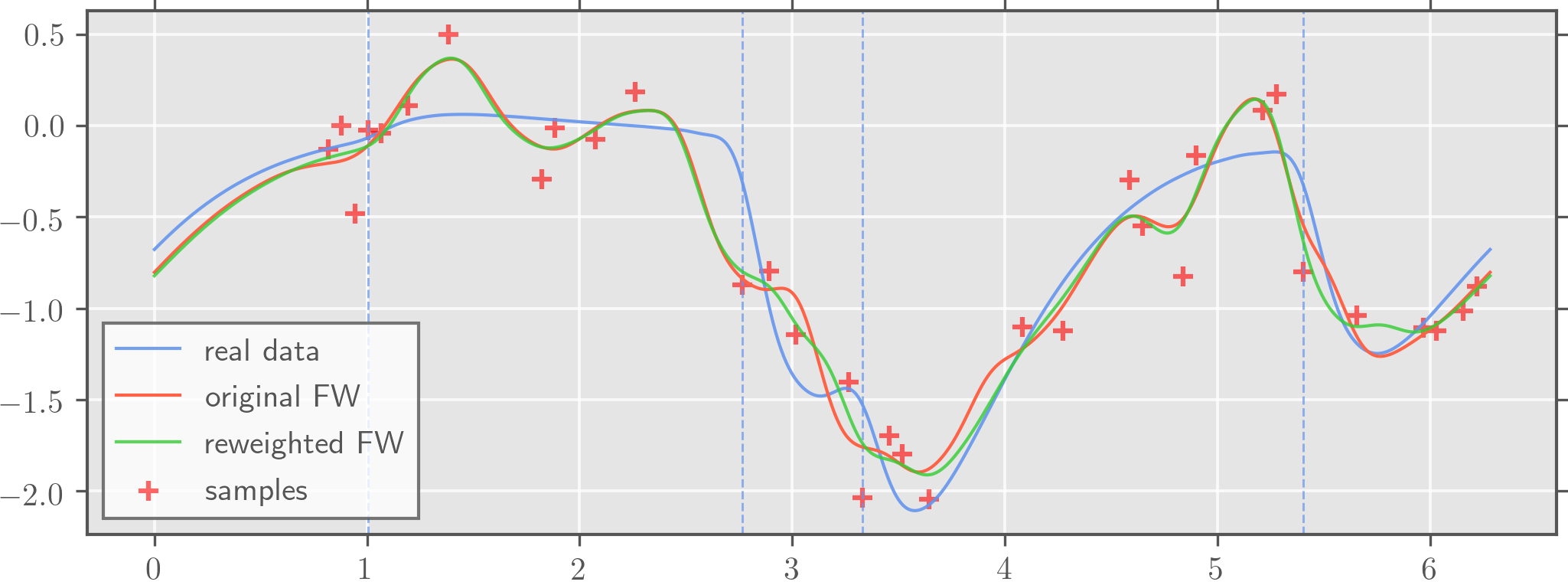}
    \caption{$\sigma$ = 0.01}
\end{subfigure}
\begin{subfigure}{0.6\textwidth}
    \centering
    \includegraphics[width=0.95\linewidth, center]{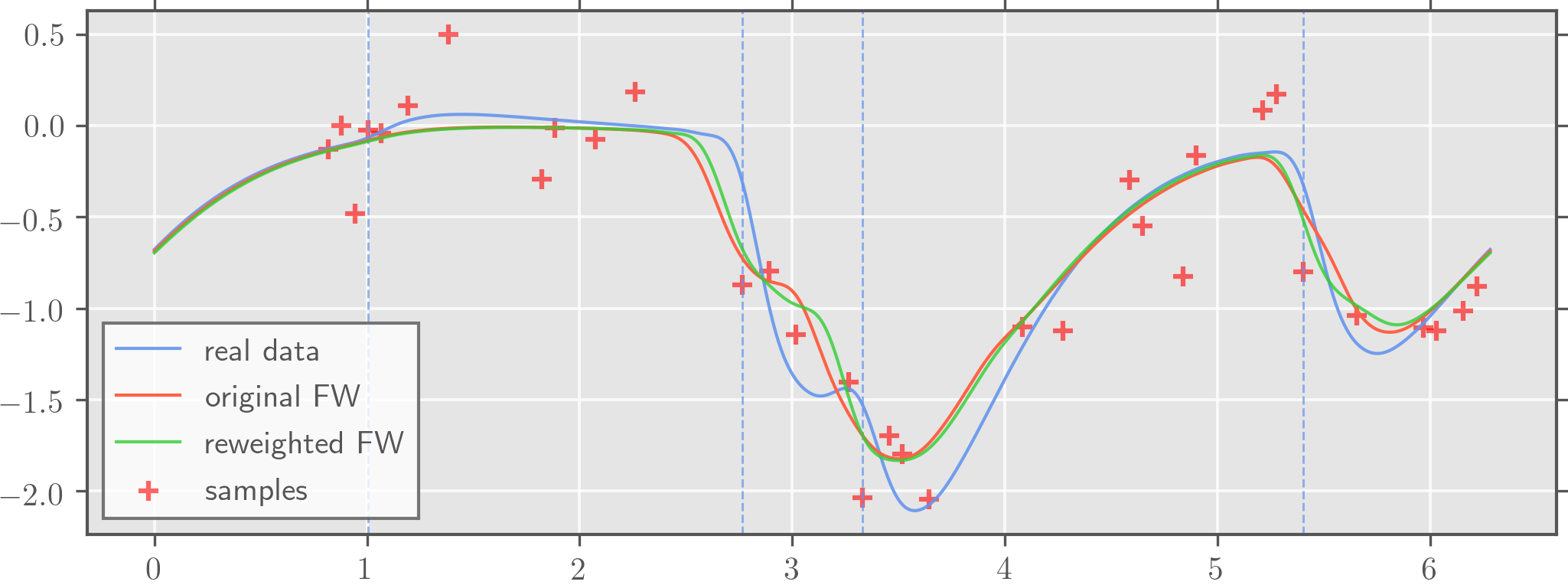}
    \caption{$\sigma$ = 0.1}
\end{subfigure}
}

\makebox[\textwidth][c]{
\begin{subfigure}{0.6\textwidth}
    \centering
    \includegraphics[width=0.95\textwidth]{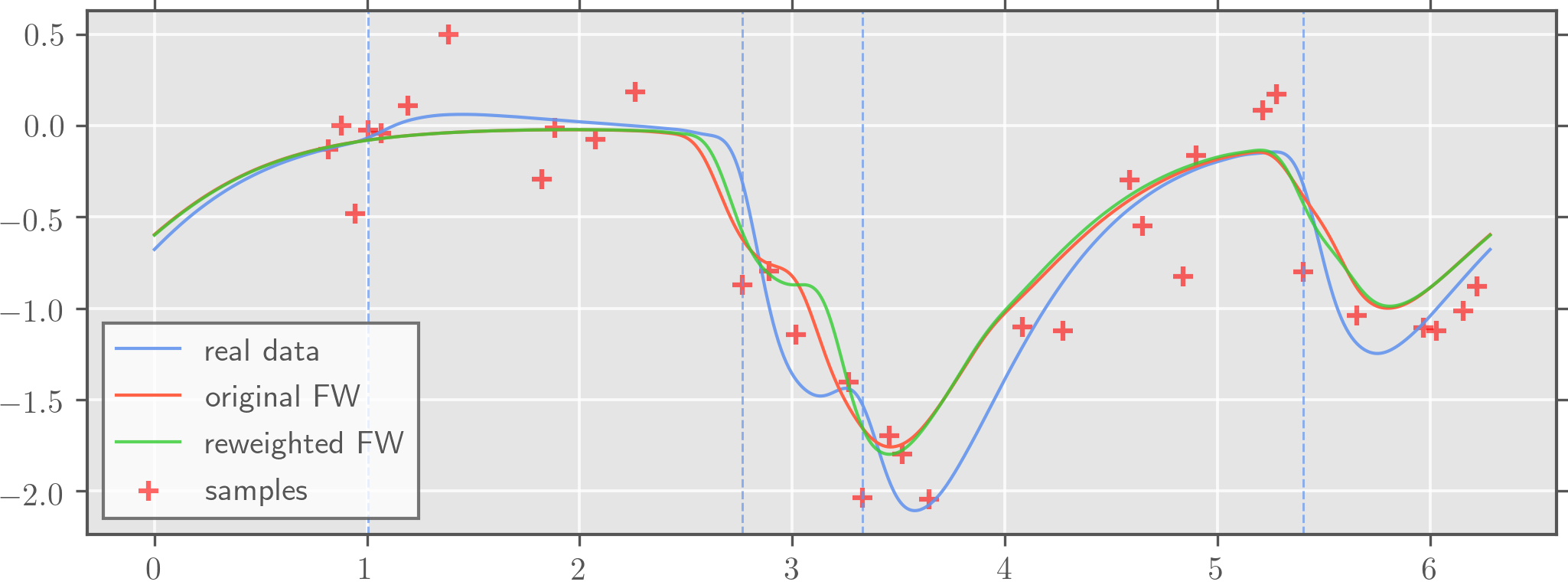}
    \caption{$\sigma$ = 0.2}
\end{subfigure}
\begin{subfigure}{0.6\textwidth}
    \centering
    \includegraphics[width=0.95\textwidth]{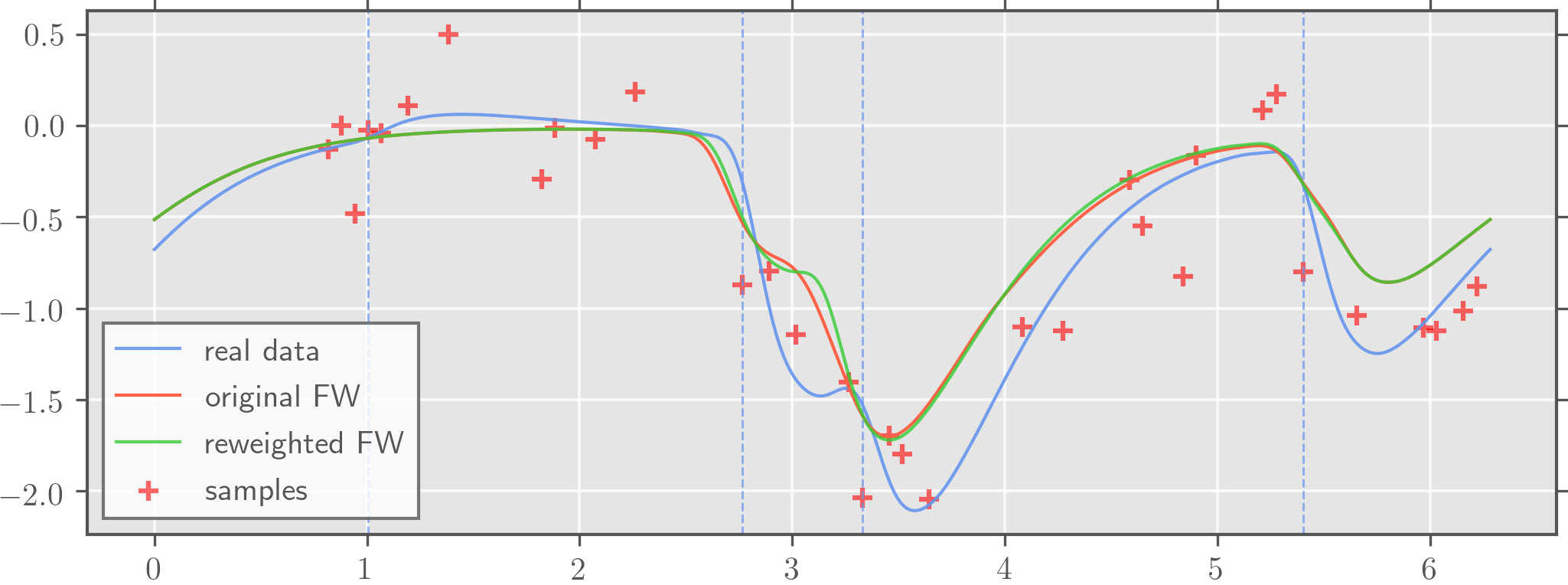}
\end{subfigure}
}
\caption{FW reconstruction with different values of $\sigma$ and the two weighting strategies}
\label{fig:fw_reco}
\end{figure}

The same way as before, we apply our implementation of FW to the same input data. In this section, we plot the results for both the original FW and our reweighting version in \cref{fig:fw_reco},  \cref{tab:fw_reco_regular} and \cref{tab:fw_reco_reweighting}.

\begin{table}[H]
    \centering
    \begin{filecontents*}{Figures/reconstruction_fw/seed9970_psnr20_comparison.csv}
factors,iterations,duration,converged,objective_fun,rrse_splines,rrse_samples,iterations_r,duration_r,converged_r,objective_fun_r,rrse_splines_r,rrse_samples_r
0.01,517,117.51474404335022,True,1.8646040305945246,0.29343314105522234,0.2409462287845656,52,11.26450514793396,True,1.386642864193172,0.2806890811124189,0.21767564498190084
0.1,161,27.188684463500977,True,5.441389955913783,0.2445737841724604,0.3414160726486718,23,4.349196195602417,True,5.262103632285674,0.22938672477971853,0.33380165283897284
0.2,107,17.101356983184814,True,8.751095175580783,0.31110122452086897,0.3935572645534868,19,3.513705015182495,True,8.613054753600592,0.31988326596545663,0.39347772830550365
0.3,263,43.48740863800049,True,11.680910735777001,0.39321319221240997,0.45900643122734247,16,2.8214662075042725,True,11.584813348948458,0.3996540528398726,0.4630076548308084
\end{filecontents*}
\csvreader[TableStyle]
               {Figures/reconstruction_fw/seed9970_psnr20_comparison.csv}{1=\factors,2=\iterations,3=\duration,4=\converged,5=\objective,6=\rrse,7=\rrsebis}%
      {\factors & \iterations & \num[round-precision=3]{\duration} & \converged & \num{\objective} & \num{\rrse} & \num{\rrsebis}}%
    \caption{Metrics for regular FW reconstruction}
    \label{tab:fw_reco_regular}
\end{table}

\begin{table}[H]
    \centering
    \csvreader[TableStyle]
               {Figures/reconstruction_fw/seed9970_psnr20_comparison.csv}{1=\factors,8=\iterations,9=\duration,10=\converged,11=\objective,12=\rrse,13=\rrsebis}%
      {\factors & \iterations & \num[round-precision=3]{\duration} & \converged & \num{\objective} & \num{\rrse} & \num{\rrsebis}}%
    \caption{Metrics for FW reconstruction with reweighting}
    \label{tab:fw_reco_reweighting}
\end{table}

In comparison with the other algorithm, FW performs well and presents same quality results. It provides really good accuracy, but at a price of longer computation times (compared to the equispaced grid algorithm). What is really striking however is the improvement obtained thanks to our reweighting strategy. When the regular convex combination weighting strategy makes the convergence really slow ($\sigma = 0.01$ or $\sigma=0.3$), the reweighting allows significant speed gains. Additionally, the final value of the objective function is either the lowest achieved among the three methods or as good as the one obtained with the first strategy. One might think it is more interesting to chose the grid approximation then, as it is still $10$ times faster, this idea needs to be mitigated. Indeed, we are working with a really simple one dimensional context. With more complex situations, like multidimensional frameworks, the grid solution might have more trouble to scale out than FW. More research is necessary in this way to tackle this question. 

\begin{figure}[H]
\makebox[\textwidth][c]{
\begin{subfigure}{0.6\textwidth}
    \centering
    \includegraphics[width=0.95\linewidth,center]{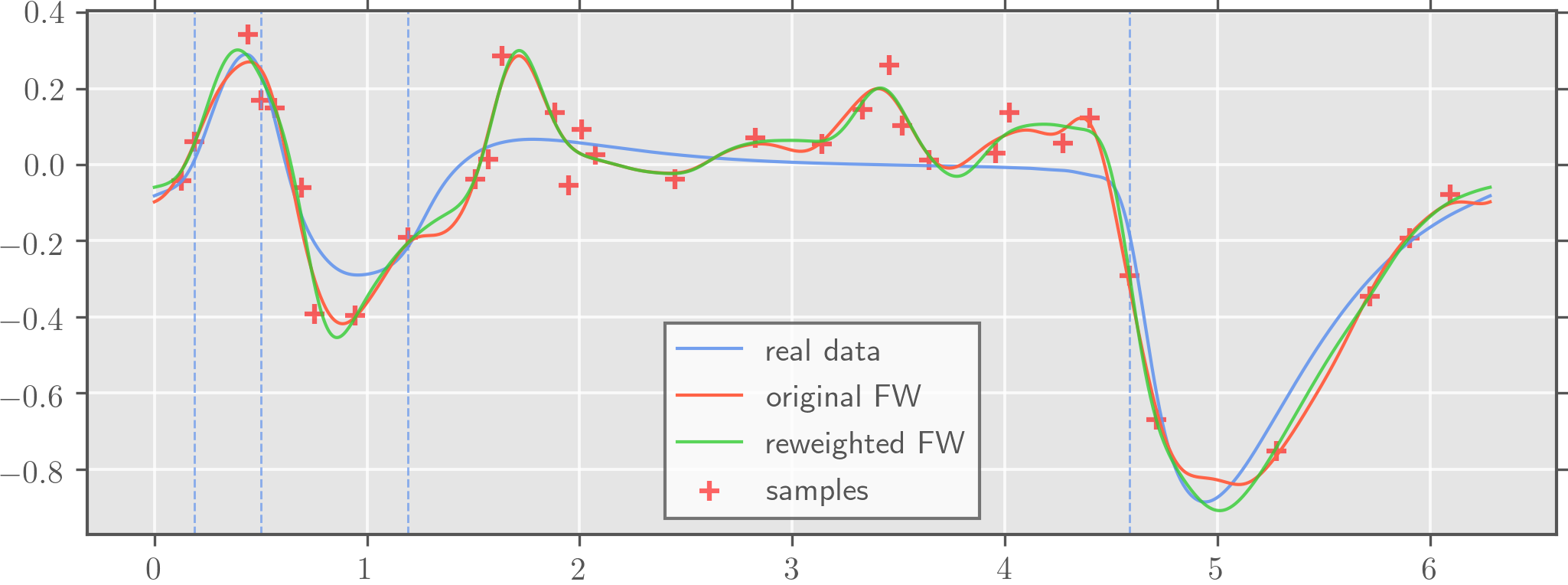}
    \caption{$\sigma$ = 0.01}
\end{subfigure}
\begin{subfigure}{0.6\textwidth}
    \centering
    \includegraphics[width=0.95\linewidth, center]{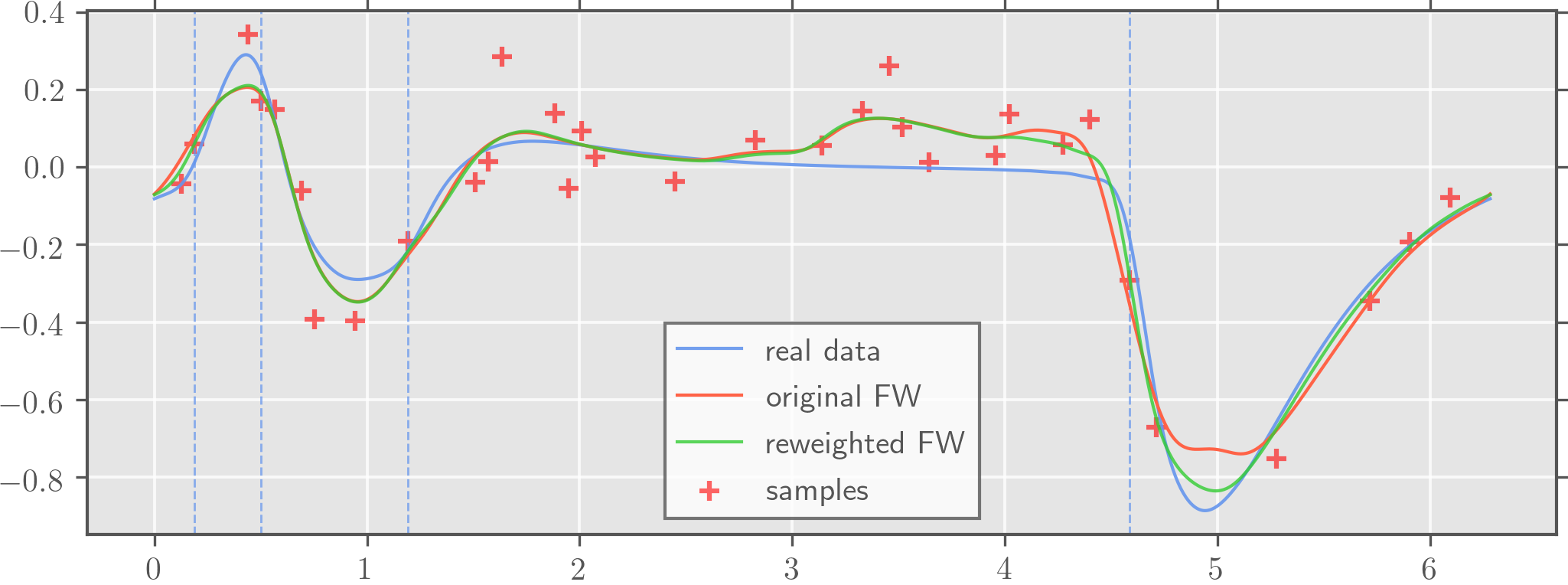}
    \caption{$\sigma$ = 0.1}
\end{subfigure}
}

\makebox[\textwidth][c]{
\begin{subfigure}{0.6\textwidth}
    \centering
    \includegraphics[width=0.95\textwidth]{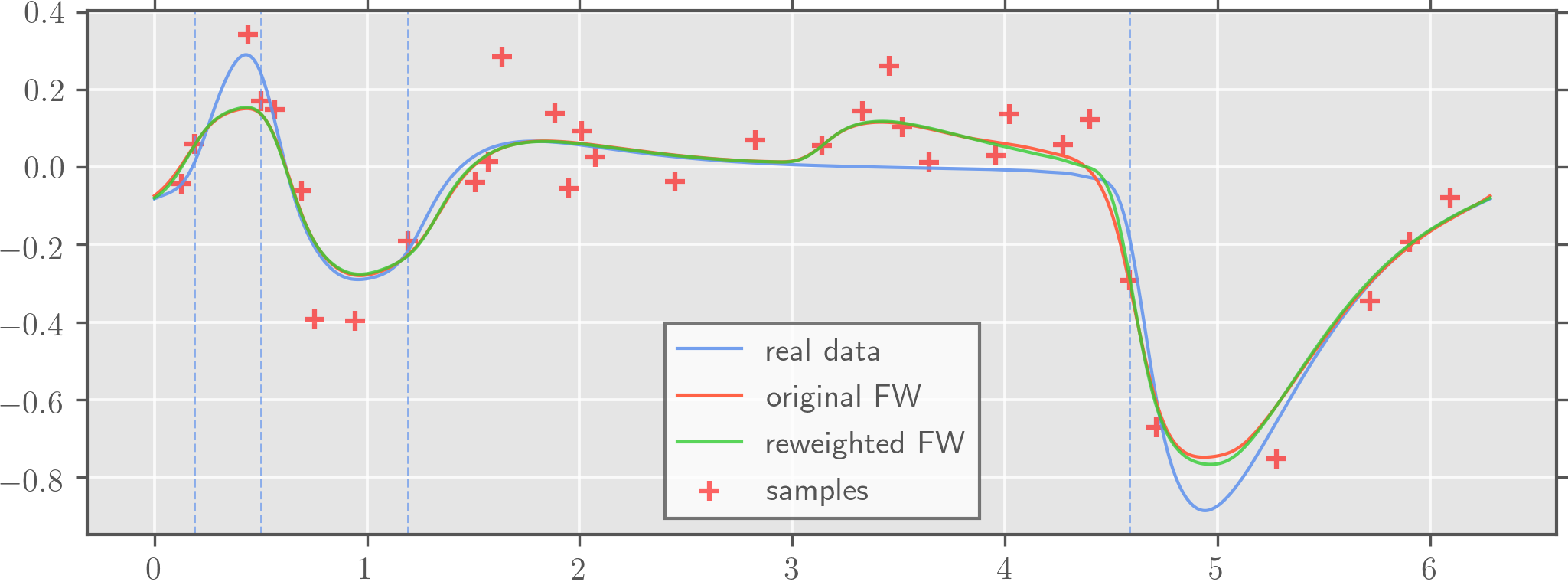}
    \caption{$\sigma$ = 0.2}
\end{subfigure}
\begin{subfigure}{0.6\textwidth}
    \centering
    \includegraphics[width=0.95\textwidth]{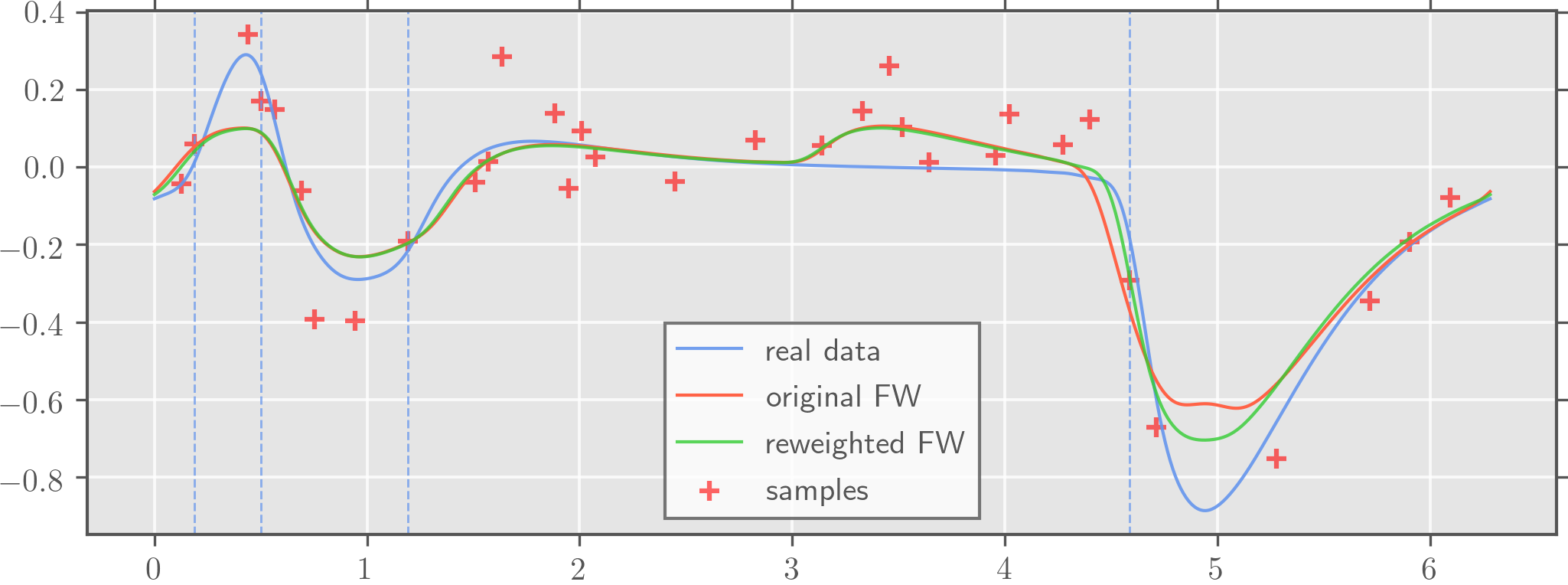}
    \caption{$\sigma$ = 0.3}
\end{subfigure}
}
\caption{FW reconstruction with different values of $\sigma$ and the two weighting strategies}
\label{fig:fw_reco_bis}
\end{figure}

We would like to emphasize the benefits of this simultaneous complete reweighting strategy within FW with a more striking example (see \cref{fig:fw_reco_bis}, \cref{tab:fw_reco_regular_bis} and \cref{tab:fw_reco_reweighting_bis}).

\begin{table}[H]
    \centering
    \begin{filecontents*}{Figures/fw_reweighting/seed4227_psnr20_comparison.csv}
factors,iterations,duration,converged,objective_fun,rrse_splines,rrse_samples,iterations_r,duration_r,converged_r,objective_fun_r,rrse_splines_r,rrse_samples_r
0.01,1957,552.350748538971,True,0.1559913566909346,0.3099941093533136,0.19513102680278765,66,15.981162071228027,True,0.12827618396128584,0.30562270055114743,0.17218804589290587
0.1,187,33.294106006622314,True,0.5142477544186126,0.24950983975540547,0.31226712304305493,29,5.829999923706055,True,0.4883227573879315,0.19726461613072077,0.29965710801387146
0.2,792,154.41544651985168,True,0.7697381992046017,0.2144165422183713,0.36978733500688676,24,4.518265247344971,True,0.7601143147865986,0.20292796146666206,0.3686206541641012
0.3,47,7.535382986068726,True,1.0164208308928484,0.3212769306202575,0.43968439559816935,19,3.5424416065216064,True,0.9855896630482037,0.2569999029062312,0.4263618381139138
\end{filecontents*}
\csvreader[TableStyle]
               {Figures/fw_reweighting/seed4227_psnr20_comparison.csv}{1=\factors,2=\iterations,3=\duration,4=\converged,5=\objective,6=\rrse,7=\rrsebis}%
      {\factors & \iterations & \num[round-precision=3]{\duration} & \converged & \num{\objective} & \num{\rrse} & \num{\rrsebis}}%
    \caption{Metrics for regular FW reconstruction}
    \label{tab:fw_reco_regular_bis}
\end{table}

\begin{table}[H]
    \centering
    \csvreader[TableStyle]
               {Figures/fw_reweighting/seed4227_psnr20_comparison.csv}{1=\factors,8=\iterations,9=\duration,10=\converged,11=\objective,12=\rrse,13=\rrsebis}%
      {\factors & \iterations & \num[round-precision=3]{\duration} & \converged & \num{\objective} & \num{\rrse} & \num{\rrsebis}}%
    \caption{Metrics for FW reconstruction with reweighting}
    \label{tab:fw_reco_reweighting_bis}
\end{table}

In this example, the improvement is considerable, regarding both the computation time and the error between the source and the approximated spline. It essentially comes from the bad estimation of the last negative bump from the regular FW method. When looking at the consecutive appearances of the Diracs along the iterations, one clearly sees that the algorithm is not able to correct the bad placements. On the other hand, the complete reweighting strategy allows to possibly remove or reduce the intensity of those problematic Diracs. Another counter effect of the convex combination weighting is that, iteration after iteration, all the previously added Diracs have their intensity being reduced. Thus, at some point, the algorithm has to reintroduce weight to those Diracs by placing spikes at nearly the exact same position of previous ones. It leads to clusters of medium intensity Diracs at the recovering position, instead of a single high intensity one, that both slows down the procedure and damages the sparsity index. This phenomenon is also present with reweighting but much less frequent. \\

To summarize this chapter, we would like to underline the fact that we are provided with three different algorithms that all seem to successfully tackle the same B-LASSO problem. While CPGD does not seem to be mature enough to demonstrate good performances in the above presented context, it is too early to conclude that it should not be considered for the problem. Indeed, one should have noticed that the evaluation of any method strongly is task-dependent. A lot of parameters are susceptible to lead towards one or another algorithm, for instance the noise level, the type of spline and operator involved, the number of available samples... A noteworthy point that arises from this remark is that this approximation framework is very versatile. It is possible to handle of large range of application situations, and the understanding of the algorithms would really benefits from real data applications. Many questions have not been studied, like the choice a reconstruction operator for an unknown input signal or the choice of a good penalization parameter $\lambda$.

%% file: Chapters/chapter4.tex
\chapter{Theoretical analysis of the solution set}

\label{chap:dual_analysis} 

In optimization theory, the convexity of the problems is a first-order property. From a numerical point of view, it makes solving the problems generally significantly easier, with the absence of local optimums. It also is of great help when theoretical analysis is involved. Indeed, convex problems are particularly well fitted for dual analysis. Under certain conditions, it is possible to obtain very insightful results as existence of solutions, optimal value of the objective function or even uniqueness of the solution. Great work has been done in this way, using duality analysis to derive properties in the general case \cite{duval2014exact} or in more particular contexts \cite{fageot2020uniqueness}. In this chapter, we first conduct the dual analysis of our problem of interest before presenting an interesting context that would lead to uniqueness in the solution set. We link the conditions of such a context to the special case of our invertible pseudo-differential operators.


\section{Duality Theory}

As explained earlier, because of the invertibility of the operator, solving the penalized problem

\begin{equation}
    \min_{f \in \mathcal{M_\mathcal{D}}(\mathbb{T})} \norm{\mathbf{y} - \mathbf{\Phi}(f)}_2^2 + \lambda \norm{\mathcal{D}f}_{TV}
    \label{eq:penalized_function}
\end{equation}
is equivalent to solving the following problem expressed over the measure space
\begin{equation}
    \min_{m\in\mathcal{M}(\mathbb{T})} \mathrm{P}_\lambda (m) := \norm{\mathbf{y} - \mathbf{\tilde{\Phi}}(m)}_2^2 + \lambda \norm{m}_{TV}.
    \label{eq:penalized_measure}
\end{equation}
The differential operator and the shape of the solution now lie in the measurement operator $\mathbf{\tilde{\Phi}}$. In this section, we will focus on the \textit{constrained} problem, which is expressed using the same measurement operator. We will call it the \textit{primal problem}.

\begin{equation}
    \min_{\substack{m\in\mathcal{M}(\mathbb{T})\\ \mathbf{y}_0 = \mathbf{\tilde{\Phi}}(m)}} \norm{m}_{TV}.
    \label{eq:primal}
\end{equation}
It is actually a reformulation of the following problem, using the invertibility of the operator.
\begin{equation}
    \min_{\substack{f\in\mathcal{M}_\mathcal{D}(\mathbb{T})\\ \mathbf{y}_0 = \mathbf{\Phi}(f)}} \norm{\mathcal{D}f}_{TV}.
    \label{eq:primal_function}
\end{equation}
The equivalence between problems \eqref{eq:penalized_measure} and \eqref{eq:primal} is actually proven \cite[Theorem~5]{gupta2018continuous} by the fact that, for any solution $m^\star$ of \eqref{eq:penalized_measure}, there exists a vector $\mathbf{y}_\lambda$ such that $m^\star$ is solution of \ref{eq:primal} with measurement vector equal to $\mathbf{y}_\lambda$. Additionally, $\mathbf{y}_\lambda$ is unique and independent from the solution $m^\star$.

We suppose that the measurement vector $\mathbf{y}_0$ is obtained by applying the measurement operator to a given measure $m_0 \in \mathcal{M}\left(\mathbb{T}\right)$ such that $\mathbf{y}_0 = \mathbf{\tilde{\Phi}}(m_0)$

\subsection{Derivation of the dual problem}

The first step of the analysis is to find the dual problem associated to \eqref{eq:primal}. By identifying the good functions and variables, it is possible to directly apply the theory from \cite[Chapter~3]{ekeland1999convex}.

\begin{proposition}
The dual problem of \eqref{eq:primal} is
 \begin{equation}
     \sup_{\substack{\norm{\mathbf{\tilde{\Phi}}^*\mathbf{p}}_\infty \leq 1 \\ \mathbf{p}\in\mathbb{R}^L}} {\langle \mathbf{y}_0, \mathbf{p} \rangle} .
     \label{eq:dual}
 \end{equation}
\end{proposition}

\begin{proof}
According to the \cite[Chapter~3.4]{ekeland1999convex}, is an optimization problem can be written as
\begin{equation*}
    \tag{P}
    \min_{x\in X}{F(x) + G(\Lambda x)},
    \label{eq:general_primal}
\end{equation*}
with $F : X \to \mathbb{R}$, $G: Y \to \mathbb{R}$ and $\Lambda : X \to Y$ linear, then one can introduce an antedual variable $y\in Y$ such that the $\Phi$ function (from Ekeland and Téman) becomes
\begin{equation*}
    \Phi(x, y) = F(x) + G(\Lambda x - y).
\end{equation*}
The dual problem of \eqref{eq:general_primal} is
\begin{equation*}
    \tag{D}
    \sup_{y^* \in Y^*}{-\left(F^*(\Lambda^*y^*) + G^*(-y^*)\right)},
    \label{eq:general_dual}
\end{equation*}
where $F^*: X^* \to \mathbb{R}$ and $G^*: Y^* \to \mathbb{R}$ denote the Fenchel conjugates of the functions $F$ and $G$, $\Lambda^*: Y^* \to X^*$ is the adjoint operator of $\Lambda$.

In our case, we shall use $F = \norm{\cdot}_{TV}$, $G = \iota_{\mathbf{y}_0}$ and $\Lambda = \mathbf{\tilde{\Phi}}$ such that $X = \mathcal{M}\left(\mathbb{T}\right)$, $X^* = \mathcal{C}$ and $Y = Y^* = \mathbb{R}^L$. $\iota$ is the indicator function over $\mathbb{R}^L$, which value is $0$ if the argument is equal to $\mathbf{y}_0$, $+\infty$ otherwise.
Thus, the dual problem of \eqref{eq:primal} becomes
\begin{equation}
    \sup_{\mathbf{p} \in \mathbb{R}^L} {-\norm{\mathbf{\tilde{\Phi}^* p}}_{TV}^* - \iota_{\mathbf{y}_0}^*(-\mathbf{p})}.
    \label{eq:intermediate_dual}
\end{equation}
Let us express the conjugate functions that play a role in the above-mentioned expression.

\begin{lemma}
    The Fenchel conjugate function of the TV-norm $\norm{\cdot}_{TV}: \mathcal{C} \to \mathbb{R} $ is given by the indicator function $\iota_{\mathcal{B}_\infty}$ of $\mathcal{B}_\infty = \left\{ \eta\in\mathcal{C}, \norm{\eta}_\infty \leq 1\right\}$ the unit ball for the infinite norm in $\mathcal{C}$.
\end{lemma}
Indeed, if $F = \norm{\cdot}_{TV}$, then for $\eta\in\mathcal{C}$:
\begin{equation*}
    F^*(\eta) = \sup_{m\in\mathcal{M}} \langle \eta, m \rangle - \norm{m}_{TV}
\end{equation*}
We have the fundamental inequality $\langle \eta, m \rangle \leq \norm{\eta}_\infty \norm{m}_{TV}$. Thus, if $\norm{\eta}_\infty \leq 1$ then $\langle \eta, m \rangle - \norm{m}_{TV} \leq \left( \norm{\eta}_\infty - 1 \right) \norm{m}_{VT} \leq 0$, and one reach $0$ for $m = 0$.
Otherwise, if $\norm{\eta}_\infty > 1$, there exists $m_0$ such that $\langle \eta, m_0 \rangle \geq a \norm{m_0}_{TV}$ with $a>1$ (for instance, $m_0 = \Sha(t-\argmax{\eta})$). Consequently, for $\lambda > 0$, $F^*(\eta) \geq \langle \eta, \lambda m_0 \rangle - \norm{\lambda m_0}_{TV} \geq \lambda a \norm{m_0}_{TV} \xrightarrow[\lambda \to + \infty]{} + \infty$.

\begin{lemma}
    The Fenchel conjugate function of $G = \iota_{\mathbf{y}_0}$ is given by $G^*(\mathbf{h}) = \langle \mathbf{y}_0, \mathbf{h} \rangle$ for $\mathbf{h}\in\mathbb{R}^L$.
\end{lemma}
Similarly, we have
\begin{equation*}
    G^*(\mathbf{h}) = \sup_{\mathbf{z}\in\mathbb{R}^L} \langle \mathbf{z}, \mathbf{h} \rangle - \iota_{\mathbf{y}_0}\left(\mathbf{z}\right)
\end{equation*}
Directly, the following holds:
\begin{equation*}
    \langle \mathbf{z}, \mathbf{h} \rangle - \iota_{\mathbf{y}_0}\left(\mathbf{z}\right)= \left\{
    \begin{aligned}
        &\langle \mathbf{y}_0, \mathbf{h} \rangle \qquad &\mathrm{if} \quad \mathbf{z} = \mathbf{y}_0 \\
        &- \infty &\mathrm{otherwise}
    \end{aligned}
    \right.
\end{equation*}
which proves the result.

Given these two lemmas, one can rewrite \eqref{eq:intermediate_dual} as
\begin{equation*}
    \sup_{\mathbf{p} \in \mathbb{R}^L} {-\iota_{\mathcal{B}_\infty}\left(\mathbf{\tilde{\Phi}^* p}\right) + \langle \mathbf{y}_0, \mathbf{p} \rangle.}
\end{equation*}
which leads to the given expression. Note that it is also possible to directly derive the dual problem without applying the previous theory (problems \eqref{eq:general_primal} and \eqref{eq:general_dual}) by introducing the equivalent of the Lagrange dual function (using duality product instead of inner product). It uses the same calculations of the Fenchel conjugate functions.
\end{proof}

The dual problem is actually much simpler than the primal one, and we will see in the next paragraphs what information it provides about the solution of the initial problem.


\subsection{Strong duality}

The key property that makes dual analysis so insightful for our problem is the strong duality. Relying on \cite{ekeland1999convex}, Duval and Peyré proved that strong duality holds between the primal and dual problems presented earlier (\cite[Proposition~13]{duval2014exact}). Let us summarize the implications here.

\begin{proposition}
Strong duality holds between \eqref{eq:primal} and \eqref{eq:dual}, such that if there exists a solution $\mathbf{p^\star}$ to the dual problem \eqref{eq:dual}, then the primal problem \eqref{eq:primal} also  admits solutions.

\noindent Additionally, the optimal values of both problems coincide, i.e.
\begin{equation*}
    \min_{\substack{m\in\mathcal{M}(\mathbb{T})\\ \mathbf{y}_0 = \mathbf{\tilde{\Phi}}(m)}} \norm{m}_{TV} = \sup_{\norm{\mathbf{\tilde{\Phi}}^*\mathbf{p}}_\infty \leq 1} {\langle \mathbf{y}_0, \mathbf{p} \rangle}
\end{equation*}

\noindent Finally, one can derive extremality conditions. If a solution $\mathbf{p^\star}$ exists, then
\begin{equation}
    \mathbf{\tilde{\Phi}^*} \mathbf{p^\star} \in \partial \norm{m^\star}_{TV}
    \label{eq:extremality_cond}
\end{equation}
where $m^\star$ is any solution of the primal \eqref{eq:primal}. Conversely, if \eqref{eq:extremality_cond} holds, then $\left(\mathbf{p}^\star, m^\star \right)$ is a joint couple of solutions to respectively the primal and dual problems.
\end{proposition}

The strong duality comes really useful for our concern. Indeed, the dual problem is actually simpler than the primal one, and we can adapt the proof from \cite[Proposition~2]{duval2014exact} to ensure existence of solutions.

\begin{proposition}
Let $m_0 \in\mathcal{M}\left(\mathbb{T}\right)$ and $\mathbf{y}_0 = \mathbf{\tilde{\Phi}}m_0 \in \mathbb{R}^L$. There exists a solution to the dual problem 
\begin{equation*}
    \sup_{\substack{\norm{\mathbf{\tilde{\Phi}}^*\mathbf{p}}_\infty \leq 1 \\ \mathbf{p}\in\mathbb{R}^L}} {\langle \mathbf{y}_0, \mathbf{p} \rangle}.
\end{equation*}
\end{proposition}
\begin{proof}
As $\mathbf{y}_0 = \mathbf{\tilde{\Phi}}m_0$ the objective function can be rewritten as
\begin{equation*}
    \sup_{\substack{\norm{\eta}_\infty \leq 1 \\ \eta\in\mathcal{R}\left(\mathbf{\tilde{\Phi}}^*\right)}} {\langle m_0, \eta \rangle} \quad \leq \quad \norm{m_0}_{TV}
\end{equation*}
where $\mathcal{R}(\cdot)$ denotes the range of an operator, and thus is bounded. We can consider $\left(\eta_n\right)_{n\in\mathbb{N}}$ a maximizing sequence of the problem. Then $\left(\eta_n\right)_{n\in\mathbb{N}}$ lives in the finite dimensional linear space $\mathcal{R}\left(\mathbf{\tilde{\Phi}}^*\right) = \mathrm{Span}\left(\mathbf{\tilde{\Phi}}^*_i\right)_{i=1, \dots, L}$. As $\left(\eta_n\right)_{n\in\mathbb{N}}$ is bounded we may extract a subsequence converging to $\eta^\star \in \mathcal{C}\left(\mathbb{T}\right)$. $\eta^\star$ satisfies both $\norm{\eta^\star}_\infty \leq 1$ and $\eta^\star \in \mathcal{R}\left(\mathbf{\tilde{\Phi}}^*\right)$ so that $\eta^\star = \mathbf{\tilde{\Phi}}^* \mathbf{p}^\star$ and $\mathbf{p}^\star \in \mathbb{R}^L$ is solution of the dual problem.
\end{proof}
Thus, as long as the data $\mathbf{y}_0$ is sampled from an actual measure $m_0 \in\mathcal{M}\left(\mathbb{T}\right)$, we have an assurance that the primal problem \eqref{eq:primal} admits solutions. What is even more insightful is that, in this case, it is possible to link the solutions of the primal and the dual problem through the \textit{extremality conditions}. Duval and Peyré have derived a closed form expression for the subdifferential of the TV-norm \cite[Proposition~12]{duval2014exact}, given by:
\begin{equation}
    \forall m \in \mathcal{M}\left(\mathbb{T}\right), \quad \partial\norm{m}_{TV} = \left\{ \eta\in\mathcal{C}\left(\mathbb{T}\right), \norm{\eta}_{TV}\leq 1 \quad \mathrm{and} \quad \int \eta\  \mathrm{d}m = \norm{m}_{TV} \right\}
    \label{eq:subdifferential_TV}
\end{equation}
Then, for problem \eqref{eq:primal}, there exists a function $\eta\in\mathcal{C}\left(\mathbb{T}\right)$ coined as \textit{dual certificate} that verify the three following properties:

\begin{enumerate}[label=(\roman*)]
    \item \label{itm:i} $\eta \in \mathcal{R}\left(\mathbf{\tilde{\Phi}}^*\right) = \mathrm{Span}\left(\mathbf{\tilde{\Phi}}^*_i\right)_{i=1, \dots, L}$
    \item \label{itm:ii} $\norm{\eta}_\infty \leq 1$
    \item \label{itm:iii} for any solution $m^\star$ of \eqref{eq:primal}, \quad $\langle \eta, m^\star \rangle = \norm{m^\star}_{TV}$
\end{enumerate}
Note that such a function $\eta$ might not be unique, and any function satisfying the three above mentioned requirements is called a \textit{dual certificate}. We could aim for uniqueness by choosing for example the certificate made from the minimal norm $\mathbf{p}$ (as made in \cite{duval2014exact}). The name comes from the fact that the existence of such a function \emph{certifies} that a solution to the primal problem exists.


\subsection{Dual certificate}

Let us consider a function $\eta$ that is a dual certificate for problem $\eqref{eq:primal}$. $\eta$ is a continuous periodic function with infinite norm being 1, so it necessarily reaches at least once $\pm 1$ over its period. We can thus consider what we call its \textit{saturation points} defined as:
\begin{equation}
    \mathrm{Sat}(\eta) = \left\{ t\in \mathbb{T}, \eta(t) = \pm 1 \right\}
\end{equation}
and its \textit{signed saturation points}:
\begin{equation}
    \mathrm{Sat}^\pm(\eta) = \left\{ (t, v) \in \mathbb{T}\times \{-1, +1\}, \eta(t) = v \right\}.
\end{equation}

\noindent Similarly, we can define the \textit{signed support} of a measure $m\in\mathcal{M}\left(\mathbb{T}\right)$ as 
\begin{equation}
    \mathrm{Supp}^\pm(m) = \left(m^{-1}(\mathbb{R}^*_+) \times \{+1\}\right) \cup \left(m^{-1}(\mathbb{R}^*_-) \times \{-1\}\right) \subset \mathbb{T} \times \{ \pm 1 \}.
\end{equation}
Then, one can state \cite[Proposition~3]{duval2014exact} :
\begin{proposition}
Let $m_0\in\mathcal{M}\left(\mathbb{T}\right)$ and $\mathbf{y}_0 = \mathbf{\tilde{\Phi}}m_0$. The following equivalence holds:
\begin{center}
\begin{tabular}{ c }
    $m_0$ is a solution to \quad $\min_{\substack{m\in\mathcal{M}(\mathbb{T})\\ \mathbf{y}_0 = \mathbf{\tilde{\Phi}}(m)}} \norm{m}_{TV}$ \\[.8cm]
    if and only if \\[.5cm]
    $\mathrm{Supp}^\pm(m_0) \subset \mathrm{Sat}^\pm(\eta)$\\
    with $\eta$ being a dual certificate of the problem.
\end{tabular}
\end{center}
\label{prop:saturation}
\end{proposition}

\begin{proof}
For the direct implication, let $\eta\in\mathcal{C}\left(\mathbb{T}\right)$ be a dual certificate and suppose there exists $t_1\in\mathrm{Supp}^\pm(m_0)$ but $t_1 \notin\mathrm{Sat}^\pm(\eta)$. Then $|\eta(t)| < 1$ and, as $\eta$ is a dual certificate, the following holds:
\begin{equation*}
    \langle \eta, m_0 \rangle = \norm{m_0}_{TV} = \sup_{\psi\in\mathcal{C}\left(\mathbb{T}\right), \norm{\psi}_{TV}\leq 1}{\int \psi \mathrm{d}m_0}
\end{equation*}
However, we could consider a function $\psi_1\in\mathcal{C}\left(\mathbb{T}\right)$ such that $\psi_1(t_1) = 1$ and $|\psi_1| \geq |\eta|$ everywhere on $\mathbb{T}$. Then $\norm{m_0}_{TV} \geq \langle m_0, \psi_1 \rangle > \langle m_0, \eta \rangle = \norm{m_0}_{TV}$ which is impossible.

On the other way, suppose $\mathrm{Supp}^\pm(m_0) \subset \mathrm{Sat}^\pm(\eta)$, then $\langle m_0, \eta \rangle = \norm{m_0}_{TV}$ which is the optimal value of the dual problem. By strong duality, is is also the optimal value of the primal problem, and it is reached with $m=m_0$, so $m_0$ is a solution of the primal.
\end{proof}
\noindent This result is crucial in our study, as any dual certificate constrains the location of the support of any solution to the primal problem.
In another manner, if we define the set of certificates
\begin{equation*}
    \mathrm{Cert}(\mathbf{y}_0) = \left\{\eta \in\mathcal{R}\left(\mathbf{\tilde{\Phi}}^*\right),\quad \eta \text{ is a dual certificate of }\eqref{eq:primal}\  \right\},
\end{equation*}
the following inclusion holds
\begin{equation}
    \left\{\bigcup_{m \text{ solution of } \eqref{eq:primal}} {\mathrm{Supp}^\pm(m)}\right\} \quad \mbox{\huge$\subset$} \quad \left\{ \bigcap_{\eta \in \mathrm{Cert}(\mathbf{y}_0)} \mathrm{Sat}^\pm(\eta)\right\}.
\end{equation} \\

We can go even further and obtain a sufficient condition that ensures uniqueness. For any dual certificate $\eta$, we can consider the set
$$\mathcal{M}_\eta\left(\mathbb{T}\right) = \left\{m \in \mathcal{M}\left(\mathbb{T}\right), \mathrm{Supp}^\pm(m) \subset \mathrm{Sat}^\pm(\eta) \right\} \subset \mathcal{M}\left(\mathbb{T}\right)$$
of the periodic Radon measures with support included in the saturation points of $\eta$ and thus the restriction $\mathbf{\tilde{\Phi}}_{\mathrm{Sat}^\pm(\eta)}$ of $\mathbf{\tilde{\Phi}}$ to $\mathcal{M}_\eta\left(\mathbb{T}\right)$. The existence of a solution to the primal problem ensures that there exists at least a measure $m_1$ that solves $\eqref{eq:primal}$, and necessarily $m_1 \in \mathcal{M}_\eta\left(\mathbb{T}\right)$. If $\mathbf{\tilde{\Phi}}_{\mathrm{Sat}^\pm(\eta)}$ would happen to be injective, then there could not be more than one solution to the problem, $m_1$ would be unique. If $\mathrm{Sat}(\eta)$ is infinite, it actually does not bring much information, as $\mathbf{\tilde{\Phi}}_{\mathrm{Sat}^\pm(\eta)}$ remains as complex as $\mathbf{\tilde{\Phi}}$. However, it becomes particularly helpful when $\mathrm{Sat}(\eta)$ is finite, as proving that the measurement function is injective reduces to proving that a certain matrix is full rank. It turns the injectivity question into a finite dimensional one.


\section{Towards Uniqueness of Solutions}

So far, with very few hypotheses, it has been possible to state strong properties that directly derive from the dual analysis of the general problem $\eqref{eq:primal}$. Indeed we did not need any assumption on the class or the type of the measurement operator $\mathbf{\tilde{\Phi}}$. A really interesting work has already been achieved when a Fourier sampling operator is concerned \cite{fageot2020uniqueness}. The authors have proven that the constrained problem (the same as the one we study) either admits only positive solutions or the solution is unique. They also have obtained uniqueness for a certain type of operator applied to the regularized problem. In this section we want to examine to what extent we can leverage the previous properties in our context of spatial sampling with invertible operators.


\subsection{Spatial sampling}

When performing spatial sampling, we obtain a vector $\mathbf{y} \in \mathbb{R}^L$ of the $L$ measurements and thus lives in finite dimensional space. It leads to a quite simple expression of the adjoint measurement operator $\mathbf{\tilde{\Phi}^*}$ as:
\begin{equation}
    \mathbf{\tilde{\Phi}^*} \colon \left\{ \begin{aligned}
        \quad \mathbb{R}^L \quad & \to \quad \mathcal{C}(\mathbb{T})\\
        \mathbf{h} \quad & \mapsto \quad \sum_{\ell=1}^L {h_\ell\ \psi_\mathcal{D}(\theta_\ell - \cdot)}
    \end{aligned} \right. .
    \label{eq:adjoint_measurement}
\end{equation}
In practice, we are able to plot (and even sometimes to give a closed-form expression) of the Green's functions involved in \eqref{eq:adjoint_measurement} for some well defined pseudo differential operators $\mathcal{D}$. On the examples we have investigated on - namely exponential, Sobolev, Wendland and Matérn operators - they behave quite regularly and admit a finite number of saturation points (see the plots from Fageot and Simeoni in \cite{fageot2020tvbased}).

Then, as any dual certificate lives in $\mathcal{R}\left(\mathbf{\tilde{\Phi}^*} \right)$, it is expressed as a linear combination of (reversed) shifted Green's function, and thus we can expect it to have a finite number of saturation points. This would be quite interesting in the perspective of finding a solution to the primal problem, as one can state the following property:
\begin{proposition}
Let $\eta\in\mathcal{C}$ be a dual certificate function associated with the primal and dual problems \eqref{eq:primal} and \eqref{eq:dual}. If $\eta$ admits a finite number of saturation points, then any solution of the primal problem is a discrete measure with a number of impulses no greater than the number of saturation points.

Even stronger, the solution set is a convex and compact subset of the finite-dimensional vector space
\begin{equation*}
    \mathrm{Span}\left( \left\{ \ \Sha(\cdot - \tau)\ \colon\quad \tau \in \mathrm{Sat}(\eta)\ \right\} \right),
\end{equation*}
and, for any solution
\begin{equation*}
    m^\star = \sum_{\tau\in\mathrm{Sat}(\eta)}{\beta_\tau \Sha(\cdot - \tau)},
\end{equation*}
the sign of the weights $\beta_\tau$ are given by the sign of the associated saturation points $\tau\in\mathrm{Sat}^\pm(\eta)$.
\label{prop:discrete}
\end{proposition}
Stated another way around, any solution of the primal problem expressed in the continuous setting is a spline with a finite rate of innovations. This property naturally derives from \cref{prop:saturation} on the inclusion of the support of the solutions. \\

The uniqueness property discussed earlier could also be made explicit in the case of spatial sampling and discrete measures. Suppose we are in the conditions of \cref{prop:discrete}, let us call $P = \left| \mathrm{Sat}(\eta) \right|$ for $\eta$ a dual certificate. The space of measures $\mathcal{M}_\eta\left(\mathbb{T}\right)$ with support included in the saturation points of $\eta$ is then directly characterized by the $P$ weights of the innovations, and thus is isomorphic to $\mathbb{R}^P$. More precisely, if $\mathrm{Sat}(\eta) = \left\{\tau_1, \dots, \tau_P\right\}$, one can write:
\begin{equation}
    \mathbf{\tilde{\Phi}}_{\mathrm{Sat}^\pm(\eta)} \colon \left\{ \begin{aligned}
        \quad \mathbb{R}^P \quad & \to \quad \mathbb{R}^L\\
        \boldsymbol{\beta} \quad & \mapsto \quad  \mathbf{\tilde{\Phi}} \left(\sum_{p=1}^P {\beta_p \ \Sha( \cdot - \tau_p)}\right) = \left(\sum_{p=1}^P {\beta_i \psi_\mathcal{D}\left(\theta_\ell - \tau_p \right)}\right)_{\ell=1, \dots, L}
    \end{aligned} \right. .
    \label{eq:measurement_sat}
\end{equation}
We may reformulate:
\begin{proposition}
Let $\eta$ be a dual certificate with a finite set of saturation points $\mathrm{Sat}(\eta) = \left\{\tau_1, \dots, \tau_P\right\}$. If $\mathbf{\tilde{\Phi}}_{\mathrm{Sat}^\pm(\eta)}$ (as defined in \eqref{eq:measurement_sat}) is injective then the solution of the primal problem is unique. In terms of matrix, it is equivalent for the following matrix $\Psi_{\boldsymbol{\tau}} \in \mathbb{R}^{L\times P}$ to have full rank:
\begin{equation}
    \Psi_{\boldsymbol{\tau}} = \left( \psi_\mathcal{D}\left(\theta_\ell - \tau_p \right) \right)_{\ell, p} = \begin{pmatrix}
    \psi_\mathcal{D}\left(\theta_1 - \tau_1 \right) & \dots & \psi_\mathcal{D}\left(\theta_1 - \tau_P \right) \\
    \vdots & \ddots & \\
    \psi_\mathcal{D}\left(\theta_L - \tau_1 \right) & & \psi_\mathcal{D}\left(\theta_L - \tau_P \right)
    \end{pmatrix}.
    \label{eq:matrix_uniqueness}
\end{equation}
\label{prop:uniqueness}
\end{proposition}
Necessarily, if the number $P$ of saturation points of a certificate is greater than $L$, the matrix \eqref{eq:matrix_uniqueness} could not satisfy the full rank condition. Thus this property is only interesting when there exists a dual certificate with $L$ or less saturation points.

The conditions stated for both the previous theorems essentially depend on the operator $\mathcal{D}$. The first one amounts to estimating the set of saturation points that a linear combination of shifted Green's function with infinite norm of $1$ can reach over a period. As the combination is finite, the actual question would be : can such a function saturate to $1$ in more than a finite number of points ? The second condition is linked to the rank of a finite number of vector evaluations of the Green's function. Thus, the central point of these properties clearly is the pseudo differential operator.


\subsection{Application to the exponential operator}

During our practical experiences, we extensively have made use of exponential splines. As Fageot and Simeoni define them \cite{fageot2020tvbased}, they are linked to the operator $\mathcal{D} = \left( \mathrm{D} + \alpha \mathrm{Id}\right)^\gamma $, with $\gamma>1$ so that the Green's function is continuous and thus the operator is \textit{sample-admissible} \cite[Section~6.2]{fageot2020tvbased}. Let us consider integer $\gamma = N \in\mathbb{N}$, such that $N \geq 2$. We shall ask if such an operator, which is invertible, satisfies the requirements of \cref{prop:discrete} and \cref{prop:uniqueness}.

For the sake of readability and without loss of generality, we will now consider that $\mathbb{T} = \left[0, 1\right]$, \textit{i.e.} the period of the considered objects is $1$. Any function $g\in\mathcal{S}(\mathbb{R)}$ over the real line can be periodized, and thus be seen as a periodic generalized function $\tilde{g}\in\mathcal{S}(\mathbb{T)}\subset \mathcal{S}'(\mathbb{T)}$, with the transformation
\begin{equation*}
    \forall r \in \mathbb{R}, \qquad \tilde{g}(r) = g(r-\floor{r}).
\end{equation*}
This operation truncates the initial functon $g$ and simply loops back to the value at the beginning of the period. It does not ensure any property as continuity or regularity for the function $\tilde{g}$. In the following, we identify both $g$ and $\tilde{g}$ as the periodized element and only note $\left\{t \mapsto g(t)\right\} \in\mathcal{S}'(\mathbb{T)}$ when the variable is in $\mathbb{T}$. If the variable $r$ is taken in $\mathbb{R}$, we specifically write $g(r-\floor{r})$.

With $N=2$, it is actually possible to derive a closed-form expression for $\psi_\mathcal{D}$ knowing that $\mathcal{D} \psi_\mathcal{D} = \Sha$ and using the derivative with respect to the periodic distributions. We obtain:
\begin{equation}
    \mathcal{D} = \left( \mathrm{D} + \alpha \mathrm{Id}\right)^2 \ \Rightarrow \ \forall r \in \mathbb{R},\ \psi_\mathcal{D}(r) = \left( \frac{1}{1-e^{-\alpha}}(r - \floor{r}) + \frac{e^{-\alpha}}{(1-e^{-\alpha})^2} \right) e^{-\alpha (r - \floor{r})}
\end{equation}
In fact, the following theorem holds:
\begin{theorem}
The Green's function of the periodic exponential operator $\mathcal{D} = \left( \mathrm{D} + \alpha \mathrm{Id}\right)^N $ of order $N\in\mathbb{N}\setminus\{0\}$ is given by
\begin{equation}
    \forall r \in \mathbb{R},\ \psi_\mathcal{D}(r) = P_{N-1}(r - \floor{r}) e^{-\alpha (r - \floor{r})}
\end{equation}
where $P_{N-1}$ is a polynomial of degree $(N-1)$.

The coefficients of $P_{N-1}$ are known and can be computed recursively, starting from the leading coefficient. Let $a_k$ be the coefficients of $P_{N-1}$ according to
\begin{equation*}
    P_{N-1}(X) = a_{N-1}X^{N-1} + \dots + a_1 X + a_0 ,
\end{equation*}
then the following holds:
\begin{equation}
    \forall k \in \left\{0, 1, \dots, N-1\right\}, \qquad a_k = \frac{1}{k!} b_k
\end{equation}
and
\begin{gather}
    b_{N-1} = \frac{1}{1-e^{-\alpha}} \\[0.2cm]
    \forall k \in \left\{ 2, 3, \dots, N\right\}, \quad b_{N-k} = A \left[ \sum_{i=1}^{k-1} \frac{1}{(k-i)!} b_{N-i} \right] = A \left[ \sum_{i=1}^{k-1} \frac{1}{i!} b_{N-k+i} \right].
\end{gather}
with
\begin{equation*}
    A = \frac{e^{-\alpha}}{1-e^{-\alpha}}.
\end{equation*}
\label{th:exp_spline}
\end{theorem}

\noindent The proof requires the statement of a short lemma.
\begin{lemma}
Let $\alpha>0$, for any polynomial $P$, the function $t \mapsto P(t)e^{-\alpha t} \in \mathcal{S}'\left(\mathbb{T}\right)$ satisfies the following equality
\begin{equation}
    \left( \mathrm{D} + \alpha \mathrm{Id}\right)\left(P(t)e^{-\alpha t}\right) = \left[P(0)- P(1)e^{-\alpha}\right] \Sha + P'(t)e^{-\alpha t}.
    \label{eq:exp_to_green}
\end{equation}
where the derivative $\mathrm{D}$ is taken in sense of the periodic distributions and $P'$ is the usual derivative of the polynomial $P$.
\end{lemma}
\begin{proof}
For any periodic test function $u\in\mathcal{C}_\infty\left(\mathbb{T}\right)$ , we can write, using integration by parts:
\begin{align*}
    \Big< \left( \mathrm{D} + \alpha \mathrm{Id}\right)\left(P(t) e^{-\alpha t}\right)&, u \Big> = \int_0^1 -P(t)e^{-\alpha t} u'(t) + \alpha P(t) e^{-\alpha t} u(t) \mathrm{d}t \\
    = -&\left\{\left[P(t) e^{-\alpha t} u(t) \right]_0^1 - \int_0^1 \left[P'(t) - \alpha P(t)\right]e^{-\alpha t}u(t) \mathrm{d}t\right\} + \int_0^1 \alpha P(t) e^{-\alpha t} u(t) \mathrm{d}t \\
    &= -\left[P(t) e^{-\alpha t} u(t) \right]_0^1 + \int_0^1 P'(t)e^{-\alpha t} u(t) \mathrm{d}t \\
    &= P(0)u(0) - P(1)e^{-\alpha}u(1)  + \int_0^1 P'(t)e^{-\alpha t} \mathrm{d}t \\
    &= \Big< \left( P(0)-P(1)e^{-\alpha}\right) \Sha + P'(t)e^{-\alpha t}, u \Big>
\end{align*}
which allows to conclude by periodicity of $u$, as $u(0) = u(1)$.
\end{proof}
Note that this formula is still valid when $P$ is of degree $0$ (a constant polynomial) with $P' = 0$. We shall now demonstrate the previous theorem.

\begin{proof}[Proof of Theorem~\ref{th:exp_spline}]
Let $\mathcal{D} = \left( \mathrm{D} + \alpha \mathrm{Id}\right)^N $ for a given $N\in\mathbb{N}$, $N\geq 2$. Let us verify that $t \mapsto P_{N-1}(t)e^{-\alpha t}$ is solution of the Green's function equation $\mathcal{D}\psi_\mathcal{D} = \Sha$. By iteratively applying $N$ times the formula \eqref{eq:exp_to_green}, one obtains
\begin{equation}
    \left( \mathrm{D} + \alpha \mathrm{Id}\right)^N \left[P_{N-1}(t)e^{-\alpha t}\right] = \sum_{k=0}^{N-1}{ \left(P_{N-1}^{(k)}(0) - P_{N-1}^{(k)}(1)e^{-\alpha}\right) \left( \mathrm{D} + \alpha \mathrm{Id}\right)^{N-1-k}\left\{\Sha\right\} }.
    \label{eq:step_1}
\end{equation}
Note that, as a polynomial, $P_{N-1}$ satisfies:
\begin{align*}
    P_{N-1}^{(k)}(0) \quad &= \quad (k!)\ a_k\quad \\
    P_{N-1}^{(k)}(1) \quad &= \sum_{i=0}^{N-1-k}{ \frac{(k+i)!}{i!} a_{k+i} }.
\end{align*}
Then, for any $k\in\left\{0, 1, \dots, N-2\right\}$,
\begin{align*}
    P_{N-1}^{(k)}(0) - P_{N-1}^{(k)}(1)e^{-\alpha} &= (k!) a_k - e^{-\alpha} \sum_{i=0}^{N-1-k}{ \frac{(k+i)!}{i!} a_{k+i} } \\
    &= b_k - e^{-\alpha} \sum_{i=0}^{N-1-k}{ \frac{1}{i!} b_{k+i} } \\
    &= (1 - e^{-\alpha})b_k - e^{-\alpha}\sum_{i=1}^{N-1-k}{ \frac{1}{i!} b_{k+i} } \\
    &=(1 - e^{-\alpha}) A \sum_{i=1}^{N-1-k}{ \frac{1}{i!} b_{k+i} } - e^{-\alpha} \sum_{i=1}^{N-1-k}{ \frac{1}{i!} b_{k+i} } \\
    &=0
\end{align*}
and for $k=N-1$,
\begin{equation*}
    P_{N-1}^{(N-1)}(0) - P_{N-1}^{(N-1)}(1)e^{-\alpha} = (N-1)!\ a_{N-1}\ ( 1-e^{-\alpha}) = 1 
\end{equation*}
Equation \eqref{eq:step_1} turns into
\begin{equation*}
    \mathcal{D} \left(P_{N-1}(t)e^{-\alpha t}\right) = \Sha.
\end{equation*}
As $\mathcal{D}$ is invertible, the solution is unique and $\psi_\mathcal{D}(t) = P_{N-1}(t)e^{-\alpha t}$ for any $t\in\mathbb{T}$.
\end{proof}

\begin{remark}
    In a matrix setting, the coefficients $b_k$ are obtained as the vector 
\begin{equation*}
    \mathbf{b} = \begin{pmatrix}
    b_{N-1} & b_{N-2} & \dots & b_1 & b_0
    \end{pmatrix}^T
\end{equation*}
solution of the eigenvector problem
\begin{equation}
    \mathrm{M}\mathbf{b} = \frac{1}{A} \mathbf{b}
\end{equation}
with $b_{N-1} = 1/(1-e^{-\alpha})$ and
\begin{equation}
    \mathrm{M} = \left[ \begin{array}{*{9}{@{}C{\mycolwd}@{}}}
    A^{-1} & 0 & & & & \dots & 0 \\
    1 & 0 & & & & & \\
    1/2! & 1 & 0 & & & & \vdots \\
     & 1/2! & \ddots & \ddots & & &  \\
    & & & & \ddots & & \\
    \frac{1}{(N-2)!} & \dots & & \ddots & 1 & 0 & 0 \\
    \frac{1}{(N-1)!} & \frac{1}{(N-2)!} & \dots & & 1/2! & 1 & 0 
    \end{array}\right]
\end{equation}
$\mathrm{M}$ is lower triangular so that the computation of $\mathbf{b}$ is actually straightforward, the coefficients being determined as soon as $b_{N-1}$ is given.
\end{remark}

\begin{remark}
    One might notice that the polynomials $P_N$ are the consecutive derivative one of the other, according to the equation $(P_N)' = P_{N-1}$. Indeed, from equation \eqref{eq:exp_to_green}, one obtains
    \begin{align*}
        \left( \mathrm{D} + \alpha \mathrm{Id}\right)^{N+1}\left(P_N(t)e^{-\alpha t}\right) = & \left( \mathrm{D} + \alpha \mathrm{Id}\right)^{N} \left\{\left[P_N(0)- P_N(1)e^{-\alpha}\right] \Sha\right\} \\
        &+ \left( \mathrm{D} + \alpha \mathrm{Id}\right)^{N}\left\{(P_N)'(t)e^{-\alpha t}\right\} \\
        = & \left( \mathrm{D} + \alpha \mathrm{Id}\right)^{N}\left\{(P_N)'(t)e^{-\alpha t}\right\} \\
        = & \ \Sha
    \end{align*}
    which allows to conclude.
\end{remark}

Such a Green's function behaves pretty well over the period. It has continuous derivative up to the order $(N-2)$ (see \cite[Section~5.1]{fageot2020tvbased}), so it is always continuous for $N \geq 2$, and it admits a unique maximum. Now if we pay attention to the expression of the adjoint operator \eqref{eq:adjoint_measurement}, we notice that any dual certificate $\eta$ would be expressed as a linear combination of such shifted Green's functions $\psi_\mathcal{D}(\theta_\ell - \cdot)$. It is actually not a spline as the Green's functions are evaluated in the "reverse" sense (it would need $(\cdot -\theta_\ell)$ to be a spline, and the Green's functions are not symmetric with respect to their knot with an exponential operator). The question that arises from now on is: can such a finite linear combination of (reversed) shifted Green's function admit an infinite number of saturation points (given that its infinite norm is equal to $1$)? For the sake of generality, we will consider the case of regular splines instead of certificates, but the conclusions would directly transfer as the reversing operation does not affect this property.

Let us analyze how such a spline behaves. For $\mathcal{D} = \left( \mathrm{D} + \alpha \mathrm{Id}\right)^N$ with $N \geq 2$, consider the function $f\in\mathcal{M}_\mathcal{D}\left(\mathbb{T}\right)$ defined as
\begin{equation*}
    \forall t \in \mathbb{T}, \quad f(t) = \sum_{\ell=1}^L {\beta_\ell \psi_\mathcal{D}(t-\theta_\ell)}
\end{equation*}
such that $\norm{f}_\infty = 1$. As $\psi_\mathcal{D}$ is continuous, $f$ is also continuous.

First let us rewrite how the shifting operation modifies the expression of the Green's function. Consider $\theta \in \mathbb{T}$, then
\begin{equation}
    \forall t \in \mathbb{T}, \quad t-\theta - \floor{t-\theta} = \left\{\begin{tabular}{l c}
        $t - \theta$ &  $\mathrm{if} \quad \theta \leq t < 1$\\
        $t + 1 - \theta$ & $\mathrm{if}\quad 0 \leq t < \theta$ 
    \end{tabular}
    \right.
\end{equation}
which leads to
\begin{align*}
    \forall t \in \mathbb{T}, \quad \psi_\mathcal{D}(t-\theta) &= P(t-\theta - \floor{t-\theta}) e^{-\alpha (t-\theta - \floor{t-\theta})} \\[0.2cm]
    &= \left\{\begin{tabular}{l c}
        $P(t-\theta)\ e^{\alpha \theta}\ e^{-\alpha t}$ &  $\mathrm{if} \quad \theta \leq t < 1$\\[0.1cm]
        $P(t + 1 - \theta)\ e^{\alpha (\theta-1)}\  e^{-\alpha t}$ & $\mathrm{if}\quad 0 \leq t < \theta$ 
    \end{tabular} \right. \\[0.2cm]
    &= \left\{\begin{tabular}{l c}
        $P_1(t)\ e^{-\alpha t}$ &  $\mathrm{if} \quad \theta \leq t < 1$\\[0.1cm]
        $P_2(t)\ e^{-\alpha t}$ & $\mathrm{if}\quad 0 \leq t < \theta$ 
    \end{tabular} \right. 
\end{align*}
with $P_1(\theta) = P_2(\theta) = P(0)$ by continuity of $\psi_\mathcal{D}$. Thus, any shifted Green's function can be expressed as a piecewise polynomial over the period multiplied by a decreasing exponential. From this statement, we can deduce the following result:
\begin{proposition}
Let $N\in\mathbb{N}\setminus\{0, 1\}$ and $\mathcal{D} = \left( \mathrm{D} + \alpha \mathrm{Id}\right)^N$. Any continuous $\mathcal{D}$-spline $f = \sum_{\ell=1}^L {\beta_\ell \psi_\mathcal{D}(\cdot - \theta_\ell)}$ with infinite norm $\norm{f}_\infty = 1$ admits a finite number P of saturation points $t\in\mathbb{T}$ such that $|f(t)| = 1$, and $P\leq LN$. 
\end{proposition}
\begin{proof}
Consider the $L$-partitioning of $\mathbb{T}$ in the sequence of the $L$ intervals $\left(\left] \theta_\ell, \theta_{\ell+1} \right[\right)_{\ell=1, \dots, L}$ with the convention $\theta_{L+1} = \theta_1$. On each subinterbval, $f$ can actually be written as $f_\ell(t) = P_\ell (t)e^{-\alpha t}$ with $P_\ell$ being a polynomial of degree $(N-1)$. That way, $f$ admits at most $(N-1)$ saturation points on each interval. Indeed, the norm constraint enforces that any saturation point is a local extremum of $f_\ell$ so that $f'_\ell$ is $0$ at these points. $f'_\ell$ is written as a product between a polynomial of degree $(N-1)$ and $ e^{-\alpha x}$ and thus cannot have more than $(N-1)$ zeros. By summing over all the subintervals and adding the $L$ innovation positions, there are at most $L + L (N-1) = LN$ saturation points over $\mathbb{T}$.
\end{proof}
Note that the proof of the previous proposition exhibits an upper bound for the number of saturation points, but it might not be sharp at all, and it is probably possible to determine a much lower bound. \\

Numerous questions remain still not answered. So far, we have not succeeded in proving Proposition \ref{prop:uniqueness} applied to the case of exponential Green's function. However, the research is ongoing and we believe that it is possible to determine a set conditions according to which an operator $\mathcal{D}$ would ensure uniqueness in the solution set. Moreover, the same analysis of the expression of periodic splines can probably be pursued with the other common differential operators (like the Sobolev operator).

The other way around, one may also ask if it is possible, starting from a function that would satisfy \cref{prop:discrete} and \cref{prop:uniqueness}, to build an associated spline-admissible operator. In that case, such an operator would ensure that the solution is unique and is always a spline.